\documentclass[
preprint, 3p, 
number, 
sort&compress,
]{elsarticle}
\pdfoutput=1

\usepackage[utf8]{luainputenc}
\usepackage[english]{babel}
\usepackage{csquotes}

\usepackage[plainpages=false,pdfpagelabels,hidelinks,unicode]{hyperref}

\usepackage{amsmath}
\allowdisplaybreaks
\usepackage{amssymb}
\usepackage{commath}
\usepackage{mathtools}
\usepackage{bbm}

\usepackage{siunitx}
\usepackage{adjustbox}

\usepackage{amsthm}
\theoremstyle{plain}
  \newtheorem{theorem}{Theorem}

\theoremstyle{definition}
  
  \newtheorem{remark}[theorem]{Remark}

\usepackage{color}
\usepackage{graphicx}
\usepackage[small]{caption}
\usepackage{subcaption}

\ifx\useTikzForPlotting\undefined
\else
  \usepackage{pgfplots}
  \pgfplotsset{compat=1.11}
  \usetikzlibrary{external}
\fi

\usepackage{booktabs}
\usepackage{rotating}
\usepackage{multirow}

\usepackage{multicol}
\usepackage{enumitem}

\usepackage{calc}
\usepackage{xparse}

\DeclareMathOperator*{\argmin}{arg\,min}
\renewcommand{\vec}[1]{\underline{#1}}
\NewDocumentCommand{\mat}{mo}{%
  \IfValueTF{#2}{%
    \underline{\underline{#1}}{#2}
  }{%
    \underline{\underline{#1}}\,
  }%
}

\newcommand{\scp}[2]{\left\langle{#1,\, #2}\right\rangle}

\renewcommand{\d}{\mathrm{d}} 
\newcommand{\intd}{\ \mathrm{d}}

\newcommand{\ubar}{\overline{u}}
\newcommand{\fnum}{f^{\mathrm{num}}}

\renewcommand{\epsilon}{\varepsilon}
\renewcommand{\phi}{\varphi}
\renewcommand{\rho}{\varrho}
\newcommand{\N}{\mathbb{N}}
\newcommand{\R}{\mathbb{R}}

\newsavebox{\DelimiterBox}
\newlength{\DelimiterHeight}
\newlength{\DelimiterDepth}
\newsavebox{\ArgumentBox}
\newlength{\ArgumentHeight}
\newlength{\ArgumentDepth}
\newlength{\ResizedDelimiterHeight}
\newlength{\ResizedDelimiterDepth}

\usepackage{lipsum}
\makeatletter
\def\ps@pprintTitle{%
 \let\@oddhead\@empty
 \let\@evenhead\@empty
 \def\@oddfoot{}%
 \let\@evenfoot\@oddfoot}
\makeatother

\begin{document}

\begin{frontmatter}

\title{Stable discretisations of high-order discontinuous Galerkin methods on equidistant and scattered points}

\author[label1,label2]{Jan Glaubitz\corref{cor1}}
\ead{j.glaubitz@tu-bs.de}

\author[label3]{Philipp \"Offner}

\cortext[cor1]{Corresponding author}

\address[label1]{Max Planck Institute for Mathematics, Bonn, Germany}
\address[label2]{Institute for Computational Mathematics, TU Braunschweig, Braunschweig, Germany}
\address[label3]{Institut f\"ur Mathematik, Universit\"at Z\"urich, Z\"urich, Switzerland}

\begin{abstract}
  In this work, we propose and investigate stable high-order collocation-type discretisations of the discontinuous 
Galerkin method on equidistant and scattered collocation points. 
We do so by incorporating the concept of discrete least squares into the discontinuous Galerkin framework. 
Discrete least squares approximations allow us to construct stable and high-order accurate approximations on arbitrary 
collocation points, while discrete least squares quadrature rules allow us their stable and exact numerical 
integration. 
Both methods are computed efficiently by using bases of discrete orthogonal polynomials. 
Thus, the proposed discretisation generalises known classes of discretisations of the discontinuous Galerkin method, 
such as the discontinuous Galerkin collocation spectral element method. 
We are able to prove conservation and linear $L^2$-stability of the proposed discretisations. 
Finally, numerical tests investigate their accuracy and demonstrate their extension to nonlinear conservation laws, 
systems, longtime simulations, and a variable coefficient problem in two space dimensions.
\end{abstract}


\end{frontmatter}

\section{Introduction}
\label{sec:introduction}

In the last decades, great efforts have been made to develop accurate and stable
numerical methods for time-dependent partial differential equations (PDEs), 
especially for hyperbolic conservation laws 
\begin{equation}\label{eq:cons-law}
  u_t + f(u)_x = 0. 
\end{equation}
The entropy solution of \eqref{eq:cons-law} satisfies the additional entropy condition 
\begin{equation}\label{eq:entropy-ineq}
  U(u)_t + F(u)_x \leq 0 
\end{equation} 
in the sense of distributions. 
$U$ is a convex entropy function and $F$ is a corresponding entropy flux satisfying ${U' f' = F'}$. 
A strict inequality in \eqref{eq:entropy-ineq} reflects the existence of physically reasonable 
shock waves.

Traditionally, low-order numerical schemes, for instance classical finite volume (FV) methods, 
have been used to solve hyperbolic conservation laws, particularly in industrial applications. 
But since they become quite costly for high accuracy or long time simulations, 
there is a rising demand of high-order (third and above) methods. 
These have the potential of providing accurate solutions at reasonable costs. 
In this work, we consider the particularly popular class of discontinuous Galerkin (DG) finite 
element methods for hyperbolic conservation laws \eqref{eq:cons-law}. 
These methods were first introduced 1973 by Reed and Hill \cite{reed1973triangular} to solve the 
hyperbolic neutron transport equation in a nuclear reactor and were put on mathematically solid 
ground by Cockburn, Shu, and co-authors in a series of papers 
\cite{cockburn1991runge,cockburn1989tvb,cockburn1989tvb2,cockburn1990runge,cockburn1998runge} 
around 1990. 
In \cite{cockburn1989tvb} and \cite{cockburn1990runge} it was proven that the resulting class of 
DG methods is (formally) high-order accurate in smooth regions, total-variation 
bounded in one space dimension, and maximum-norm bounded in any number of space dimensions. 
Further, Jiang and Shu proved in \cite{jiang1994cell} a cell entropy inequality for the square 
entropy 
\begin{equation}
  U(u) = \frac{u^2}{2}
\end{equation} 
of linear as well as nonlinear scalar conservation laws for the DG method. 
It should be noted that this result does not need any nonlinear limiting as introduced in 
\cite{cockburn1989tvb} and \cite{cockburn1990runge}. 
The cell entropy inequality makes the DG method consistent with the entropy 
condition \eqref{eq:entropy-ineq} and implies $L^2$-stability of the scheme. 
Yet, Jiang and Shu's proof relies on exact evaluation of integrals and thus only applies to the 
\emph{analytical} DG method, where all arising integrals are evaluated exactly. 
Unfortunately, exact evaluation of integrals is often computationally impractical or, depending on 
the nonlinearity of the flux function $f$, even impossible and is hence usually replaced by 
numerical quadrature rules.  
Then, exactness (at least up to machine precision) can be guaranteed by a sufficiently great number 
of integration points. 
This was, for instance, investigated by Kirby and Karniadakis in \cite{kirby2003aliasing}. 
Another alternative  are so-called DG collocation spectral element methods (DGSEMs) \cite{hesthaven2007nodal}.
In these methods, the solution $u$ as well as the flux $f$ are approximated by interpolation 
polynomials in each element and the corresponding interpolation points are further matched with the 
integration points, which results in highly efficient operators \cite{kopriva2009implementing}. 
A typical problem of such discretisation of the DG method is their 
instability, especially for nonlinear conservation laws, due to the reduced accuracy of the 
integrals. 
There are several possible stabilisation methods in the literature, such as 
minmod-type limiting 
\cite{cockburn1989tvb,cockburn1990runge},
artificial viscosity methods 
\cite{persson2006sub,klockner2011viscous,ranocha2018stability,glaubitz2019smooth,offner2019stability}, 
modal filtering 
\cite{vandeven1991family,hesthaven2008filtering,glaubitz2018application,ranocha2018stability}, 
finite volume sub-cells 
\cite{huerta2012simple,dumbser2014posteriori,sonntag2014shock,meister2016positivity}, 
and many more. 
Yet, in \cite{gassner2013skew,kopriva2014energy,gassner2016well} 
Gassner, Kopriva, and co-authors have been able to construct 
a DGSEM which is $L^2$-stable for 
certain linear (variable coefficient) as well as nonlinear conservation laws by utilising 
skew-symmetric formulations of the conservation law \eqref{eq:cons-law} and Summation-by-Parts 
(SBP) operators, which were first used and investigated in finite difference (FD) methods  
\cite{kreiss1974finite,strand1994summation,olsson1995summation,olsson1995summation2,nordstrom1999boundary}.
It should be stressed that the theoretical stability (in the sense of a provable $L^2$-norm inequality) as 
well as the numerical stability (in the sense 
of stable interpolation polynomials and quadrature rules) of the DGSEM heavily rely on the usage 
of Gauss--Lobatto points and quadrature weights, which include the boundary nodes and are more 
dense there, see \cite{gassner2013skew,kopriva2014energy}. 
In \cite{ranocha2016summation} these results have been extended to 
Gauss--Legendre points and quadrature weights, which do not include the boundary nodes but are still more 
dense there. 
A more general approach, ensuring entropy stability by using decoupled SBP operators was recently proposed in 
\cite{chan2019efficient} by Chen and coauthors.
Their approach is quite general and will be addressed in greater detail in \S \ref{sub:entropy}.

In this work, we describe a first step towards $L^2$-stable DGSEM on equidistant and even scattered collocation points.
Yet, when using equidistant points for the (interpolation) polynomials approximating 
$u$ and $f$ in the DGSEM, we run into the Runge phenomenon \cite{runge1901empirische}. 
Adapting an idea of Gelb et al. \cite{gelb2008discrete} from spectral collocation methods, we 
tackle this problem by making a somewhat maverick generalisation. 
Instead of usual polynomial approximations by interpolation on $K+1$ points, we build the polynomial 
approximations by the method of discrete least squares, where the data at a greater number of 
$N > K+1$ points is used. 
The method thus utilises more information from the underlying function. 
By going over to a higher number of nodal values than technically needed, polynomial discrete least squares 
approximations are known to provide high accuracy while also successfully suppressing Runge 
oscillations, even on equidistant points. 
Care also has to be taken when performing numerical integration on equidistant or even 
scattered collocation points. 
To prove conservation as well as $L^2$-stability, quadrature rules which exactly 
evaluate polynomials of degree $2K$ are needed. 
While interpolatory quadrature rules of arbitrary high order of exactness can be constructed even 
on equidistant (and scattered) points, they quickly become unstable, 
see \cite{huybrechs2009stable}. 
Thus, for an exact as well as stable evaluation of integrals, we propose to use 
so-called \emph{least squares quadrature rules}. 
These have been introduced by Wilson in \cite{wilson1970necessary,wilson1970discrete} and since then only revisited by 
Huybrechs in \cite{huybrechs2009stable}. 
To the best of our knowledge, in particular, least squares quadrature rules have never been investigated in a 
PDE solver before. 
Utilising discrete least squares approximations and stable high-order quadrature rules, 
we are able to prove conservation and linear $L^2$-stability of the resulting discretisation 
of the DG method on equidistant and scattered points.

The rest of this work is organised as follows. 
In \S \ref{sec:DG}, we revisit DG methods (without additional stabilisation) 
and their usual collocation-type discretisation on Gauss--Lobatto points. 
\S \ref{sec:DLS} introduces the concepts of discrete least squares approximations and least squares quadrature 
rules, which will provide stable and high-order numerical integration on equidistant and even scattered points. 
Their efficient computation is based on bases of discrete orthogonal polynomials. 
Building up on these concepts, we propose our stable high-order discretisation of the DG method on 
equidistant and scattered points in \S \ref{sec:proposed}. 
The resulting discretisation of the DG method, which is based on the concept of discrete least squares, is referred to 
as the \emph{DG discrete least squares (DGDLS) method} and generalises the usual discretisation of 
the DG method by the DGSEM. 
Conservation as well as linear $L^2$-stability of the DGDLS method are proven in \S \ref{sec:cons_stab}. 
In \S \ref{sec:num}, numerical tests demonstrate conservation and $L^2$-stability of the proposed DGDLS method for the 
linear advection equation and the nonlinear inviscid Burgers' equation. 
Further, we investigate accuracy of the method and address the extension to systems of conservation laws, longtime 
simulations, and a variable coefficient problem in two spatial dimensions.  
We close this work with a summary and outlook in \S \ref{sec:summary}.

\section{The discontinuous Galerkin method and its discretisation}
\label{sec:DG}

Decoupling space and time by the method of lines \cite{leveque2002finite}, DG 
methods are designed as semidiscretisations of hyperbolic conservation laws, 
\begin{equation}\label{eq:cons-law2}
    u_t + f(u)_x = 0, 
\end{equation}
on a computational domain $\Omega$. 
In the following, we will only consider conservation laws in one dimension, i.\,e.\ $\Omega \subset 
\R$. 
For sake of simplicity, we further assume \eqref{eq:cons-law2} to be a scalar conservation law for the moment. 
The extension of the later proposed discretisation will be addressed in \S \ref{sub:system}. 
We refer to the function $u=u(t,x)$ as the \emph{conservation variable} and to $f$ as 
the \emph{flux function}.
The resulting system of ordinary differential equations (ODEs)
\begin{equation}\label{eq:ODE}
  \frac{\d}{\d t} u = L(u),  
\end{equation} 
where $L(u)$ is a discretisation of the spatial operator, can be solved by any time integration 
method. 
In this work, we use the explicit strong stability preserving (SSP) Runge--Kutta (RK)
method\footnote{SSP-RK methods are a generalisation of TVD-RK methods, where instead of the total 
variation of the solution $u$ every convex functional of the solution $u$ is ensured to decrease if 
this holds for a step of the simple explicit Euler method.} 
of third order using three stages (SSPRK(3,3)) given in \cite{gottlieb1998total} by Gottlieb and 
Shu: 
Let $u^n$ be the solution at time $t^n$, then the solution $u^{n+1}$ at time $t^{n+1}$ is obtained 
by 
\begin{equation}
\begin{aligned}
  u^{(1)} & = u^n + \Delta t L\left( u^n \right), \\ 
  u^{(2)} & = \frac{3}{4} u^n + \frac{1}{4} u^{(1)} + \frac{1}{4} \Delta t L\left( u^{(1)} \right), 
\\ 
  u^{n+1} & = \frac{1}{3} u^n + \frac{2}{3} u^{(2)} + \frac{2}{3} \Delta t L\left( u^{(2)} \right). 
\end{aligned}
\end{equation}
For this method, in particular, $L^2$-stability will also hold in time if it is ensured for 
the simple explicit Euler method.
For more details about (TVD/SSP) RK methods see the extensive literature 
\cite{gottlieb1998total,levy1998semidiscrete,gottlieb2001strong,ketcheson2008highly,gottlieb2011strong}. 
The later numerical tests in \S \ref{sec:num} have all been performed using the timestep size 
\begin{equation}
  \Delta t = \frac{C}{I (K+1) \lambda}
\end{equation}
with $C = 0.1$ and where $\lambda = \max_u |f'(u)|$ is the fastest propagation speed. 
In the following, we discuss the discretisation $L(u)$ of the spatial operator given by the 
DG method.

\subsection{The analytical discontinuous Galerkin method} 
\label{sub:anal-DG}

DG methods are obtained in the following way. 
First, the computational domain $\Omega$ is subdivided into smaller elements $\Omega_i$ for 
$i=1,\dots,I$ with $\Omega = \bigcup_{i=1}^I \Omega_i$. 
In one dimension, the elements are given by simple subintervals 
\begin{equation}
  \Omega_i = \left( x_{i-\frac{1}{2}}, x_{i+\frac{1}{2}} \right), 
\end{equation} 
which are typically mapped to a single reference element $\Omega_{\mathrm{ref}} = (-1,1)$ where 
all computations are performed. 
In two or more space dimensions, the elements might be given by triangles, tetrahedrons or other 
simple geometric objects.
The conservation law \eqref{eq:cons-law2} is solved in a weak form for every 
element then. 
On an element $\Omega_i$, equation \eqref{eq:cons-law2} is multiplied by a test function $v$ and 
integrated in space. 
Integration by parts yields the weak form 
\begin{equation}
  \int_{\Omega_i} u_t v \intd x - \int_{\Omega_i} f(u) v_x \intd x 
    + f\left( u(x_{i+\frac{1}{2}}) \right) v(x_{i+\frac{1}{2}}) 
    - f\left( u(x_{i-\frac{1}{2}}) \right) v(x_{i-\frac{1}{2}}) 
    = 0.
\end{equation}
Now assume that both the solution $u$ as well as the test function $v$ come from a finite 
dimensional approximation space $V_h$, which is usually chosen to be the space of piecewise 
polynomials of degree at most $K$, i.\,e.\ 
\begin{equation}
  V_h = \left\{ v:\Omega \to \R \ \Big| \ v^i := v|_{\Omega_i} \in \mathbb{P}_K(\Omega_i) \right\}.
\end{equation}
It should be stressed that for this choice $u,v \in V_h$ might be discontinuous at element 
interfaces.
Thus, the boundary terms $f\left(u(x_{i \pm \frac{1}{2}})\right)$ and $v(x_{i \pm \frac{1}{2}})$ are not well defined. 
We distinguish between values from inside $\Omega_i$, e.\,g.\ 
\begin{equation}
  u^+_{i - \frac{1}{2}} = u^i( x_{i - \frac{1}{2}} ), \quad
  u^-_{i + \frac{1}{2}} = u^i( x_{i + \frac{1}{2}} ), \quad
  v^+_{i - \frac{1}{2}} = v^i( x_{i - \frac{1}{2}} ), \quad
  v^-_{i + \frac{1}{2}} = v^i( x_{i + \frac{1}{2}} )
\end{equation}
and values from neighbouring elements $\Omega_{i-i},\Omega_{i+i}$, e.\,g.\ 
\begin{equation}
  u^-_{i - \frac{1}{2}} = u^{i-1}( x_{i - \frac{1}{2}} ), \quad
  u^+_{i + \frac{1}{2}} = u^{i+1}( x_{i + \frac{1}{2}} ), \quad
  v^-_{i - \frac{1}{2}} = v^{i-1}( x_{i - \frac{1}{2}} ), \quad
  v^+_{i + \frac{1}{2}} = v^{i+1}( x_{i + \frac{1}{2}} ).
\end{equation}
From conservation and stability (upwinding) considerations, we take a single valued 
\emph{numerical flux} 
\begin{equation}
  \fnum_{i + \frac{1}{2}} = \fnum\left( u^-_{i + \frac{1}{2}}, u^+_{i + \frac{1}{2}} \right)
\end{equation}
to replace $f\left( u(x_{i + \frac{1}{2}}) \right)$. 
Further, the numerical flux is consistent ($\fnum(u,u) = f(u)$), Lipschitz continuous, and monotone 
($\fnum$ is nondecreasing in the first argument and nonincreasing in the second argument).
Examples of commonly used numerical fluxes can be found in \cite{cockburn1989tvb} and 
\cite{toro2013riemann}. 
Hence, the DG method is: 
Find $u \in V_h$ such that 
\begin{equation}\label{eq:anal-DG}
  \int_{\Omega_i} u_t v \intd x - \int_{\Omega_i} f(u) v_x \intd x 
    + \fnum_{i + \frac{1}{2}} v^-_{i + \frac{1}{2}} 
    - \fnum_{i - \frac{1}{2}} v^+_{i - \frac{1}{2}} 
    = 0
\end{equation} 
for all $v \in V_h$ and $i=1,\dots,I$. 
Note that in \eqref{eq:anal-DG} all integrals are assumed to be evaluated exactly. 
We thus refer to \eqref{eq:anal-DG} as the \emph{analytical DG method}.
Finally, the analytical DG method \eqref{eq:anal-DG} can be rewritten as a system of ODEs 
\begin{equation}
  \frac{\d}{\d t} u = L(u)   
\end{equation} 
and the approximation $u \in V_h$ is evolved over time by some explicit TVD-RK method, as 
discussed before. 
Also see \cite{cockburn1989tvb} and \cite{hesthaven2007nodal}.

\subsection{The discontinuous Galerkin collocation spectral element method} 
\label{sub:DGSEM}

This subsection revisits the DGSEM in which $u$ and $f$ are both approximated by interpolation 
polynomials in each element and the corresponding interpolation points are matched with the 
integration points. 
We assume $u,v, f(u) \in V_h$, i.\,e.\ all functions are approximated by piecewise polynomials of 
degree up to $K$. 
Further, it is convenient to transform equation \eqref{eq:cons-law2} for each element 
$\Omega_i = \left( x_{i-\frac{1}{2}}, x_{i+\frac{1}{2}} \right)$ to the reference element 
$\Omega_{\mathrm{ref}} = (-1,1)$ via the linear map 
\begin{equation}\label{eq:transformation}
  x_i(\xi) = \overline{x}_i + \frac{\Delta x_i}{2} \xi,
\end{equation}
where $\overline{x}_i$ is the centre of the element $\Omega_i$ and $\Delta x_i$ is its length.
This results in the transformed equation 
\begin{equation}\label{eq:trans_eq}
    \frac{\Delta x_i}{2} u_t + f(u)_\xi = 0 
\end{equation} 
on $\Omega_{\mathrm{ref}} = (-1,1)$.
By defining a set of $K+1$ interpolation points $-1 \leq \xi_0 < \xi_1 < \dots < \xi_K \leq 1$ in 
the reference element, $u^i = u\big|_{\Omega_i}$ and 
$f^i = f(u)\big|_{\Omega_i}$ are computed by polynomial interpolation as 
\begin{equation}
  u^i(t,\xi) = \sum_{k=0}^K u^i_k(t) \ell_k(\xi) 
  \quad \text{and} \quad 
  f^i = \sum_{k=0}^K f^i_k \ell_k(\xi)
\end{equation} 
with time depended nodal values $u^i_k(t) = u^i(t,\xi_k)$, $f^i_k = f(u^i_k)$ at the 
interpolation points and corresponding Lagrange basis functions $\ell_k$, which are defined by 
\begin{equation}
    \ell_k(\xi) = \prod_{j=0, j \neq k}^K \frac{ \xi - \xi_j }{ \xi_k - \xi_j }
\end{equation}
and satisfy the cardinal property $\ell_k(x_j) = \delta_{kj}$.
For sake of simplicity, we will just focus on the Gauss--Lobatto points for the interpolation as well 
as quadrature points $\{ \xi_j \}_{j=0}^K$ and denote the associated quadrature weights by $\{ 
\omega_j \}_{j=0}^K$. 
Analytic integration is replaced by the resulting quadrature rule 
\begin{equation}
    \int_{-1}^1 g(\xi) \intd \xi \approx \sum_{j=0}^K \omega_j g(\xi_j)
\end{equation}
then.
Besides the Gauss--Lobatto points, the Gauss--Legendre points are another typical choice for 
the collocation approach. 
Inserting these approximations into the DG formulation \eqref{eq:anal-DG} and choosing $v^i = 
\ell_k$, we get the \emph{DGSEM}: 
Solve 
\begin{equation}\label{eq:DGSEM}
  \omega_k \frac{\Delta x_i}{2} \frac{\d}{\d t} u^i_k 
    - \sum_{j=0}^K \omega_j f_j^i \ell_k'(\xi_j) 
    + \fnum_{i+\frac{1}{2}} \ell_k(1) - \fnum_{i-\frac{1}{2}} \ell_k(-1) 
  = 0
\end{equation} 
for all $k=0,\dots,K$ and $i=1,\dots,I$.
Note that \eqref{eq:DGSEM} is directly solved for the nodal degrees of freedom 
$\{u^i_k\}_{k=0,i=1}^{K,I}$ of the solution $u \in V_h$. 
 
\section{Discrete least squares approximations and quadrature rules}
\label{sec:DLS}

In this section, we revisit the concept of discrete least squares (DLS), which is closely related to 
(continuous) weighted inner products. 
The later proposed discretisation of the DG method on equidistant and scattered points is based on piecewise 
polynomial approximations of $u$ and $f$. 
In every element, these approximations are computed by DLS approximations, which generalise polynomial interpolation. 
Further, exact integrals (in the analytical DGM) and usual Gauss--Lobatto quadrature rules 
(in the DGSEM), respectively, is replaced by stable high-order least squares quadrature rules (LS-QRs). 
The efficient implementation of both methods relies on bases of discrete orthogonal polynomials (DOPs).  

\subsection{Bases of discrete orthogonal polynomials}
\label{sub:DOP}

Let $N \in \N$, $\{\xi_n\}_{n=0}^N$ be a set of $N+1$ distinct points in $[-1,1]$, and 
$\vec{\omega} = [ \omega_0, \dots, \omega_N ]^T \in \R^{N+1}$ be a vector of corresponding positive 
weights. 
Then, in analogy to usual (continuous) inner products, 
\begin{equation}\label{eq:disc-scp}
  \scp{u}{v}_{\vec{\omega}} 
    := \sum_{n=0}^N \omega_n u(\xi_n) v(\xi_n)
\end{equation}
is called a \emph{discrete inner product}. 
An accompanying discrete norm is induced by 
\begin{equation}
  \norm{u}_{\vec{\omega}}^2 
    = \scp{u}{u}_{\vec{\omega}} 
    = \sum_{n=0}^N \omega_n u(\xi_n)^2.
\end{equation}
Now considering a basis $\{ \phi_k \}_{k=0}^K$ of $\mathbb{P}_K([-1,1])$ with $K \leq N$, we call 
the basis elements \emph{DOPs} if 
\begin{equation}
  \scp{\phi_k}{\phi_l}_{\vec{\omega}} = \delta_{kl} 
\end{equation}
holds for $k,l=0,\dots,K$, i.\,e.\ if they are orthogonal w.\,r.\,t.\ to the discrete inner product 
\eqref{eq:disc-scp}. 
For more details, we refer to the book \cite{gautschi2004orthogonal} of Gautschi.
It should be stressed that for (nonclassical) DOPs often no explicit formula is 
known.\footnote{At least for classical DOPs on equidistant points some formulas can, for 
instance, be found in the book of Gautschi \cite[Chapter 1.5.2]{gautschi2004orthogonal}. 
These examples include the discrete Chebyshev polynomials, the Krawtchouk polynomials, the Charlier polynomials, the 
Meixner polynomials, and the Hahn polynomials.}
Thus, we utilise numerical algorithms to construct bases of 
DOPs, such as the Stieltjes procedure 
\cite{gautschi2004orthogonal} and the Gram-Schmidt process \cite{trefethen1997numerical}. 
In this work, we construct bases of DOPs 
by the numerical stable modified Gram--Schmidt process \cite{trefethen1997numerical} applied to an initial basis of 
Legendre polynomials.

\subsection{Discrete least square approximations}
\label{sub:DLS-approx}

Originally, DLS approximations are born from the wish to fit a linear mathematical model to given observations. 
In contrast to polynomial interpolation \cite{platte2011impossibility}, DLS approximations can be stable and highly 
accurate on equidistant and even scattered points. 
In this work, we use DLS approximations to fit polynomials $f_{K,N}$ of degree at most $K$ to a greater number of 
$N$ observations given by nodal values at the collocation points.

Let $f$ be a function on $[-1,1]$ which we only know at a set of $N+1$ distinct points $\{ \xi_n \}_{n=0}^N$ in 
$[-1,1]$. 
The problem is to find a polynomial $f_{K,N} \in \mathbb{P}_K$ such that 
${f_{K,N} - f}$
is minimised. 
Note that when we have a basis $\{ \varphi_k \}_{k=0}^K$ of $\mathbb{P}_K([-1,1])$, the polynomial $f_{K,N}$ can be 
represented as 
\begin{equation}
  f_{K,N}(\xi) = \sum_{k=0}^K \hat{f}_{k,N} \varphi_k(\xi), 
\end{equation}
where the $\hat{f}_{k,N}$ are the modal coefficients. 
Since we only know $f$ at a set of $N+1$ distinct points $\{ \xi_n \}_{n=0}^N$, it is convenient to compare $f_{K,N}$ 
and $f$ only at the points $\{ \xi_n \}_{n=0}^N$. 
Thus, let us denote the vector of the known nodal values of $f$ by $\vec{f} = \left( f(\xi_1), \dots, f(\xi_N) 
\right)^T$, the vector of the modal coefficients of $f_{K,N}$ by $\vec{\hat{f}} = \left( \hat{f}_{0,N}, \dots, 
\hat{f}_{K,N} \right)^T$, and the matrix which contains the values of the basis functions $\phi_k$ at the 
collocation points by
\begin{equation}
    \mat{A} = \left( \phi_{k}(\xi_{n}) \right)_{k,n=0}^{K,N}.  
\end{equation}
Then, the problem can be reformulated as follows:  
Find $\vec{\hat{f}} \in \R^{K+1}$ such that
${\vec{r} := \mat{A} \vec{\hat{f}} - \vec{f}}$
is minimised, where $\vec{r}$ is called the \emph{residual vector}.
In fact, there are many possible ways of defining a best approximation.
%
In this work, we minimise $f_{K,N}-f$ w.\,r.\,t.\ discrete inner products and thus choose $f_{K,N}$ as the 
\emph{DLS approximation} of $f$ characterised by 
\begin{equation}\label{eq:DLS-approx-char}
  \norm{ f_{K,N}-f }_{\vec{\omega}} = \min_{v \in \mathbb{P}_K} \norm{ v-f }_{\vec{\omega}}. 
\end{equation} 
Since the discrete norm $\norm{\cdot}_{\vec{\omega}}$ is induced by the discrete inner product \eqref{eq:disc-scp}, 
relation \eqref{eq:DLS-approx-char} is equivalent to 
\begin{equation}
  \scp{ f_{K,N} - f }{ v }_{\vec{\omega}} = 0 \quad \forall v \in \mathbb{P}_K
\end{equation}
and, when $f_{K,N}$ is represented w.\,r.\,t.\ a basis $\{ \varphi_k \}_{k=0}^K$ of $\mathbb{P}_K$, yields a system of 
linear equations 
\begin{align}\label{eq:DLS-LES}
    \sum_{k=0}^K \hat{f}_{k,N} \scp{\varphi_k}{\varphi_i}_{\vec{\omega}} 
      = \scp{f}{\varphi_i}_{\vec{\omega}}, 
    \quad i=0,\dots,K,
\end{align}
which can, for instance, be solved by Gaussian elimination. 
Yet, highest efficiency is obtained by choosing a basis of DOPs as discussed in \S \ref{sub:DOP}. 
Then, the sum in \eqref{eq:DLS-LES} reduces to a single entry and the coefficients are given by 
\begin{align}\label{eq:DLS_coeff}
    \hat{f}_{k,N} = \scp{f}{\varphi_k}_{\vec{\omega}}
\end{align}
for $k = 0,\dots,K$. 
When computed by using bases of DOPs, DLS approximations are sometimes referred to as 
\emph{DOP-LS approximations}, see \cite{gelb2008discrete}.

\subsection{Stable high-order quadrature rules on equidistant and scattered points}
\label{sub:QRs} 

\begin{figure}[!htb]
  \centering
    \includegraphics[width=0.5\textwidth]{%
        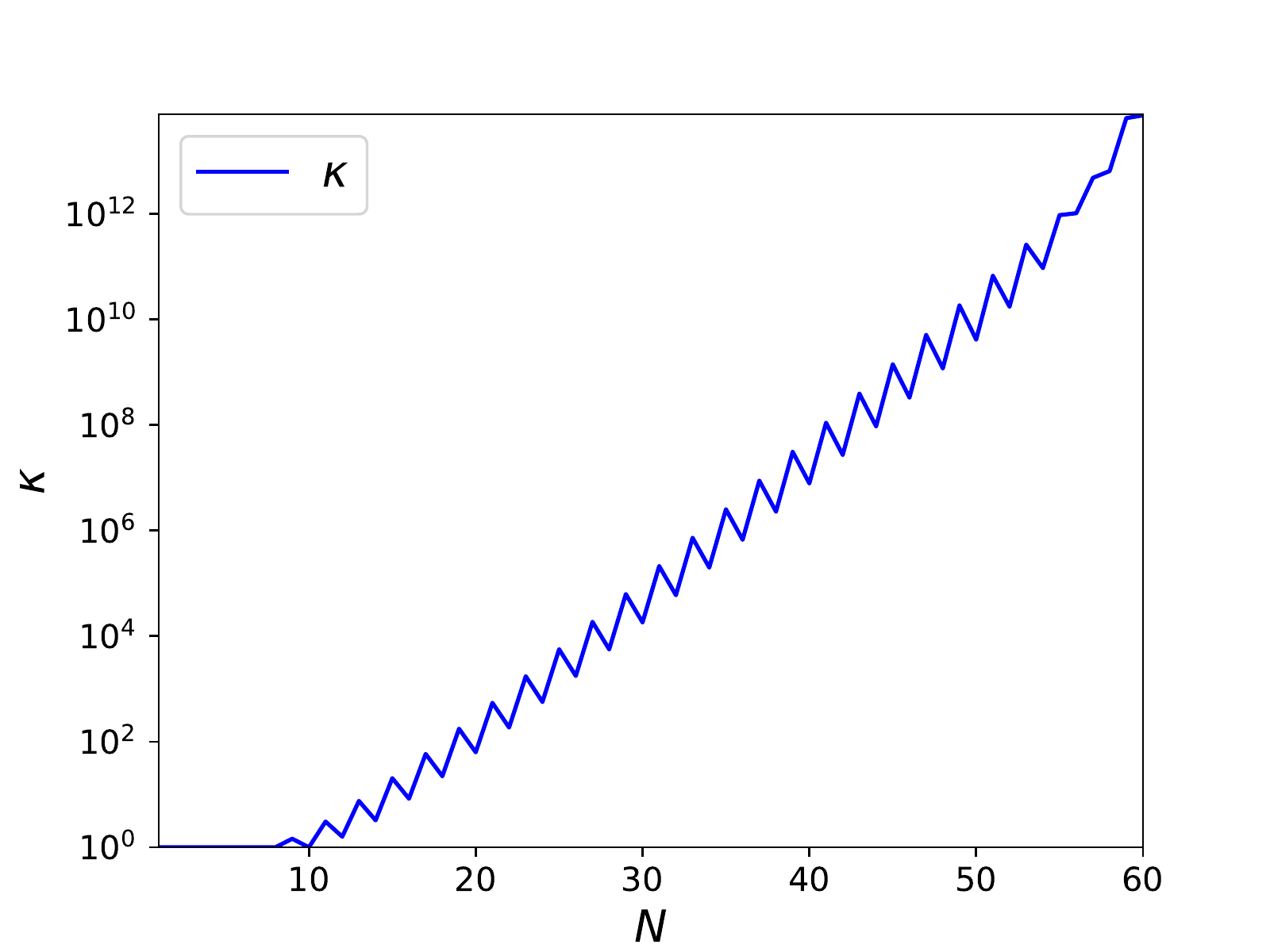}
    \caption{Stability values $\kappa(\vec{\omega})$ for Newton--Cotes QRs}
   \label{fig:kappa-NC}
\end{figure}

In numerical integration, one approximates the continuous integral of a function $g$ by finite sums over weighted 
nodal values of $g$ at a set of distinct points, i.\,e.\ 
\begin{equation}
  I[g] := \int_{-1}^1 g(\xi) \intd \xi \approx \sum_{n=0}^N \omega_n g(\xi_n) =: Q_N[g], 
\end{equation} 
where the QR $Q_N$ is uniquely defined by the quadrature points $\{ \xi_n \}_{n=0}^N$ and the quadrature 
weights $\{ \omega_n \}_{n=0}^N$. 
In many applications, such as the discretisation of exact integrals in DG discretisations, the quadrature points are 
chosen as Gauss--Lobatto points and the quadrature weights are obtained by an interpolatory approach, resulting 
in the well-known Gauss--Lobatto QR; see \cite{krylov2006approximate,gautschi2011numerical}. 
The Gauss--Lobatto QR provides a stable and highly accurate approximation and, when available, is often recommended. 
In particular, when at least $N+1$ quadrature points are used, they provide order of exactness $d=2N-1$, i.\,e.\  
polynomials of degree up to $d=2N-1$ are treated exactly and 
\begin{equation}
  Q_N[g] = I[g]
\end{equation}
holds for all $g \in \mathbb{P}_{2N-1}$. 
Choosing $d=2K-1$ allows us to mimic integration by parts on a discrete level for our numerical solution $u \in 
\mathbb{P}_{K}$, since 
\begin{equation}
  Q_N[ u u' ] 
    = \int_{-1}^1 u u' \intd \xi 
    = u^2 \Big|_{-1}^1 - \int_{-1}^1 u' u \intd \xi 
    = u^2 \Big|_{-1}^1 - Q_N[ u' u]
\end{equation}
is satisfied then. 
In this work, we are in need of QRs on equidistant or even scattered points. 
In this case, a first choice are composite Newton--Cotes QRs, such as the trapezoidal QR. 
Unfortunately, such QRs only provide small orders of exactness and will not allow us to construct conservative 
and stable high-order DG discretisations on equidistant and scattered points. 
Another option are general interpolatory QRs, where the QR is obtained by exactly integrating the interpolation 
polynomial $g_N$ corresponding to the data set $\left\{ (\xi_n,g(\xi_n)) \ | \ n=0,\dots,N \right\}$, i.\,e.\ 
\begin{equation}
  Q_N[g] 
    = \int_{-1}^1 g_N(\xi) \intd \xi
\end{equation}
with $g_N \in \mathbb{P}_N$ such that $g_N(\xi_n) = g(\xi_n)$ for $n=0,\dots,N$. 
Interpolatory QRs naturally yield at least order of exactness $N$, since $g_N = g$ for $g \in \mathbb{P}_N$. 
Yet, interpolatory QRs are known to become unstable on equidistant points for $N \to \infty$, see 
\cite{huybrechs2009stable}. 
Let $\vec{\omega} = (\omega_0,\dots,\omega_N)^T$ be a vector of quadrature weights. 
A common measure of stability of a QR is given by 
\begin{equation}
  \kappa(\vec{\omega}) := \sum_{n=0}^N |\omega_n|.
\end{equation} 
We call a QR with weights $\vec{\omega}$ \emph{stable}, if $\kappa$ is uniformly bounded w.\,r.\,t.\ $N$, i.e. 
\begin{equation}
  \sup_{n \in \N} \kappa(\vec{\omega}) = C < \infty. 
\end{equation}
The idea behind this concept is that for a perturbed input $\tilde{g}$ with $|g(\xi)-\tilde{g}(\xi)| \leq \varepsilon$, 
the error of the QR can be estimated by 
\begin{equation}
  \left| Q_N[g] - Q_N[\tilde{g}] \right| 
    \leq \sum_{n=0}^N |\omega_n| \cdot \left| g(\xi_n)-\tilde{g}(\xi_n) \right| 
    \leq \kappa(\vec{\omega}) \varepsilon.
\end{equation}
Thus, round-off errors due to inexact arithmetics are bounded by the factor $\kappa(\vec{\omega})$. 
The best possible stability value is given by $\kappa(\vec{\omega}) = I[1]$ for a consistent QR, i.\,e.\ 
$Q_N[1] = I[1]$, and is obtained when all weights are nonnegative. 
Interpolatory (Newton--Cotes) QRs feature negative weights and their instability intensifies for $N \to \infty$. 
Figure \ref{fig:kappa-NC} illustrates the rising instability of Newton--Cotes QRs as the number of equidistant 
quadrature points $N+1$ increases. 
Here, integration is performed on $[0,1]$.

As a result, Newton--Cotes QRs suffer from a major loss of precision when inexact arithmetics are used, due to 
round-off errors to heavily pollute the computations then. 
To overcome this problem and construct high-order collocation based DG discretisations on equidistant and scattered 
points, we propose to use LS-QRs.  
LS-QRs have been introduced by Wilson in \cite{wilson1970necessary,wilson1970discrete} and --- to the best of our 
knowledge --- only have been revisited by Huybrechs in \cite{huybrechs2009stable}. 
The idea behind these methods is to allow the number of quadrature points, $N+1$, to be greater than the order of 
exactness, $d$. 
Then, the resulting linear system of equations that describe the exactness conditions, 
\begin{equation}\label{eq:ex_cond_system}
  \mat{A} \vec{\omega} = \vec{m} 
  \quad \text{with} \quad 
  \mat{A} = \left( \phi_{k}(\xi_{n}) \right)_{k,n=0}^{d,N} 
  \quad \text{and} \quad 
  \vec{m} = \left( I[\varphi_k] \right)_{k=0}^d,
\end{equation}
becomes underdetermined.  
The set of solutions forms a $(N-d)$-dimensional affine subspace 
\begin{equation}
    W = \left\{ \vec{\omega} \in \R^{N+1} \ | \ \mat{A} \vec{\omega} = \vec{m} \right\},
\end{equation}
which is generated by $(N-d)$ vectors $\vec{\omega}$ corresponding to distinct interpolatory QRs 
utilising only $d+1$ out of the $N+1$ points. 
Further note that for every $\vec{\omega} \in W$, the corresponding QR 
\begin{equation}
    Q_N(g) = \sum_{n=0}^N \omega_n g(\xi_n)
\end{equation}
provides order of exactness $d$. 
We now seek to determine $\vec{\omega}^* \in W$ such that the Euclidean norm 
\begin{equation}
  \norm{ \vec{\omega} }_2 := \sqrt{ \sum_{n=0}^N |\omega_n|^2 } 
\end{equation}
is minimised. 
Thus, $\vec{\omega}^*$ is given by 
\begin{equation}
  \vec{\omega}^* = \argmin_{\vec{\omega} \in W} \norm{ \vec{\omega} }_2, 
\end{equation}
i.\,e.\ as the LS solution of the underdetermined system \eqref{eq:ex_cond_system}. 
This solution is obtained by solving 
\begin{equation}\label{eq:omega*}
  \vec{\omega}^* = \mat{A}^T \vec{v}, 
\end{equation}
where $\vec{v}$ is the unique solution of the normal equation 
\begin{equation}
  \mat{A}^T \mat{A} \vec{v} = \vec{m}. 
\end{equation}
See for instance the book \cite{golub2012matrix} of Golub and Van Loan. 
Note that when we choose a basis of DOPs w.\,r.\,t.\ the discrete inner product 
\begin{equation}
  \scp{f}{g} = \sum_{n=0}^N f(\xi_n) g(\xi_n)
\end{equation}
in \eqref{eq:ex_cond_system}, we get $\mat{A}^T \mat{A} = \mat{I}$ and \eqref{eq:omega*} reduces to 
\begin{equation}
  \vec{\omega}^* = \mat{A}^T \vec{m}. 
\end{equation}
Hence, the weights are explicitly given by 
\begin{equation}
  \omega_n^* = \sum_{k=0}^d \varphi_k(\xi_n) I[\varphi_k]
\end{equation}
for $n=0,\dots,N$. 
It was proven in \cite{wilson1970necessary} that the weights $\vec{\omega}^*$ will be nonnegative if a sufficiently 
large number of quadrature points is used. 
Thus, the above procedure results in a stable\footnote{In fact, the QR will even have the optimal stability value of 
$\kappa = I[1]$.} QR 
\begin{equation}\label{eq:DLS-QR1}
  Q_N[g] = \sum_{n=0}^N \omega_n^* g(\xi_n)
\end{equation} 
with order of exactness $d$ on any set of collocation points. 
To such QRs \eqref{eq:DLS-QR1}, we refer to as \emph{LS-QRs}. 
In the later proposed discretisation of the DG method, we will typically choose $d = 2K$. 
\section{Proposed discretisation of the discontinuous Galerkin method on equidistant and scattered points}
\label{sec:proposed}

\begin{figure}[!htb] 
  \centering
  \begin{subfigure}[b]{0.45\textwidth}
    \includegraphics[width=\textwidth]{%
      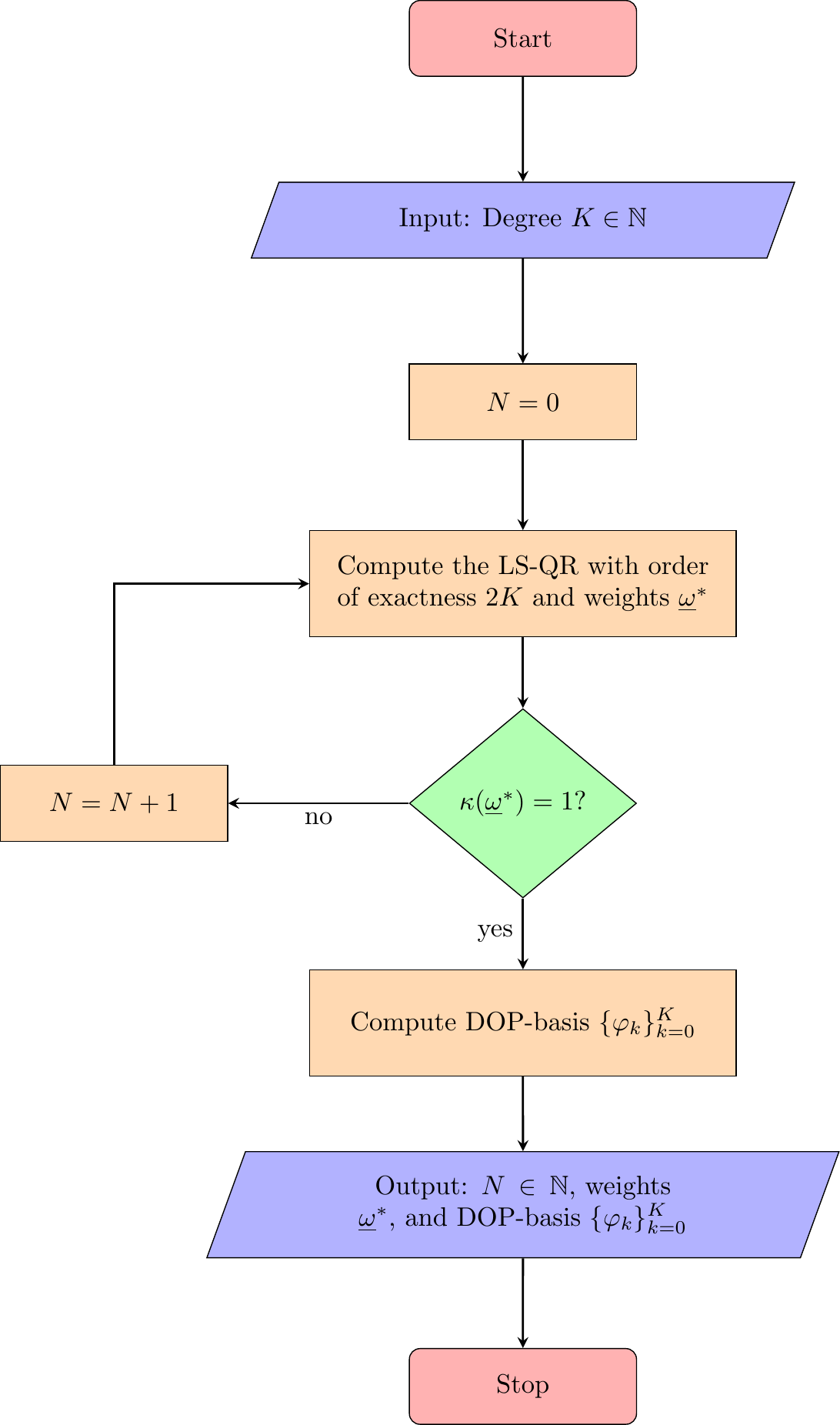}
  \end{subfigure}%
  \caption{Flowchart describing the initial construction of suitable weights $\vec{\omega}^*$ and a DOP-basis $\{ 
\varphi_k \}_{k=0}^K$}
  \label{fig:flowchart}
  \end{figure}

In this section, we propose a stable discretisation for (high-order) DG methods on equidistant and scattered points. 
Therefor, we utilise the techniques discussed in \S \ref{sec:DLS}, involving DLS approximations, LS-QRs, and bases of 
DOPs. 
Starting point of our discretisation is the analytical DG method \eqref{eq:anal-DG}. 
Transformed to the reference element $\Omega_{\mathrm{ref}} = (-1,1)$, the task is to find 
$u^i \in \mathbb{P}_K$ such that 
\begin{equation}\label{eq:trans_an_DG}
  \frac{\Delta x_i}{2} \int_{-1}^1 \dot{u}^i v \intd \xi 
    - \int_{-1}^1 f(u^i) v' \intd \xi 
    + \fnum_R v(1) - \fnum_L v(-1) 
    = 0
\end{equation} 
for all $v \in \mathbb{P}_k$. 
Here, $u^i$ denotes the transformation of the numerical solution $u$, consisting of piecewise polynomials, on the 
element $\Omega_i$ to the reference element. 
Further, $\dot{u}^i$ denotes the temporal derivative $\partial_t u^i$ and $v'$ the spatial derivative $\partial_x v$. 
%
We start our discretisation by replacing $f(u^i)$ in \eqref{eq:trans_an_DG} with a polynomial $f^i \in \mathbb{P}_K$ 
given by a DLS approximation 
\begin{equation}\label{eq:DLS-approx}
  f(u^i) \approx f^i := \sum_{k=0}^K \hat{f}^i_{k,N} \varphi_k,
\end{equation} 
where $\{ \varphi_k \}_{k=0}^K$ is a basis of $\mathbb{P}_K$ and the modal coefficients are given by a simple finite 
sum 
\begin{equation}\label{eq:coeffs}
  \hat{f}^i_{k,N} 
    = \scp{f(u^i)}{\varphi_k}_{\vec{\omega}} 
    = \sum_{n=0}^N \omega_n f(u^i(\xi_n)) \varphi_k(\xi_n).
\end{equation} 
For $N=K$, this is the usual polynomial interpolation. 
Yet, in our discretisation, $N$ might be chosen greater than $K$. 
Next, the involved integrals are replaced by LS-QRs 
\begin{equation}\label{eq:DLS-QR}
  Q_N[g] 
    = \sum_{n=0}^N \omega_n^* g(\xi_n) 
    \approx \int_{-1}^1 g(\xi) \intd \xi
\end{equation}
as discussed in \S \ref{sub:QRs}. 
Utilising the discrete inner product $\scp{\cdot}{\cdot}_{\vec{\omega}^*}$, this results in the discretisation 
\begin{equation}\label{eq:DG-DLS}
  \frac{\Delta x_i}{2} \scp{\dot{u}^i}{v}_{\vec{\omega}^*} 
    = \scp{f^i}{v'}_{\vec{\omega}^*} 
     - \left[ \fnum_R v(1) - \fnum_L v(-1) \right]. 
\end{equation}
$N \in \N$ and $\vec{\omega}^* \in \R^{N+1}$ are chosen such that the resulting LS-QR is stable with 
$\kappa(\vec{\omega}^*) = 1$ and provides order of exactness $2K$. 
This will be crucial to prove conservation and linear stability of the resulting discretisation. 
Further, we follow the idea of collocation and match the quadrature points with the points at which the nodal values of 
$f(u^i)$ are used to construct the DLS approximation $f^i \in \mathbb{P}_K$. 
This results in a more efficient implementation of the proposed discretisation. 
For sake of simplicity, we also use the same weights $\vec{\omega}^*$ for the LS-QR \eqref{eq:DLS-QR} and the DLS 
approximation \eqref{eq:DLS-approx}. 
Different choices are possible but will not be investigated here.
Finally, we include bases of DOPs. 
Note that we can avoid computing a mass matrix on the left hand side of \eqref{eq:DG-DLS} by utilising such a basis of 
DOPs w.\,r.\,t.\ the discrete inner product $\scp{\cdot}{\cdot}_{\vec{\omega}^*}$. 
Since, when choosing $v = \varphi_l$, we have 
\begin{equation}\label{eq:mass_matrix}
  \scp{\dot{u}^i}{v}_{\vec{\omega}^*} 
    = \sum_{k=0}^K \frac{\d}{\d t} \hat{u}_k^i \scp{\varphi_k}{\varphi_l}_{\vec{\omega}^*} 
    = \frac{\d}{\d t} \hat{u}^i_l
\end{equation}
and \eqref{eq:DG-DLS} becomes a system of $K+1$ ODEs 
\begin{equation}\label{eq:DG-DLS-ODE}
  \frac{\Delta x_i}{2} \frac{\d}{\d t} \hat{u}^i_l 
    = \scp{f^i}{\varphi_l'}_{\vec{\omega}^*} 
    - \left[ \fnum_R \varphi_l(1) - \fnum_L \varphi_l(-1) \right]
\end{equation}
for $l=0,\dots,K$. 
The discrete inner product on the right hand side of \eqref{eq:DG-DLS-ODE} is given by 
\begin{equation}\label{eq:v-term}
  \scp{f^i}{\varphi_l'}_{\vec{\omega}^*} 
    = \sum_{k=0}^K \hat{f}^i_{k,N} \scp{\varphi_k}{\varphi_l'}_{\vec{\omega}^*}.
\end{equation}
For the sake of brevity, we will refer to the DLS based discretisation \eqref{eq:DG-DLS-ODE} of the DG method as the 
\textit{discontinuous Galerkin discrete least squares (DGDLS) method}.
A main part of the DGDLS method is to initially determine a suitable vector of weights $\vec{\omega}^*$ and 
a corresponding basis of DOPs. 
Figure \ref{fig:flowchart} provides a flowchart which summarises this procedure. 

\begin{remark}
  Computing the derivative $\frac{\d}{\d t} \hat{u}^i_l$ for a fixed $l \in \{ 0,\dots,K \}$ in the DGDLS method has 
the following complexity: 
  First note that the inner products $\scp{\varphi_k}{\varphi_l'}_{\vec{\omega}^*}$ and 
$\scp{\varphi_k}{\varphi_l'}_{\vec{\omega}^*}$ respectively in \eqref{eq:v-term} and \eqref{eq:coeffs} are computed 
once-for-all a priori. 
  Thus, the computation of a flux coefficients $\hat{f}^i_{k,N}$ is performed in $\mathcal{O}(N)$ operations and the 
whole set $\{\hat{f}^i_{k,N}\}_{k=0}^K$, corresponding to a fixed element $\Omega_i$, is computed in $\mathcal{O}(NK)$ 
operations. 
  This also yields $\frac{\d}{\d t} \hat{u}^i_l$ in the DGDLS method \eqref{eq:DG-DLS-ODE} to be computed in 
$\mathcal{O}(NK)$ operations. 
  In the subsequent numerical tests, we found the choice $N=2K$ to be sufficient for a stable computation on 
equidistant points by the DGDLS method. 
  From point of complexity, this compares to a DG method using over-integration with $2K$ Gauss--Lobatto or 
Gauss--Legendre points and is by a factor $2$ less efficient than a collocation-type DGSEM method using $K+1$ 
Gauss--Lobatto or Gauss--Legendre points. 
\end{remark}

\begin{remark} 
  In many methods it is desirable to avoid mass matrices and, in particular, their inversion, see 
  \cite{abgrall2017high, abgrall2016avoid}. 
  In the proposed discretisation of the DG method this was possible by using bases of DOPs. 
  Continuous Galerkin methods might benefit even more from the combination of DLS approximations and QRs with 
bases of DOPs. 
  Of course, the restriction to a continuous approximation space has to be regarded. 
  Nevertheless, we believe that the application of DOPs could have a positive impact on these schemes. 
  Future work will investigate this possibility.
\end{remark} 
\section{\texorpdfstring{Conservation and $L^2$-stability}{Conservation and L2-stability}}
\label{sec:cons_stab}

In this section, we prove conservation and linear $L^2$-stability of the discretisation proposed in \S 
\ref{sec:proposed}. 
Finally, we also give an outlook how entropy stability can be guaranteed.

\subsection{Conservation}
\label{sub:cons}

Let $u \in V_h$ be a numerical solution consisting of a piecewise polynomial of degree at most $K$ of the scalar 
conservation law 
\begin{equation}\label{eq:cl_cons}
  u_t + f(u)_x = 0
\end{equation}
on $\Omega = [a,b]$. 
Note that the total amount of the conserved variable $u$ only changes due to the flux across the element boundaries, 
i.\,e.\ 
\begin{equation}\label{eq:cons-anal}
  \frac{\d}{\d t} \int_{\Omega} u \intd x = - \left[ f\left(u(b)\right) - f\left(u(a)\right) \right]
\end{equation}
holds for solutions of \eqref{eq:cl_cons}. 
We now show that the numerical solution $u \in V_h$ produced by our discretisation \eqref{eq:DG-DLS} of the DG methods 
also fulfils this property of \emph{conservation}. 

The contribution due to a single element $\Omega_i$, transformed to the reference element, is given 
by 
\begin{equation}
  \frac{\d}{\d t} \int_{\Omega_i} u \intd x
    = \frac{\Delta x_i}{2} \int_{-1}^1 \dot{u}^i \intd \xi, 
\end{equation}
where $u^i = u^i(t,\xi)$ is the transformation of the numerical solution $u = u(t,x)$ on 
${\Omega_i = (x_{i-1/2},x_{i+1/2})}$ transformed to $\Omega_{\mathrm{ref}} = (-1,1)$. 
Further, note that by choosing $v = 1$ in our discretisation \eqref{eq:DG-DLS}, we have 
\begin{equation}
  \frac{\Delta x_i}{2} \scp{\dot{u}^i}{1}_{\vec{\omega}^*} 
    = - \left[ \fnum_{i+1/2} - \fnum_{i-1/2} \right]
\end{equation}
and since the weights $\vec{\omega}^*$ are related to a QR on the collocation points $\{\xi_n\}_{n=0}^N$ 
with order of exactness $2K$, 
\begin{equation}
  \scp{\dot{u}^i}{1}_{\vec{\omega}^*} = Q_N[\dot{u}^i] = \int_{-1}^1 \dot{u}^i \intd \xi 
\end{equation}
holds. 
Thus, we get 
\begin{equation}
  \frac{\d}{\d t} \int_{\Omega_i} u \intd x 
    = \frac{\Delta x_i}{2} \int_{-1}^1 \dot{u}^i \intd \xi 
    = \frac{\Delta x_i}{2} \scp{\dot{u}^i}{1}_{\vec{\omega}^*} 
    = - \left[ \fnum_{i+1/2} - \fnum_{i-1/2} \right].
\end{equation}
Finally, summing up over all elements, the rate of change of the total amount of the numerical solution $u \in V_h$ is 
given by 
\begin{equation}\label{eq:cons-num}
  \frac{\d}{\d t} \int_{\Omega} u \intd x 
    = - \left[ \fnum_b - \fnum_a \right], 
\end{equation}
where $\fnum_b = \fnum_{I+1/2}$ is the numerical flux at the right boundary of $\Omega$ and $\fnum_a = \fnum_{1/2}$ 
is the numerical flux at the left boundary of $\Omega$.
For a consistent single valued numerical flux, \eqref{eq:cons-num} is consistent with \eqref{eq:cons-anal}. 
Note that for periodic boundary conditions, in particular, we have 
\begin{equation}
  \frac{\d}{\d t} \int_{\Omega} u \intd x 
    = 0. 
\end{equation}
This means that for an isolated system, the numerical method is able to exactly mimic the conservation of mass.

\subsection{\texorpdfstring{$L^2$-stability}{L2-stability}}
\label{sub:stab}

Another fundamental design principle for numerical methods is stability. 
It is well-known that the squared $L^2$-norm of physically reasonable solutions of \eqref{eq:cl_cons} with periodic 
boundary conditions does not increase over time, i.\,e.\ 
\begin{equation}\label{eq:stab}
  \frac{\d}{\d t} \norm{ u }_{L^2}^2 \leq 0,
\end{equation}
where a strict inequality reflects the presence of shock waves. 
For scalar conservation laws in one dimension, the squared $L^2$-norm is also referred to as the \emph{energy} or 
the \emph{square entropy}. 
We now prove that for the linear advection equation 
\begin{equation}
  u_t + u_x = 0 
  \quad \text{on} \quad 
  \Omega = [a,b]
\end{equation}
\emph{(strong) $L^2$-stability} \eqref{eq:stab} is also fulfilled by our discretisation \eqref{eq:DG-DLS} of the DG 
method. 
Therefor, let us assume periodic boundary conditions and let $u \in V_h$ be the numerical solution consisting of 
a piecewise polynomial of degree at most $K$. 
We start by noting that the rate of change of the squared $L^2$-norm can be expressed as 
\begin{equation}
  \frac{\d}{\d t} \norm{ u }_{L^2}^2 
    = 2 \int_{\Omega} u_t u \intd x. 
\end{equation}
In a single element, transformed to the reference element, we have 
\begin{equation}\label{eq:L2-1}
  \int_{\Omega_i} u_t u \intd x = \frac{\Delta x_i}{2} \int_{-1}^1 \dot{u}^i u^i \intd \xi.
\end{equation}
This time choosing $v = u^i$ in \eqref{eq:DG-DLS}, results in 
\begin{equation}
  2 \int_{\Omega_i} u_t u \intd x
    = \Delta x_i \scp{\dot{u}^i}{u^i}_{\vec{\omega}^*} 
    = 2\scp{u^i}{(u^i)'}_{\vec{\omega}^*} 
      - 2\left[ \fnum_{i+1/2} u^i(1) - \fnum_{i-1/2} u^i(-1) \right], 
\end{equation}
since $\vec{\omega}^*$ provides a QR with order of exactness $2K$. 
Further, we have 
\begin{equation}
  2\scp{u^i}{(u^i)'}_{\vec{\omega}^*} 
    = 2\int_{-1}^1 u^i \left(u^i\right)_\xi \intd \xi 
    = (u^i)^2(1) - (u^i)^2(-1)
\end{equation}
and thus 
\begin{equation}
  2 \int_{\Omega_i} u_t u \intd x
    = u^i(1) \left[ u^i(1) - 2\fnum_{i+1/2} \right] 
      - u^i(-1) \left[ u^i(-1) - 2\fnum_{i-1/2} \right].
\end{equation}
Finally, summing up over all elements, the global rate of change cuts down to a sum of local contributions 
\begin{equation}
  u_- \left[ u_- - 2 \fnum \right] - u_+ \left[ u_+ - 2 \fnum \right],
\end{equation}
where the interface between two neighbouring elements is considered. 
Here, $u_-$ denotes the value $u^i(1)$ from the left element, $u_+$ denotes the value $u^{i+1}(-1)$ from the 
right element, and $\fnum = \fnum(u_-,u_+)$ denotes the single valued numerical flux at the interface. 
Using a usual full upwind numerical flux, i.e. 
\begin{equation}\label{eq:upwind-flux}
  \fnum(u_-,u_+) = u_-,
\end{equation}
we get 
\begin{equation}
  u_- \left[ u_- - 2 \fnum \right] - u_+ \left[ u_+ - 2 \fnum \right] 
    = -(u_- - u_+)^2 
    \leq 0.
\end{equation}
Hence, the proposed discretisation \eqref{eq:DG-DLS} of the DG method is $L^2$-stable for the linear advection 
equation.

\subsection{An outlook on entropy stability} 
\label{sub:entropy}

Besides $L^2$-stability, often methods are desired that also satisfy the more general property of entropy stability. 
For the DGSEM on Gauss--Lobatto points (including the boundary nodes) a skew-symmetric formulation and SBP operators 
have been used to prove entropy stability 
\cite{gassner2013skew, gassner2016well, gassner2016split}. 
Recently, Chan, Fernandez, and Carpenter \cite{chan2019efficient} proved entropy stability also for a collocation based 
DG method on Gauss--Legendre nodes (not including the boundary nodes) using a decoupled SBP formulation. 
They were further able to extend their approach to general selections of quadrature points and various bases. 
Thus, once we can ensure the SBP property for the proposed DLS discretisation, we could follow the results of Chan, 
Fernandez, and Carpenter to ensure entropy stability also for the proposed DLS based discretisation of the DG method.

Another option to address entropy stability could be the introduction of entropy correction terms as proposed in 
\cite{abgrall2018general}. 
For degree of freedom $k \in \{0,\dots,K\}$ and element $i \in I$, the correction term is given by
\begin{equation}\label{eq:correction_term}
  r_k^i = \alpha ( \hat{u}^i _k- \ubar^i ), 
\end{equation}
where
\begin{equation*}
  \ubar:= \frac{1}{K+1}\sum_{k=0}^K \hat{u}_k^i 
  \quad \text{and} \quad  
  \alpha:= \frac{E}{\sum_{k=0}^N (\hat{u}^i_k-\ubar^i )^2}.
\end{equation*} 
Here, $E$ is the so-called entropy error, which can be calculated using an entropy numerical flux 
$\hat{g}^{\operatorname{num}}$ and an entropy variable $v^i$, which is set equal to the solution 
$u^i$ for sake of simplicity.  
Yet, the correction term \eqref{eq:correction_term} is also valid for general entropy variables $v^i$, see 
\cite{abgrall2018general}.
Using \eqref{eq:DG-DLS-ODE}, $E$ is given by
\begin{equation}\label{eq:entropy_error}
E:= \  \left[ \hat{g}^{\operatorname{num}}_R(1) - \hat{g}^{\operatorname{num}}_L (-1)) \right] -   \frac{2}{\Delta x_i} 
\sum_{k=0}^K \hat{u}_k^i \left(\scp{f^i}{(u^i)'}_{\vec{\omega}^*} 
    - \left[ \fnum_R u^i(1)) - \fnum_L u^i(-1) \right]\right). 
\end{equation} 
The correction term is consistent with zero and does not effect the conservation relation since 
\begin{equation}
\sum_{k=0}^Kr^i_k = \alpha\left( \sum_{k=0}^K (\hat{u}_k^i -\ubar) \right)=0 
\end{equation} 
holds. 
The correction term \eqref{eq:correction_term} is added to the right hand side of the DGDLS scheme 
\eqref{eq:DG-DLS-ODE} and results in an entropy stable scheme. 
This idea was already applied in \cite{abgrall2018connection} to construct entropy stable flux reconstruction schemes on 
polygonal meshes.
Future works will address both approaches to construct entropy stable (DLS based) discretisations of DG methods as well 
as their comparison.
\section{Numerical results}
\label{sec:num}

In the subsequent numerical tests, we investigate conservation, $L^2$-stability, and approximation properties of the 
DGDLS methods on equidistant and scattered points. 
We will do so for a linear advection equation in \S \ref{sub:lin_adv} and a nonlinear inviscid Burgers' equation in \S 
\ref{sub:Burgers}. 
Finally, \S \ref{sub:system} demonstrates the extension to systems of conservation laws 
 and \S \ref{sub:2d} addresses a variable coefficients problem in two space dimensions.

\subsection{Linear advection equation}
\label{sub:lin_adv}

Let us consider the linear advection equation 
\begin{equation}\label{eq:lin-adv2}
  u_t + u_x = 0
\end{equation}
on $\Omega=[0,1]$ with a smooth initial condition $u_0(x) = \sin( 4 \pi x)$ and periodic boundary conditions. 
By the method of characteristics, the (entropy) solution is simply given by $u(t,x) = u_0(x-t)$. 
In the following we evolve this solution until time $t=1$, so that the solution $u(1,x)$ is equal to the initial 
condition $u_0(x)$.

\subsubsection{\texorpdfstring{Conservation and $L^2$-stability}{Conservation and L2-stability}}
\label{subsub:lin_adv_cons}

\begin{figure}[!htb]
  \centering
  \begin{subfigure}[b]{0.33\textwidth}
    \includegraphics[width=\textwidth]{%
      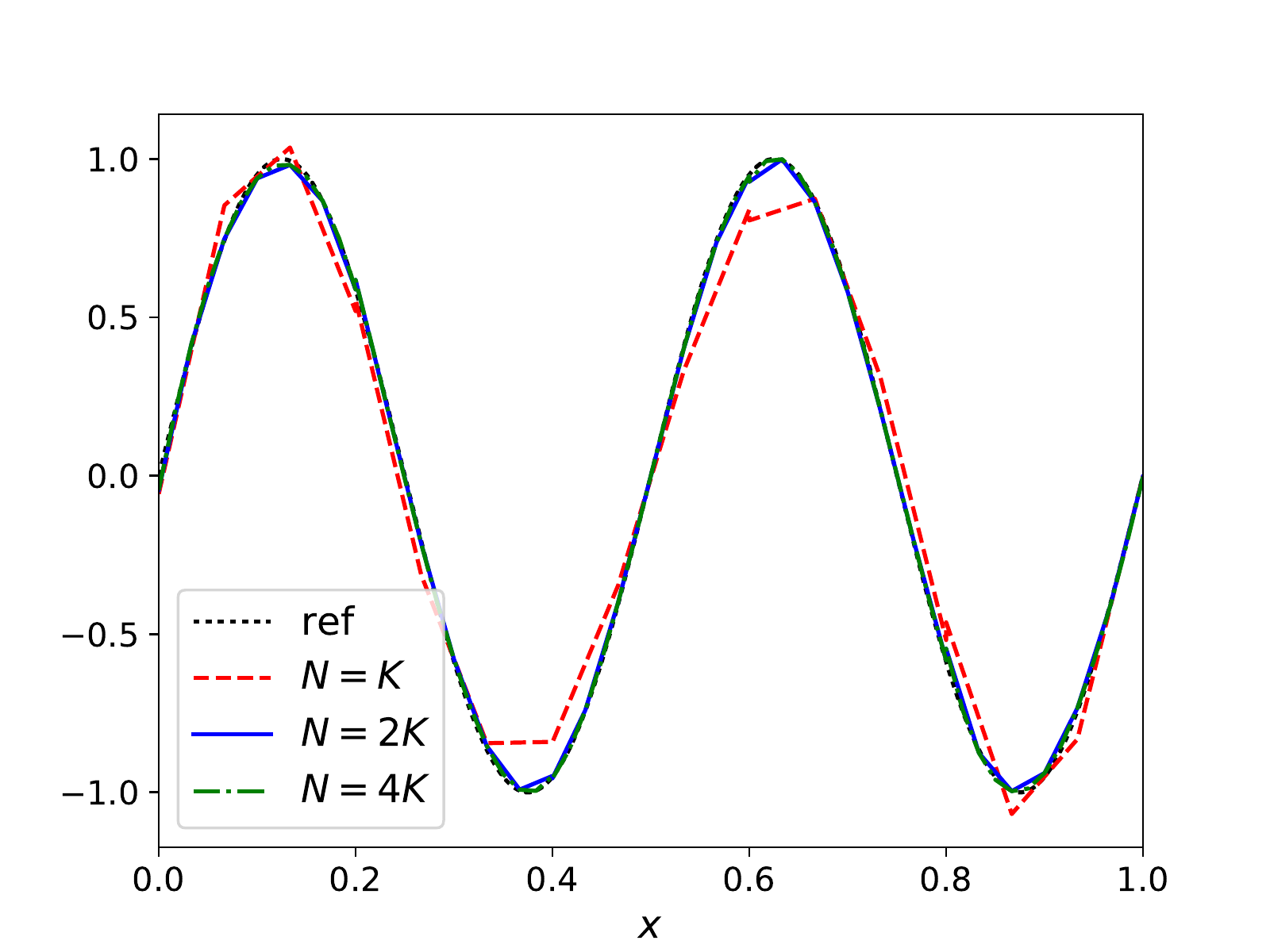}
    \caption{Solution at time $t=1$}
    \label{fig:lin_ad_sol_eq}
  \end{subfigure}%
  ~
  \begin{subfigure}[b]{0.33\textwidth}
    \includegraphics[width=\textwidth]{%
      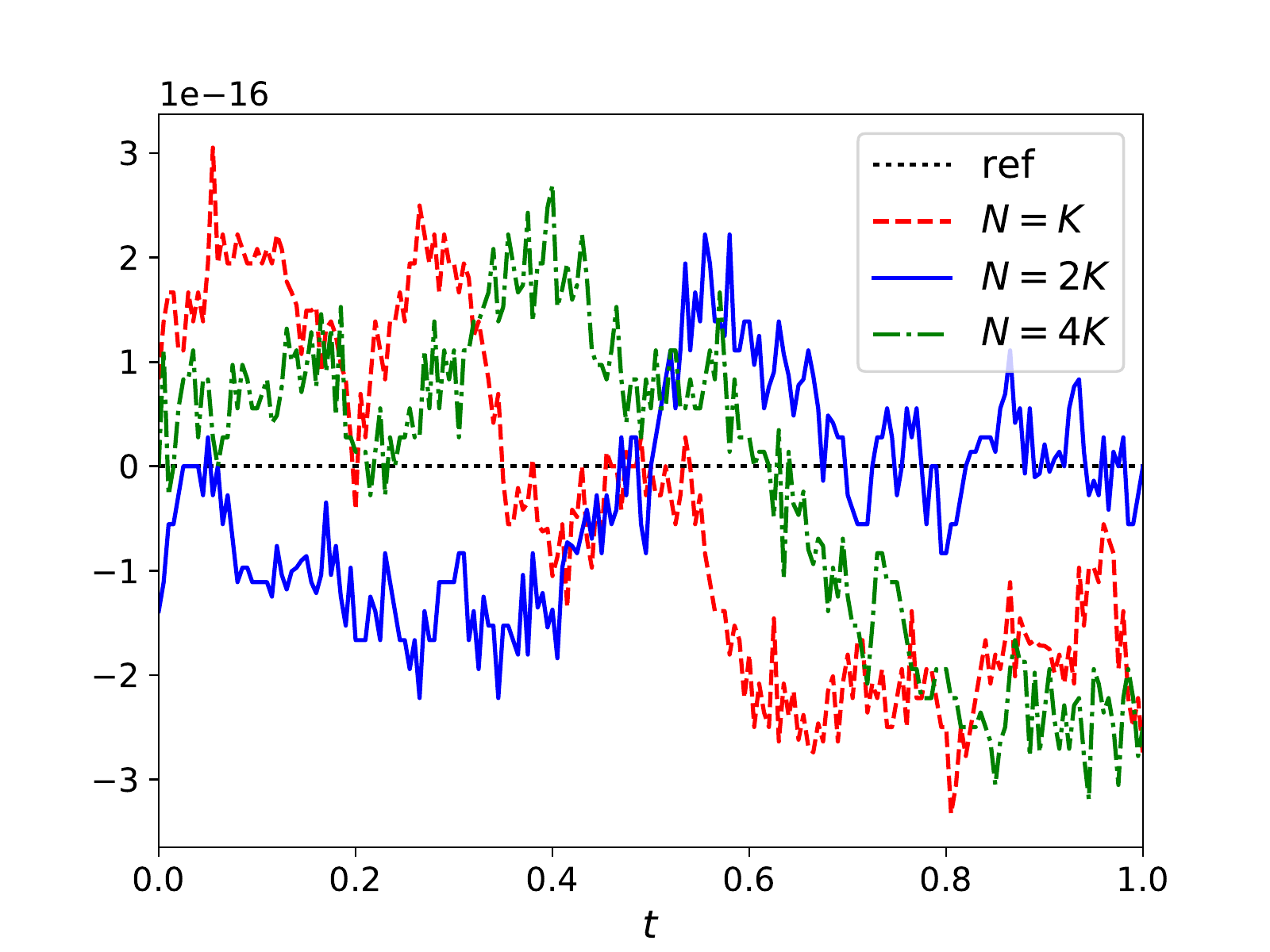}
    \caption{Mass $\int u \intd x$ over time}
    \label{fig:lin_ad_mass_eq}
  \end{subfigure}%
  ~ 
  \begin{subfigure}[b]{0.33\textwidth}
    \includegraphics[width=\textwidth]{%
      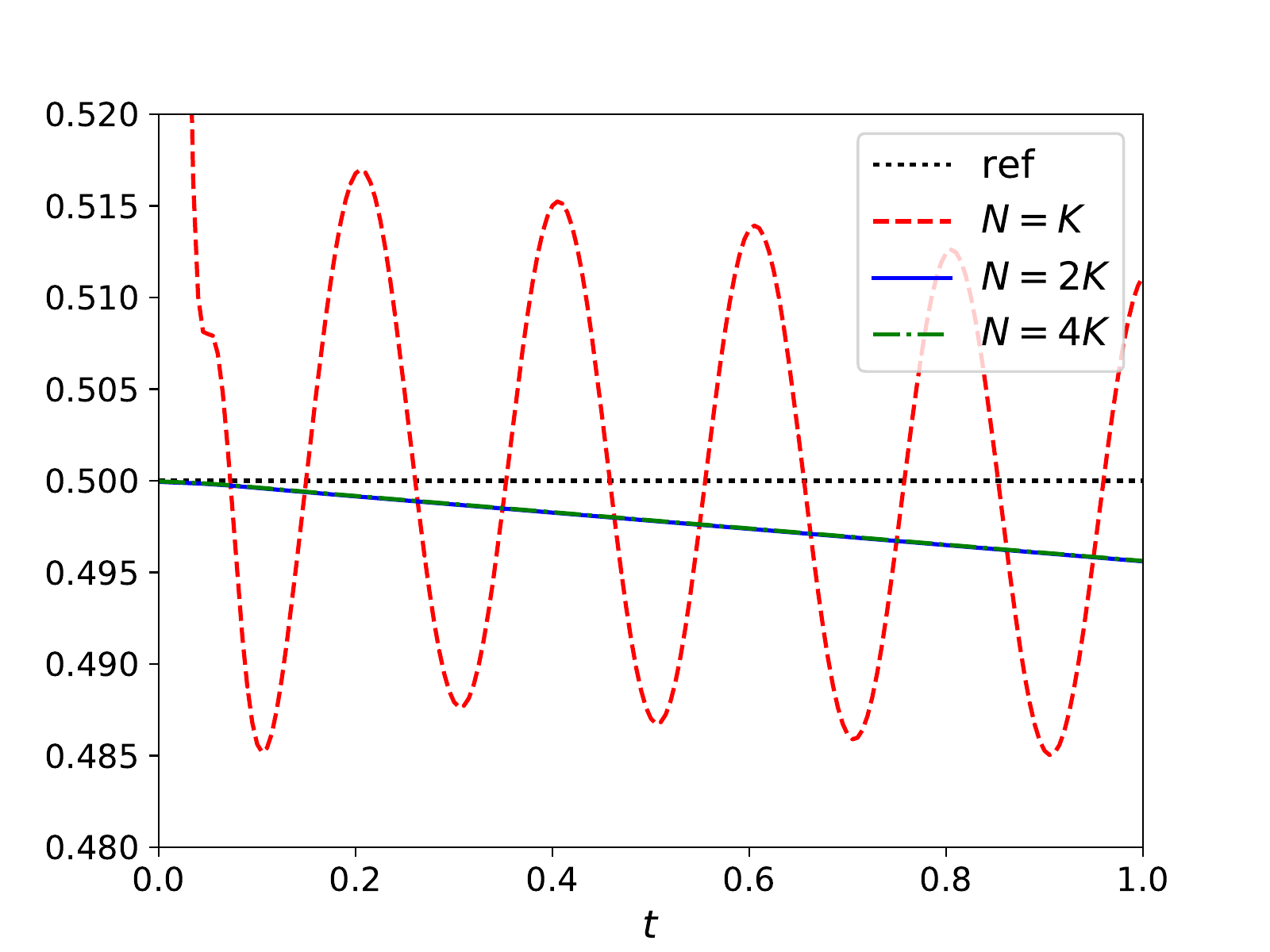}
    \caption{Energy $\int u^2 \intd x$ over time}
    \label{fig:lin_ad_energy_eq}
  \end{subfigure}%
  \caption{Numerical solution, mass, and energy for $I=5$, $K=3$ and $N=K,2K,4K$. Linear advection equation 
and equidistant points.}
  \label{fig:lin_ad_eq}
\end{figure}

We start with a numerical demonstration of conservation and $L^2$-stability of the DGDLS method on equidistant and 
scattered points.
Figure \ref{fig:lin_ad_eq} illustrates the behaviour of the solution, including its mass and energy over time when 
equidistant collocation points are used. 
For this test, we have chosen $I=5$ equidistant elements, a polynomial degree $K=3$, and the full upwind numerical flux 
\eqref{eq:upwind-flux}.

\begin{figure}[!htb]
  \centering
  \begin{subfigure}[b]{0.33\textwidth}
    \includegraphics[width=\textwidth]{%
      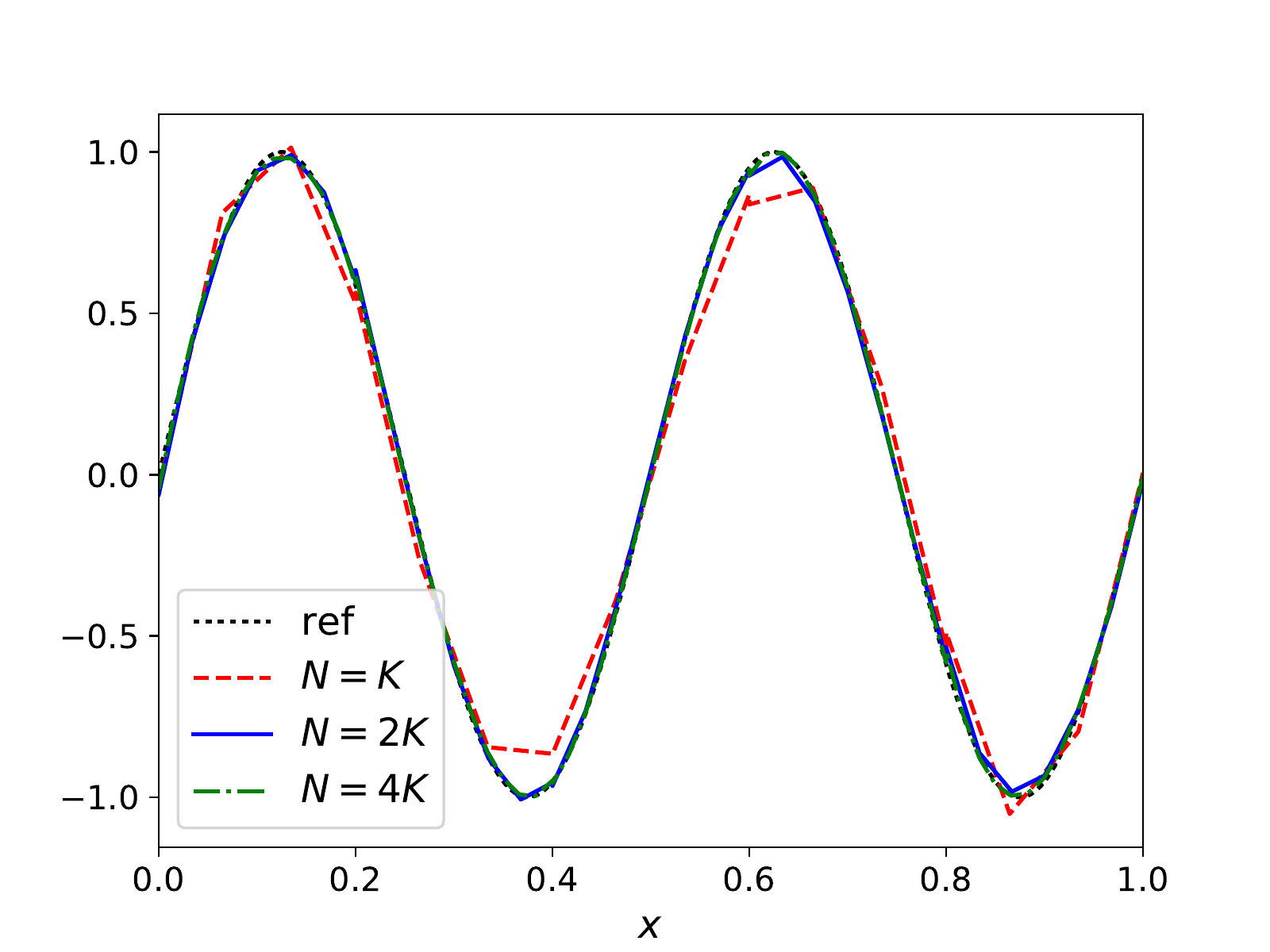}
    \caption{Solution at time $t=1$.}
    \label{fig:lin_ad_sol_rand}
  \end{subfigure}%
  ~
  \begin{subfigure}[b]{0.33\textwidth}
    \includegraphics[width=\textwidth]{%
      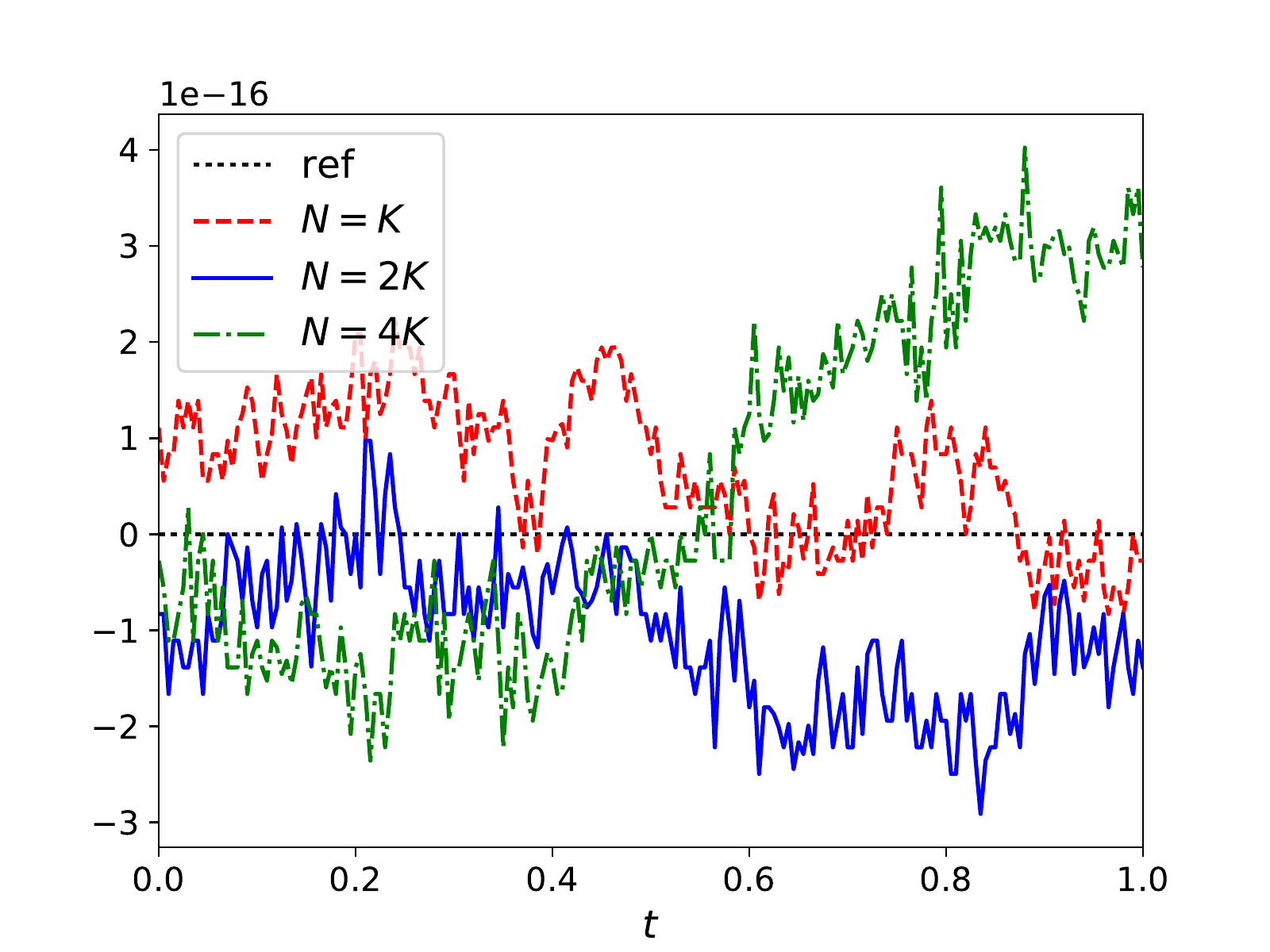}
    \caption{Mass $\int u \intd x$ over time.}
    \label{fig:lin_ad_mass_rand}
  \end{subfigure}%
  ~ 
  \begin{subfigure}[b]{0.33\textwidth}
    \includegraphics[width=\textwidth]{%
      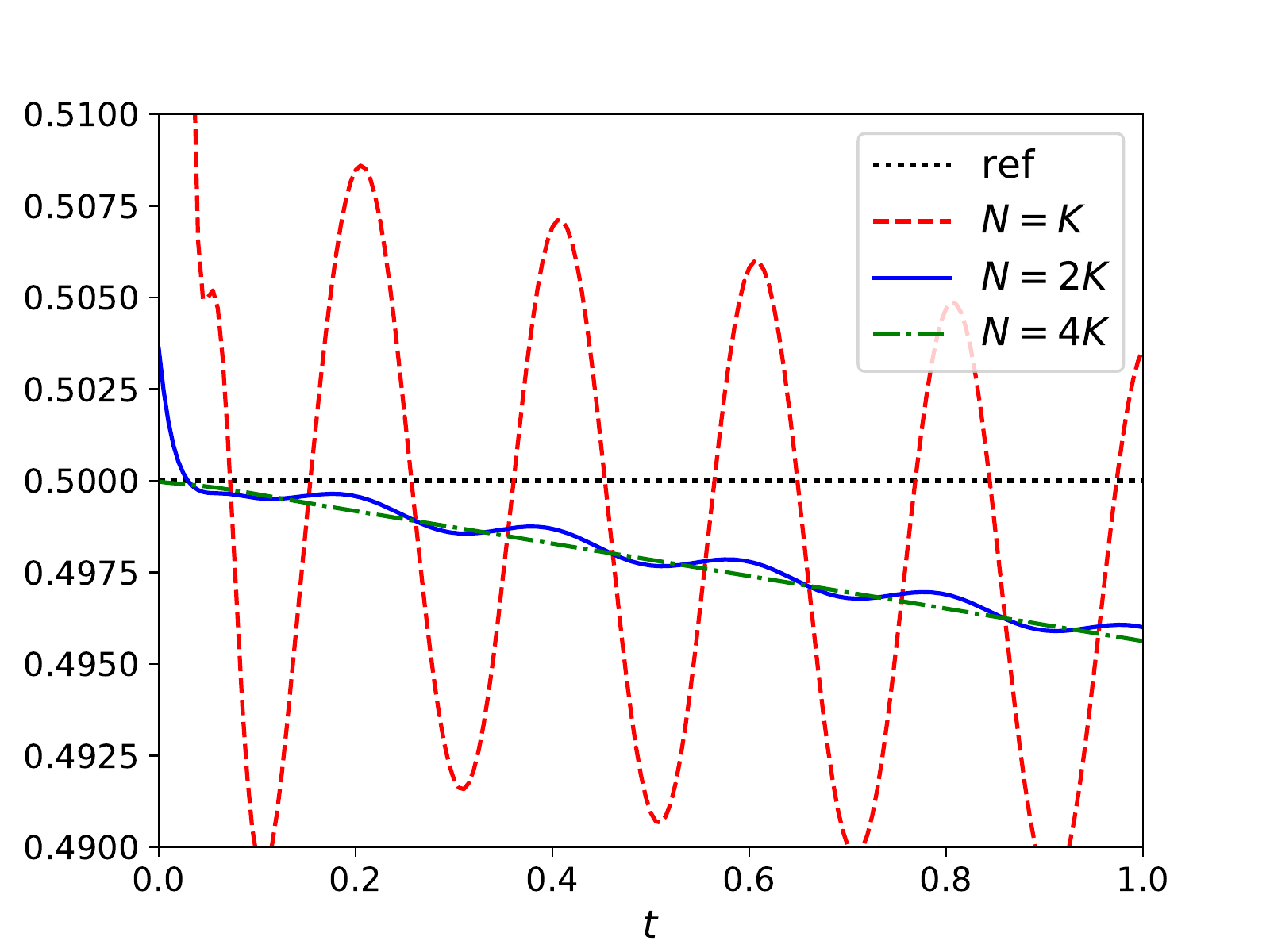}
    \caption{Energy $\int u^2 \intd x$ over time.}
    \label{fig:lin_ad_energy_rand}
  \end{subfigure}%
  \caption{Numerical solution, mass, and energy for $I=5$, $K=3$, and $N=K,2K,4K$. Linear advection equation 
and scattered points.}
  \label{fig:lin_ad_rand}
\end{figure} 

Besides the reference solution, the figures show the numerical solutions for $N=K,2K,4K$, i.\,e.\ for 
an increasing number of equidistant points used in the DLS approximation and the LS-QR. 
All choices yield a conservative method, which can be observed in Figure \ref{fig:lin_ad_mass_eq}. 
Yet, in Figure \ref{fig:lin_ad_energy_eq}, we observe the discretisation to be (strongly) $L^2$-stable only for $N \geq 
2K$. 
Due to an insufficiently high order of exactness of the LS-QR, the numerical solution for $N=K$ yields 
spurious oscillation in the energy of the solution over time. 
As a consequence, the discretisation for $N=2K$ and $N=4K$ can be observed to provide more accurate numerical solutions 
in Figure \ref{fig:lin_ad_sol_eq}. 
We only observe slight differences between the numerical solutions for $N=2K$ and $N=4K$.
Figure \ref{fig:lin_ad_rand} provides a similar demonstration for a set of scattered collocation points.
Here, the scattered collocation points are obtained by adding white uniform noise to the set of equidistant collocation 
points. 
Thus, the scattered collocation points are given by 
\begin{equation}\label{eq:noise}
  \tilde{\xi}_0 = -1, \quad 
  \tilde{\xi}_N = 1, \quad
  \tilde{\xi}_n = \xi_n +  Z_n 
  \quad \text{with} \quad 
  \xi_n = -1 + \frac{2n}{N} 
  \quad \text{and} \quad 
  Z_n \in \mathcal{U}\left(-\frac{1}{40N},\frac{1}{40N}\right),  
\end{equation}
for $n=1,\dots,N-1$, where the $Z_n$ are independent, identically distributed, and further assumed to not be correlated 
with the $\xi_n$. 
Using $N=K,2K,4K$ scattered collocation points, again, all DLS discretisations yield a 
conservative method, see Figure \ref{fig:lin_ad_mass_rand}. 
This time however, slight oscillations can be observed in the energy profile displayed in Figure 
\ref{fig:lin_ad_energy_rand} even for $N=2K$ collocation points. 
These oscillations only vanish when going over to a greater set of $N=4K$ scattered collocation points. 
This behaviour is again caused by an insufficiently high order of exactness of the QR when scattered 
points are used. 
Yet, $N=2K$ and $N=4K$ provide similar accurate numerical solutions, as can be observed in Figure 
\ref{fig:lin_ad_sol_rand}.

\subsubsection{Accuracy and convergence}
\label{subsub:lin_adv_accuracy}

Next, we investigate the approximation properties of the proposed discretisation by an error analysis. 
We consider the same problem as before. 
Table \ref{tab:linear_sin_equid} lists the $L^2$-errors for the DGDLS  
method on equidistant collocation points and an increasing number of degrees of freedom. 
Table \ref{tab:linear_sin_rand} provides the same analysis for sets of scattered collocation points, which are again 
constructed by adding white uniform noise; see \eqref{eq:noise}.
The experimental orders of convergence (EOCs) have been computed by performing a least squares fit for the 
parameters $C$ and $s$ in the model $y = C \cdot N^{-s}$, where $y$ denotes the $L^2$-error for a fixed $N$.

\begin{table}[htb!]
\parbox{.49\textwidth}{
    \centering 
    \begin{adjustbox}{width=0.49\textwidth}
    \begin{tabular}{c c c c c c} 
        \hline
        \multicolumn{2}{c}{}  
        & \multicolumn{4}{c}{$L^2$-errors}
        \\ \hline
	\multicolumn{2}{c}{} & DGSEM & \multicolumn{3}{c}{DGDLS}
        \\ \hline
        $K$ & $I$ & (GL points) & $N=K$ \ & $N=2K$ \ & $N=4K$ \ \,
        \\ \hline
        1 &  5 & 5.8E-1  
	       & 5.8E-1  
               & 6.4E-1  
               & 6.4E-1 \\
          & 10 & 1.0E-1  
	       & 1.0E-1  
               & 2.0E-1  
               & 1.9E-1 \\
          & 20 & 2.6E-2  
	       & 2.6E-2  
               & 3.5E-2  
               & 3.3E-2 \\
          & 40 & 9.6E-3  
	       & 9.6E-3  
               & 6.4E-3  
               & 5.9E-3 \\ 
        \multicolumn{2}{c}{EOC:} & 2.5 & 2.5 & 1.7 & 1.8 \\
        \hline
        2 &  5 & 6.6E-2  
	       & 6.6E-2  
               & 1.0E-1  
               & 9.9E-2 \\
          & 10 & 1.0E-2  
	       & 1.0E-2  
               & 8.7E-3  
               & 7.9E-3 \\
          & 20 & 1.3E-3  
	       & 1.3E-3  
               & 1.0E-3  
               & 9.0E-4 \\
          & 40 & 1.6E-4  
	       & 1.6E-4  
               & 1.2E-4  
               & 1.1E-4 \\ 
        \multicolumn{2}{c}{EOC:} & 2.7 & 2.7 & 3.5 & 3.6 \\ 
        \hline
        3 &  5 & 1.1E-2  
	       & 8.1E-2  
               & 1.0E-2  
               & 8.9E-3 \\
          & 10 & 7.6E-4  
	       & 2.0E-2  
               & 6.3E-4  
               & 5.4E-4 \\
          & 20 & 4.9E-5  
	       & 3.8E-4  
               & 4.0E-5  
               & 3.4E-5 \\
          & 40 & 2.9E-6  
	       & 3.7E-6  
               & 2.5E-6  
               & 2.1E-6 \\ 
        \multicolumn{2}{c}{EOC:} & 3.8 & 2.1 & 3.9 & 4.0 \\ 
        \hline
        4 &  5 & 1.3E-3  
	       & 1.1E-2  
               & 1.2E-3  
               & 1.0E-3 \\
          & 10 & 5.1E-5  
	       & 4.3E-3  
               & 4.2E-5  
               & 3.4E-5 \\
          & 20 & 2.3E-6  
	       & 9.0E-4  
               & 1.5E-6  
               & 1.2E-6 \\
          & 40 & 1.1E-7  
	       & 9.1E-5  
               & 1.0E-7  
               & 9.7E-8 \\ 
        \multicolumn{2}{c}{EOC:} & 4.1 & 1.5 & 4.1 & 4.1 \\ 
        \hline
    \end{tabular} 
    \end{adjustbox}
    \caption{Linear advection equation \& equidistant points.}
    \label{tab:linear_sin_equid} 
}
\hfill
\parbox{.49\textwidth}{
    \centering 
    \begin{adjustbox}{width=0.49\textwidth}
    \begin{tabular}{c c c c c c c c c c} 
        \hline
        \multicolumn{2}{c}{}  
        & \multicolumn{4}{c}{$L^2$-errors}
        \\ \hline
	\multicolumn{2}{c}{} & DGSEM & \multicolumn{3}{c}{DGDLS}
        \\ \hline
        $K$ & $I$ & (GL points) & $N=4K$ & $N=16K$ & $N=64K$
        \\ \hline
        1 &  5 & 5.8E-1 
               & 6.4E-1   
               & 6.4E-1
               & 6.4E-1 \\
          & 10 & 1.0E-1
               & 1.9E-1  
               & 1.9E-1
               & 1.9E-1 \\
          & 20 & 2.6E-2 
               & 3.3E-2  
               & 3.3E-2
               & 3.3E-2 \\
          & 40 & 9.6E-3
               & 5.8E-3  
               & 5.7E-3 
               & 5.6E-3 \\
        \multicolumn{2}{c}{EOC:} & 2.4 & 1.8 & 1.8 & 1.8 \\ 
        \hline
        2 &  5 & 6.6E-2 
               & 9.9E-2  
               & 9.8E-2
               & 9.8E-2 \\
          & 10 & 1.0E-2 
               & 7.9E-3  
               & 7.7E-3 
               & 7.6E-3 \\
          & 20 & 1.3E-3 
               & 9.0E-4   
               & 8.7E-4 
               & 8.6E-4 \\
          & 40 & 1.6E-4
               & 1.1E-4   
               & 1.0E-4 
               & 1.0E-4 \\
        \multicolumn{2}{c}{EOC:} & 2.7 & 3.6 & 3.6 & 3.6 \\ 
        \hline
        3 &  5 & 1.1E-2 
               & 8.9E-3   
               & 8.7E-3 
               & 8.7E-3 \\
          & 10 & 7.6E-4 
               & 5.6E-4   
               & 5.2E-4 
               & 5.2E-4 \\
          & 20 & 4.9E-5  
               & 3.7E-5   
               & 3.3E-5 
               & 3.3E-5 \\
          & 40 & 2.9E-6
               & 2.1E-6   
               & 2.0E-6 
               & 2.0E-6 \\
        \multicolumn{2}{c}{EOC:} & 3.8 & 3.9 & 4.0 & 4.0 \\ 
        \hline
        4 &  5 & 1.3E-3 
               & 1.1E-3   
               & 1.0E-3  
               & 1.0E-3 \\
          & 10 & 5.1E-5
               & 3.8E-4   
               & 4.3E-5  
               & 3.3E-5 \\
          & 20 & 2.3E-6
               & 3.9E-4   
               & 5.8E-6 
               & 2.0E-6 \\
          & 40 & 1.1E-7
               & 4.6E-6   
               & 6.4E-6 
               & 2.7E-7 \\
        \multicolumn{2}{c}{EOC:} & 4.1 & 1.1 & 4.0 & 4.1 \\ 
        \hline
    \end{tabular}  
    \end{adjustbox}
    \caption{Linear advection equation \& scattered points.}
    \label{tab:linear_sin_rand}
}
\end{table}

We note from both tables that increasing the number of collocation points does not always yield more accurate numerical 
solutions when using a fairly small number of degrees of freedom, i.\,e.\ polynomials of degree $K=1$ or only $I=5$ 
elements. 
Yet, for all higher degrees $K \geq 2$ (and $I>5$ for $K=2$), we observe the numerical solutions to become more 
accurate when the number of collocation points is increased. 
Further, both tables provide a comparison of the DGDLS method with the usual DGSEM on a set of Gauss--Lobatto points. 
The DGSEM can be considered as a special case of the proposed DGDLS method when Gauss--Lobatto points are used and 
$N=K$ is chosen, i.\,e.\ using polynomial 
interpolation as well as an interpolatory (Gauss--Lobatto) QR. 
Of course, Gauss--Lobatto points and their corresponding QR are known to be superior to equidistant or 
even scattered points, at least when the same number of points is used. 
Yet, when sufficiently increasing the number of collocation points, we often observe the DGDLS method to provide 
more accurate results on equidistant and even scattered collocation points than the DGSEM on Gauss--Lobatto 
points.

\begin{remark}
  In some cases, we observe the error to increase even though the number of elements $I$ is increased, e.\,g.\ for 
$K=4$ and $N=4K$ when going over from $I=10$ to $I=20$ as well as for $K=4$ and $N=16K$ when going over from $I=20$ to 
$I=40$ in Table \ref{tab:linear_sin_rand}. 
  The same observation can be made for the subsequent tables \ref{tab:burgers_t1_equid} and \ref{tab:burgers_t1_rand} 
concerning the inviscid Burgers' equation. 
  This behaviour is caused by the instability of the LS-QR when an insufficiently great number of collocation points 
$N$ is used.  
  Hence, we observe this problem to vanish when $N$ is increased; see the case $N=64$ in tables 
\ref{tab:linear_sin_rand} and \ref{tab:burgers_t1_rand}.
\end{remark}

\subsection{Inviscid Burgers' equation}
\label{sub:Burgers} 

Let us now consider the nonlinear inviscid Burgers' equation 
\begin{equation}
  u_t + \left( \frac{u^2}{2} \right)_x = 0 
\end{equation}
on $\Omega=[0,1]$ with smooth initial condition ${u_0(x) = 1+ \frac{1}{4 \pi} \sin( 2 \pi x )}$
and periodic boundary conditions. 
For this problem a shock develops in the solution when the wave breaks at time 
\begin{equation}
  t_b = - \frac{1}{\min_{0 \leq x \leq 1} u_0'(x)} = 2.
\end{equation}
In the subsequent numerical tests we consider the solution at times $t=1$ and $t=3$. 
This time, a local Lax--Friedrichs (LLF) numerical flux 
\begin{equation}\label{eq:LLF-flux}
  \fnum_{LLF}(u_-,u_+) 
    = \frac{1}{2} \left[ f(u_-) + f(u_+) \right] 
    - \frac{\lambda}{2} \cdot (u_+ - u_-) 
    \quad \text{with} \quad 
    \lambda = \max\{ |u_-|, |u_+| \}
\end{equation}
is used at the element interfaces. 
The reference solutions have been computed using characteristic tracing, solving the implicit equation $u(t,x) = 
u_0(x - tu)$ in smooth regions. 
The jump location, separating these regions, can be determined by the Rankine--Hugoniot condition.

\subsubsection{\texorpdfstring{Conservation and $L^2$-stability}{Conservation and L2-stability}}
\label{subsub:Burgers_cons}

\begin{figure}[!htb]
  \centering
  \begin{subfigure}[b]{0.24\textwidth}
    \includegraphics[width=\textwidth]{%
      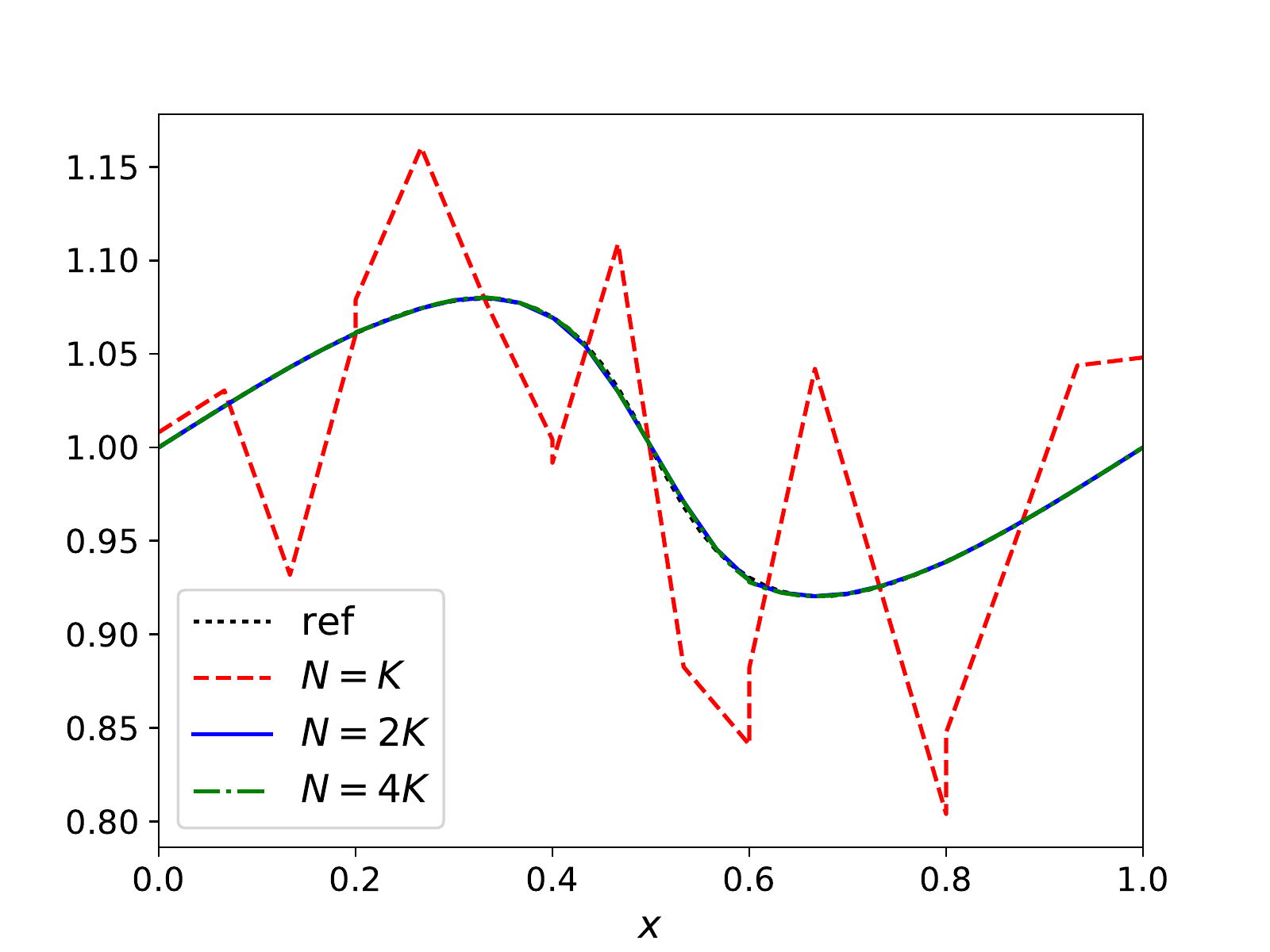}
    \caption{Solution at time $t=1$.}
    \label{fig:burgers_t1_sol_eq}
  \end{subfigure}%
  ~
  \begin{subfigure}[b]{0.24\textwidth}
    \includegraphics[width=\textwidth]{%
      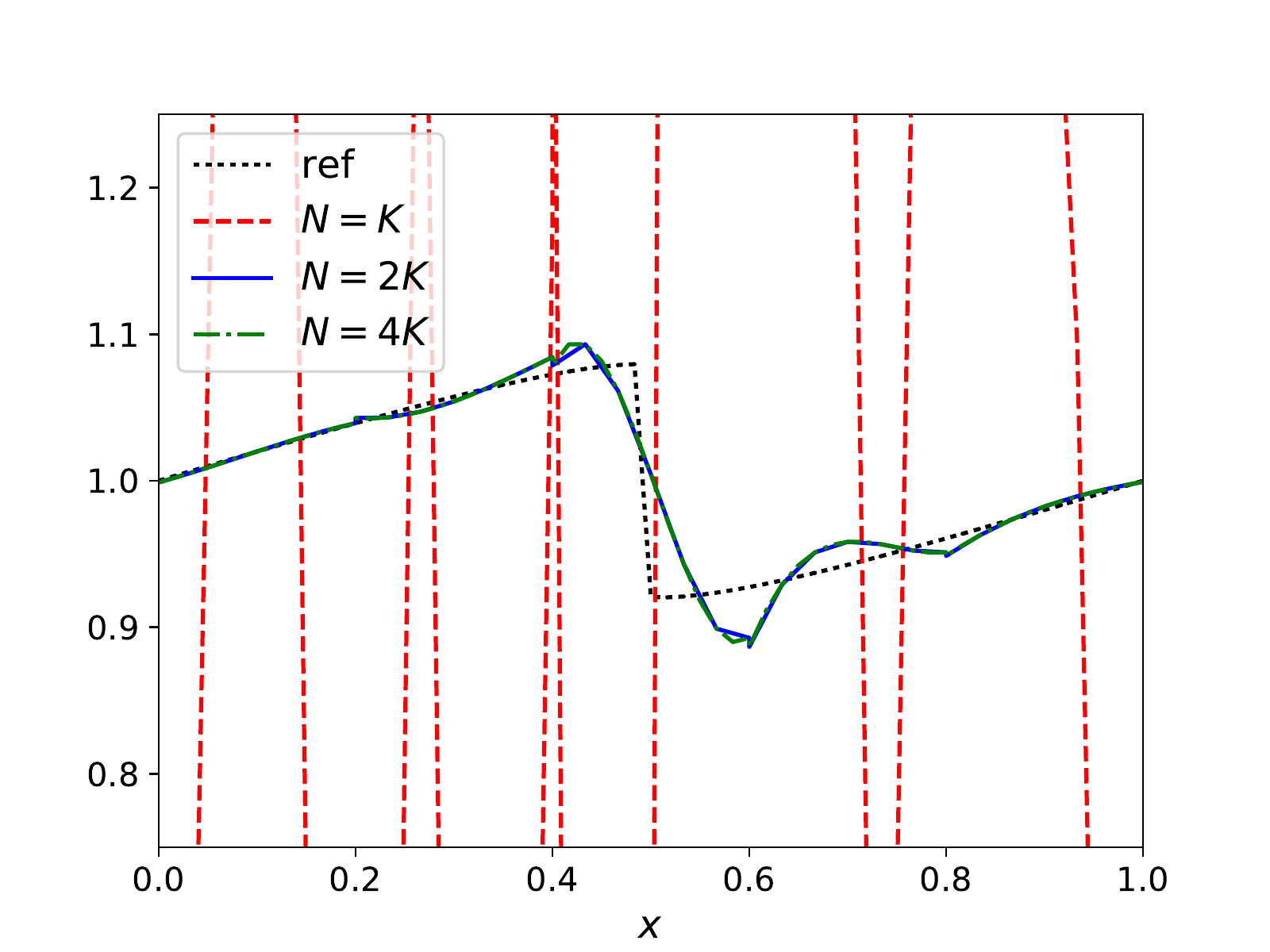}
    \caption{Solution at time $t=3$.}
    \label{fig:burgers_t3_sol_eq}
  \end{subfigure}%
  ~
  \begin{subfigure}[b]{0.24\textwidth}
    \includegraphics[width=\textwidth]{%
      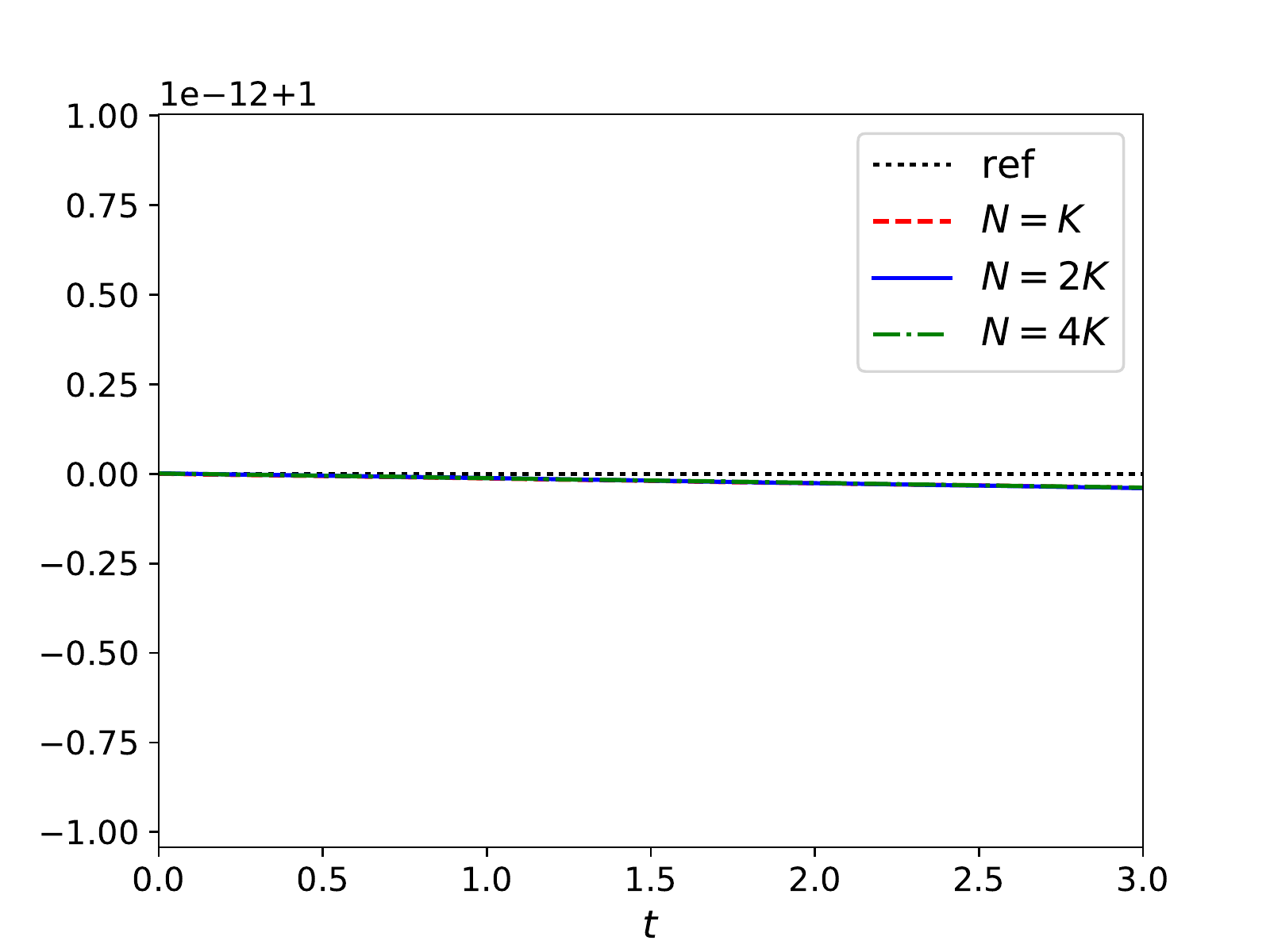}
    \caption{Mass over time.}
    \label{fig:burgers_t3_ad_mass_eq}
  \end{subfigure}%
  ~ 
  \begin{subfigure}[b]{0.24\textwidth}
    \includegraphics[width=\textwidth]{%
      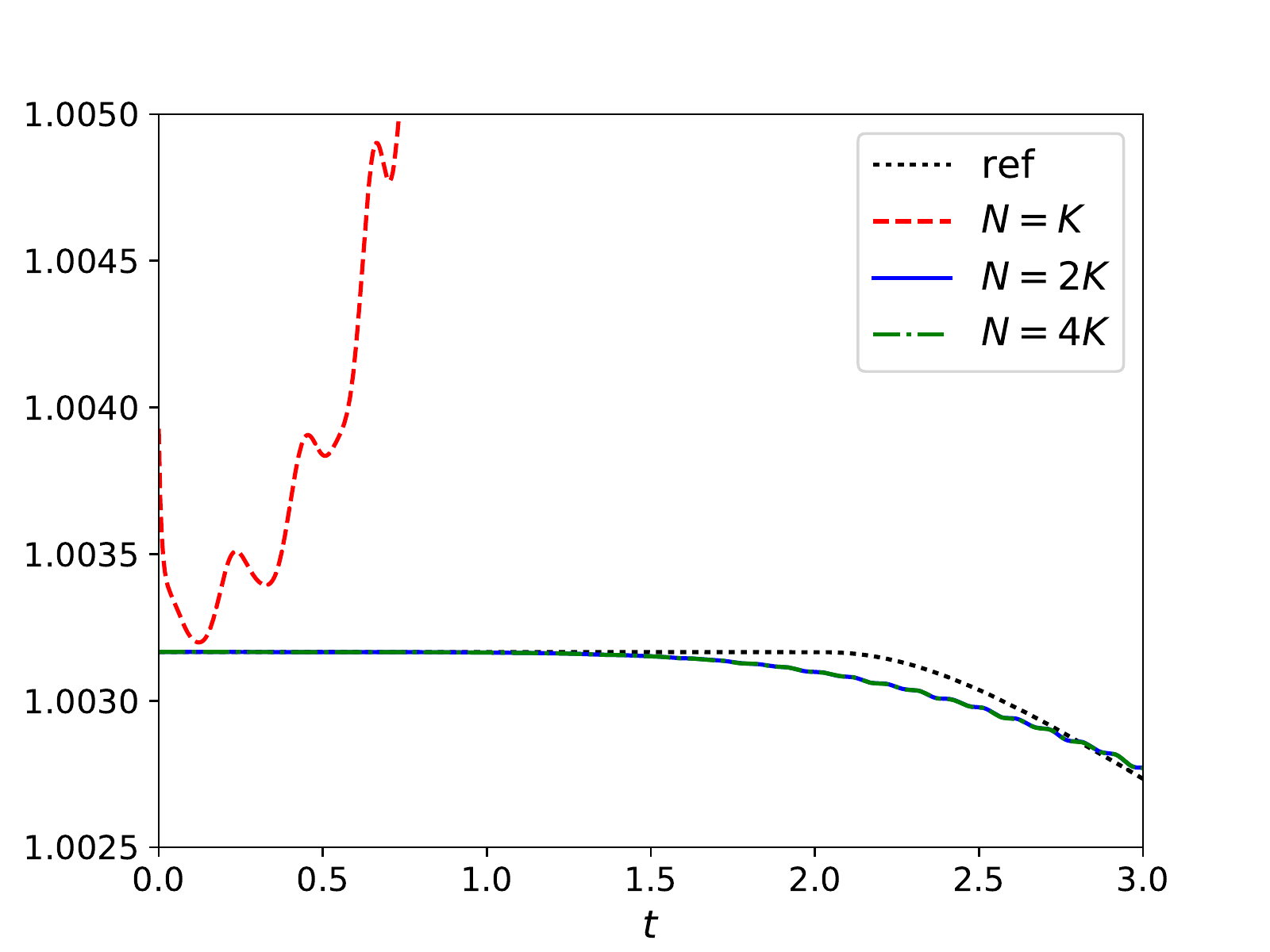}
    \caption{Energy over time.}
    \label{fig:burgers_t3_ad_energy_eq}
  \end{subfigure}%
  \caption{Numerical solution, mass, and energy for $I=5$, $K=3$, and $N=K,2K,4K$. Burgers' equation and equidistant 
points.}
  \label{fig:burgers_ad_eq}
\end{figure}

\begin{figure}[!htb]
  \centering
  \begin{subfigure}[b]{0.24\textwidth}
    \includegraphics[width=\textwidth]{%
      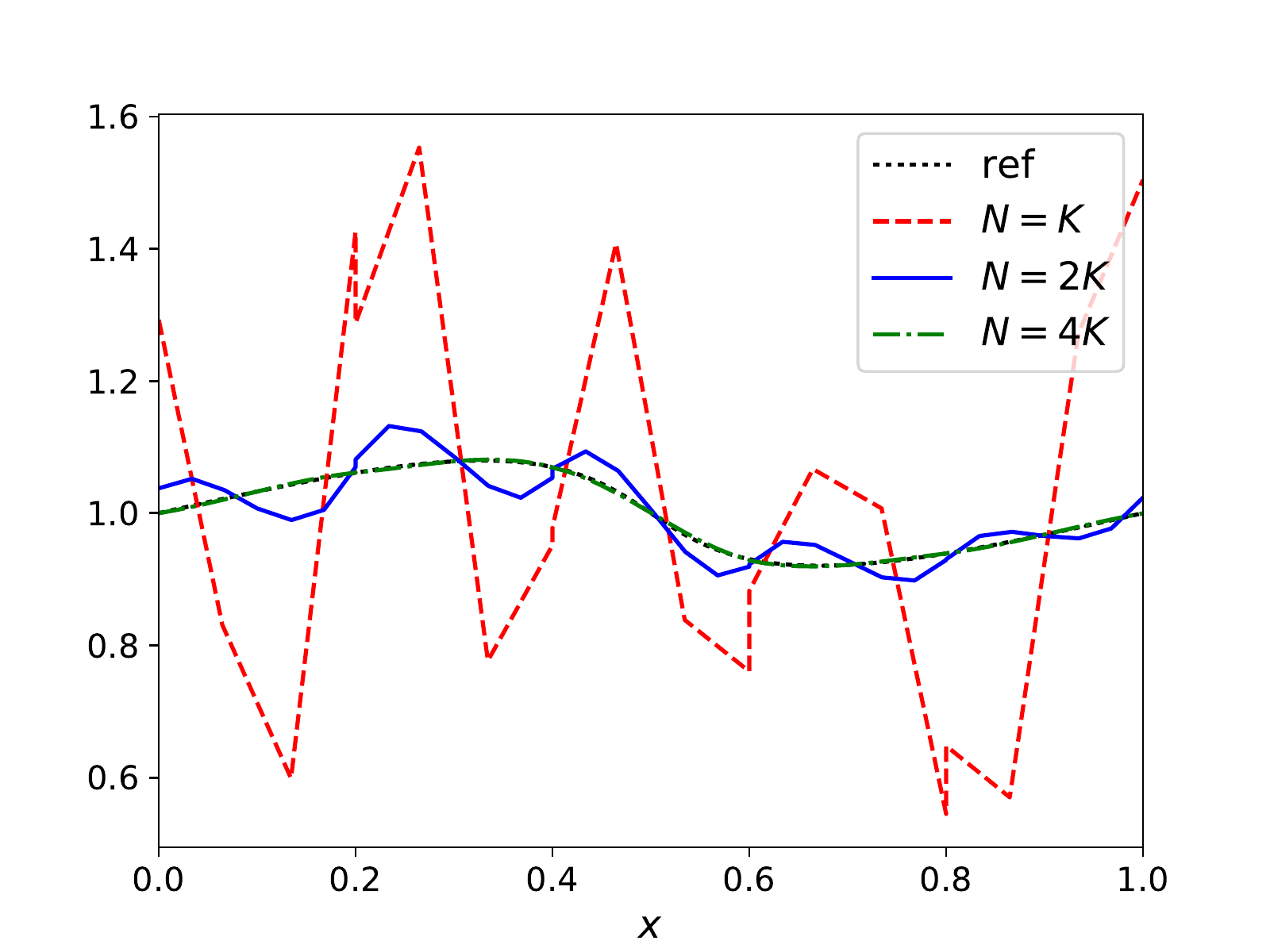}
    \caption{Solution at time $t=1$.}
    \label{fig:burgers_t1_sol_rand}
  \end{subfigure}%
  ~
  \begin{subfigure}[b]{0.24\textwidth}
    \includegraphics[width=\textwidth]{%
      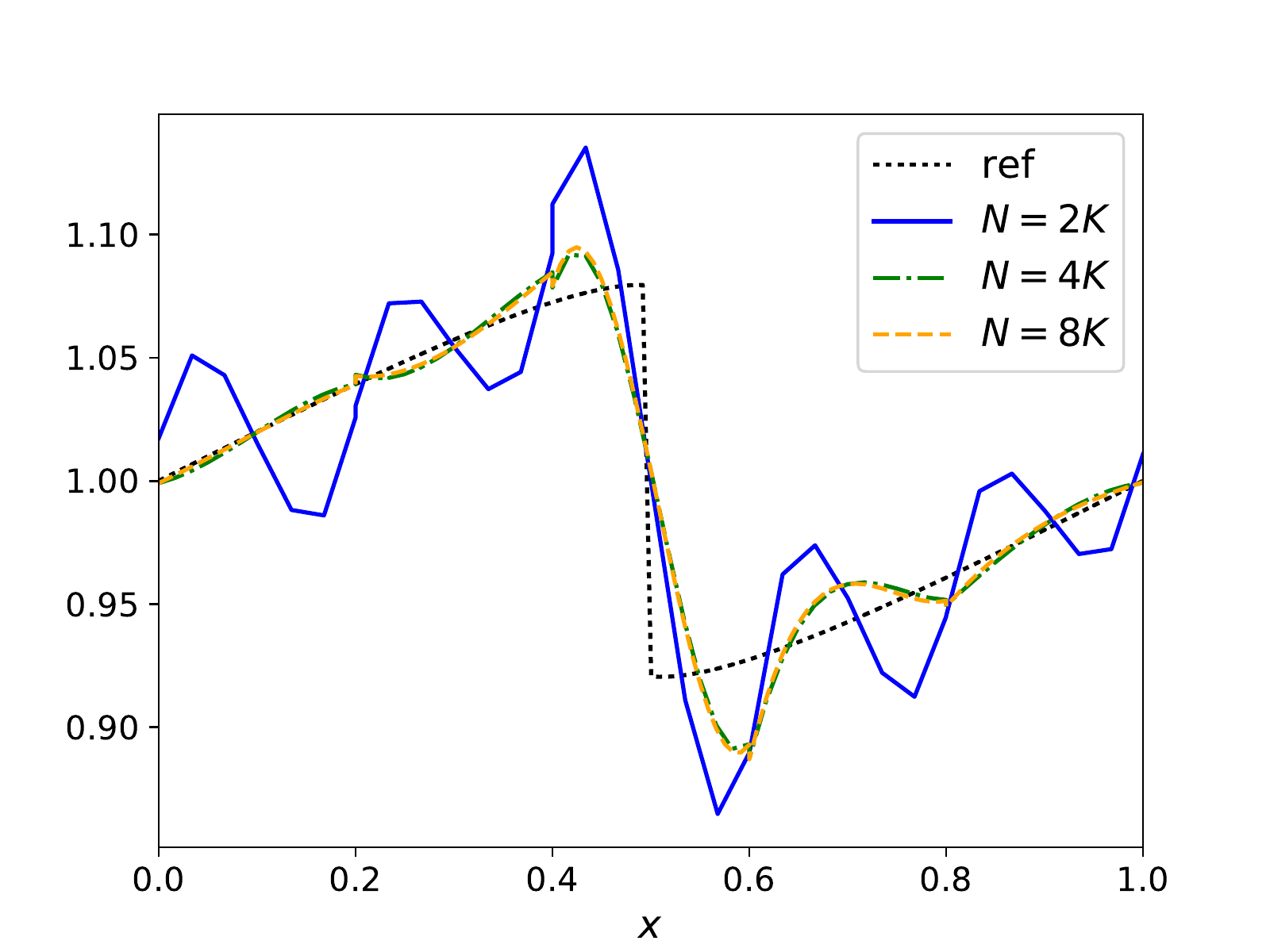}
    \caption{Solution at time $t=3$.}
    \label{fig:burgers_t3_sol_rand}
  \end{subfigure}%
  ~
  \begin{subfigure}[b]{0.24\textwidth}
    \includegraphics[width=\textwidth]{%
      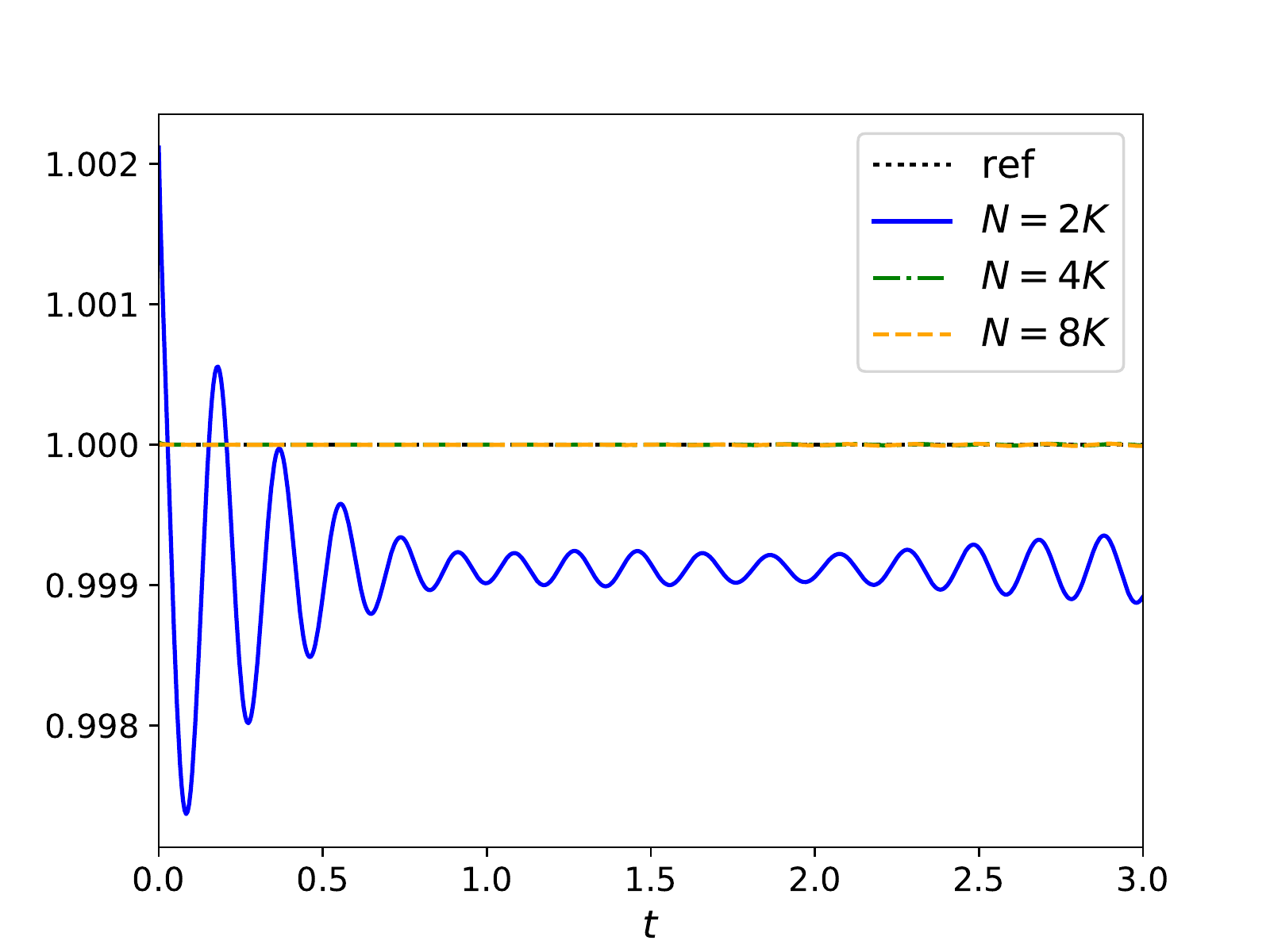}
    \caption{Mass over time.}
    \label{fig:burgers_t3_ad_mass_rand}
  \end{subfigure}%
  ~ 
  \begin{subfigure}[b]{0.24\textwidth}
    \includegraphics[width=\textwidth]{%
      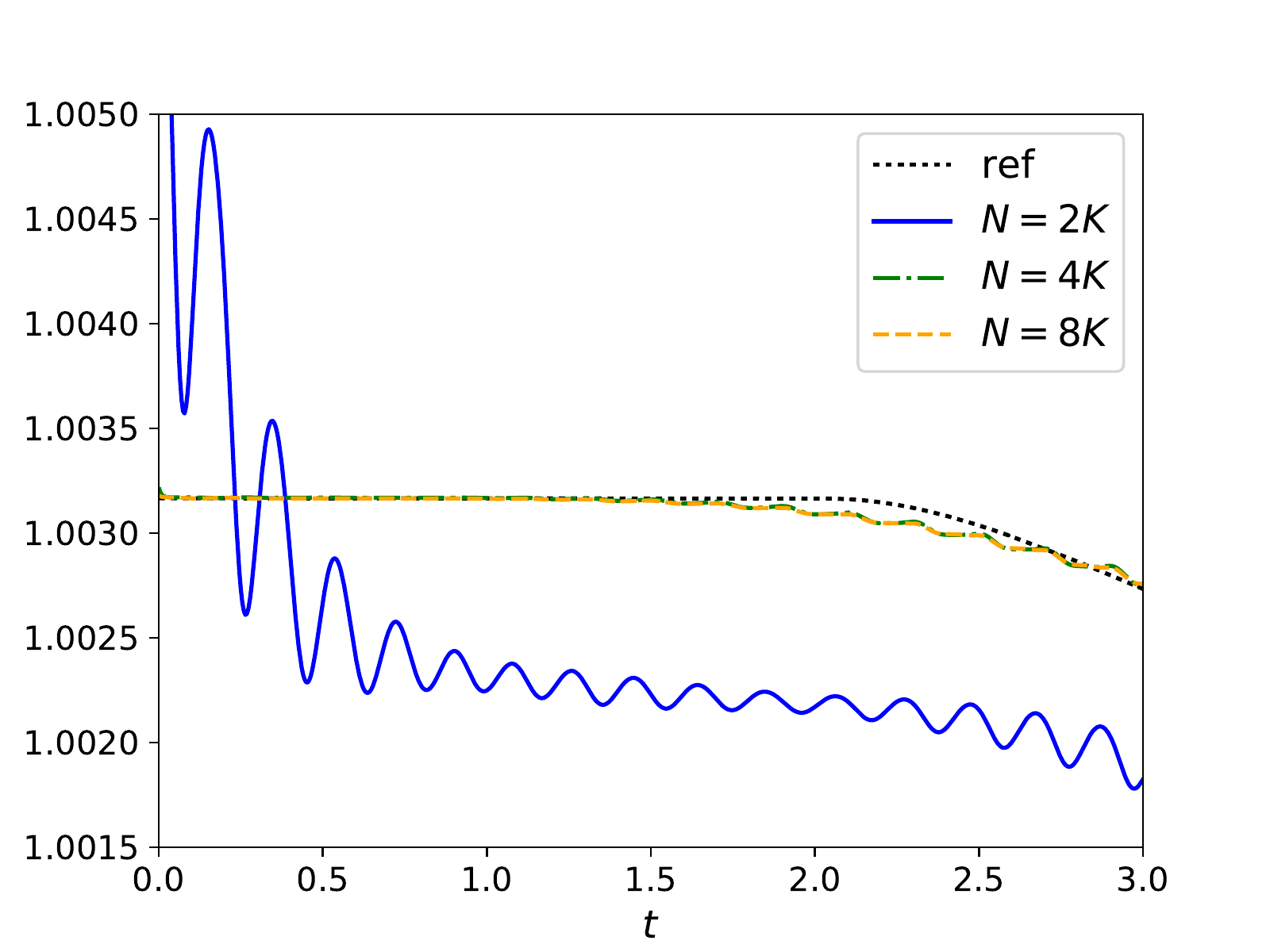}
    \caption{Energy over time.}
    \label{fig:burgers_t3_ad_energy_rand}
  \end{subfigure}%
  \caption{Numerical solution, mass, and energy for $I=5$, $K=3$, and $N=K,2K,4K$. Burgers' equation and scattered 
points.}
  \label{fig:burgers_ad_rand}
\end{figure}

Again, we start with a numerical investigation of conservation and $L^2$-stability. 
Note that our proof of $L^2$-stability in \S \ref{sub:stab} only addresses the linear advection equation. 
Still, we observe similar results for the nonlinear Burgers' equations. 
These results are illustrated in Figure \ref{fig:burgers_ad_eq} for a set of equidistant collocation points and in 
Figure \ref{fig:burgers_ad_rand} for a set of scattered collocation points, constructed by adding white uniform noise 
as described in \eqref{eq:noise}. 

Figures \ref{fig:burgers_t1_sol_eq} and \ref{fig:burgers_t1_sol_rand} show the solutions of the DGDLS 
method at time $t=1$. 
At this time, no discontinuity has developed and the solution is still smooth. 
For the case of equidistant collocation points in Figure \ref{fig:burgers_t1_sol_eq}, the DGDLS method for $N=2K,4K$ is 
observed to provide reasonable numerical solutions. 
However, the numerical solution for $N=K$ shows heavy oscillations. 
The same can be observed for numerical solutions for $N=K$ as well as $N=2K$ when scattered collocation points are used 
in Figure \ref{fig:burgers_t1_sol_rand}. 
Here, only the numerical solution for $N=4K$ can be considered as reasonable. 
A similar observation can be made for both kinds of collocation points at time $t=3$. 
In the case of scattered collocation points, the computation even broke down for $N=K$. 
Thus, Figure \ref{fig:burgers_t3_sol_rand} instead illustrates the results for $N=2K,4K,8K$. 
Note that all numerical results for $t=3$, also the ones for equidistant collocation points in Figure 
\ref{fig:burgers_t3_sol_eq}, show at least some minor oscillations. 
This is a common problem for high-order methods and might be overcome by post-processing (assuming a stable 
computation until the final time has been reached) or additional shock-capturing. 
Shock-capturing might be performed, for instance, by artificial viscosity methods 
\cite{persson2006sub,klockner2011viscous,glaubitz2019smooth,ranocha2018stability}, 
modal filtering 
\cite{vandeven1991family,hesthaven2008filtering,meister2013extended,glaubitz2018application}, 
finite volume subcells 
\cite{huerta2012simple,sonntag2014shock,dumbser2014posteriori,meister2016positivity}, 
or other methods 
\cite{glaubitz2019high,glaubitz2019shock}.
Shock capturing in DLS based high-order methods might be investigated in future works. 
Finally, when a sufficiently great number of collocation points is used, we again observe the mass of the 
numerical solutions to nearly remain constant, see Figures \ref{fig:burgers_t3_ad_mass_eq} and 
\ref{fig:burgers_t3_ad_mass_rand}, and the energy to nearly remain constant until the discontinuity occurs at time 
$t=2$ and to decrease afterwards, see Figures \ref{fig:burgers_t3_ad_energy_eq} and  
\ref{fig:burgers_t3_ad_energy_rand}.

\subsubsection{Accuracy and convergence}
\label{subsub:Burgers_accuracy}

We continue the above investigation of the DGDLS method by providing an additional error analysis. 
Tables \ref{tab:burgers_t1_equid} and \ref{tab:burgers_t1_rand} respectively list the errors for DGDLS method at time 
$t=1$ on equidistant and scattered points.
\begin{table}[htb!]
\parbox{.49\textwidth}{
    \centering 
    \begin{adjustbox}{width=0.49\textwidth}
    \begin{tabular}{c c c c c c} 
        \hline
        \multicolumn{2}{c}{}  
        & \multicolumn{4}{c}{$L^2$-errors}
        \\ \hline
	\multicolumn{2}{c}{} & DGSEM & \multicolumn{3}{c}{DGDLS}
        \\ \hline
        $K$ & $I$ & (GL points) & $N=K$ \ & $N=2K$ \ & $N=4K$ \ \, 
        \\ \hline
        1 &  5 & 1.3E-2
	       & 1.3E-2
               & 1.1E-2
               & 1.2E-2 \\
          & 10 & 3.8E-3
	       & 3.8E-3
               & 4.1E-3
               & 3.8E-3 \\
          & 20 & 1.1E-3
	       & 1.1E-3
               & 9.3E-4
               & 8.7E-4 \\
          & 40 & 2.8E-4
	       & 2.8E-4
               & 2.0E-4
               & 1.8E-4 \\
        \multicolumn{2}{c}{EOC:} & 1.8 & 1.8 & 1.5 & 1.7 \\ 
        \hline
        2 &  5 & 1.7E-3
	       & 1.7E-3
               & 3.4E-3
               & 3.0E-3 \\
          & 10 & 5.9E-4
	       & 5.9E-4
               & 3.4E-4
               & 3.5E-4 \\
          & 20 & 6.7E-5
	       & 6.7E-5
               & 5.0E-5
               & 4.5E-5 \\
          & 40 & 8.0E-6
	       & 8.0E-6
               & 6.2E-6
               & 5.5E-6 \\
        \multicolumn{2}{c}{EOC:} & 1.7 & 1.7 & 3.3 & 3.0 \\ 
        \hline
        3 &  5 & 1.0E-3
	       & 7.1E-2
               & 4.9E-4
               & 6.7E-4 \\
          & 10 & 9.4E-5
	       & 4.2E-1
               & 8.3E-5
               & 7.5E-5 \\
          & 20 & 5.8E-6
	       & NaN
               & 5.6E-6
               & 4.7E-6 \\
          & 40 & 3.6E-7
	       & NaN
               & 3.0E-7
               & 2.6E-7 \\
        \multicolumn{2}{c}{EOC:} & 3.4 & - & 2.6 & 3.1 \\ 
        \hline
        4 &  5 & 1.9E-4
	       & NaN
               & 4.3E-4
               & 3.1E-4 \\
          & 10 & 1.0E-5
	       & NaN
               & 1.3E-5
               & 1.1E-5 \\
          & 20 & 3.7E-7
	       & NaN
               & 2.2E-7
               & 2.2E-7 \\
          & 40 & 1.9E-8
	       & NaN
               & 1.8E-8
               & 1.4E-8 \\
        \multicolumn{2}{c}{EOC:} & 4.2 & - & 4.3 & 4.4 \\ 
        \hline
    \end{tabular} 
    \end{adjustbox}
    \caption{Burgers' equation at $t=1$ \& equidistant points.}
    \label{tab:burgers_t1_equid}
}
\hfill
\parbox{.49\textwidth}{
    \centering 
    \begin{adjustbox}{width=0.49\textwidth}
    \begin{tabular}{c c c c c c} 
        \hline
        \multicolumn{2}{c}{}  
        & \multicolumn{4}{c}{$L^2$-errors}
        \\ \hline
	\multicolumn{2}{c}{} & DGSEM & \multicolumn{3}{c}{DGDLS}
        \\ \hline
        $K$ & $I$ & (GL points) & $N=4K$ & $N=16K$ & $N=64K$
        \\ \hline
        1 &  5 & 1.3E-2
               & 1.2E-2  
               & 1.2E-2 
               & 1.2E-2 \\
          & 10 & 1.3E-2
               & 3.9E-3  
               & 3.7E-3 
               & 3.7E-3 \\
          & 20 & 1.1E-3
               & 8.6E-4  
               & 8.5E-4 
               & 8.4E-4 \\
          & 40 & 2.8E-4
               & 1.9E-4  
               & 1.8E-4 
               & 1.7E-4 \\
        \multicolumn{2}{c}{EOC:} & 0.8 & 1.7 & 1.7 & 1.7 \\ 
        \hline
        2 &  5 & 1.7E-3
               & 3.0E-3  
               & 2.9E-3 
               & 2.9E-3 \\
          & 10 & 5.9E-4
               & 1.3E-3  
               & 3.5E-4 
               & 3.5E-4 \\
          & 20 & 6.7E-5
               & 1.3E-4  
               & 4.9E-5 
               & 4.4E-5 \\
          & 40 & 8.0E-6
               & 3.6E-4  
               & 6.9E-5 
               & 9.9E-6 \\
        \multicolumn{2}{c}{EOC:} & 1.7 & 1.3 & 3.0 & 3.0 \\ 
        \hline
        3 &  5 & 1.0E-3
               & 9.1E-4  
               & 7.1E-4 
               & 6.9E-4 \\
          & 10 & 9.4E-5
               & 5.0E-4  
               & 9.9E-5 
               & 7.2E-5 \\
          & 20 & 5.8E-6
               & 4.4E-4  
               & 1.1E-5 
               & 4.9E-6 \\
          & 40 & 3.6E-7
               & 3.1E-4  
               & 7.7E-6 
               & 3.0E-6 \\
        \multicolumn{2}{c}{EOC:} & 3.4 & 0.5 & 2.8 & 3.2 \\ 
        \hline
        4 &  5 & 1.9E-4
               & 7.5E-4  
               & 3.5E-4 
               & 3.0E-4 \\
          & 10 & 1.0E-5
               & 9.1E-4  
               & 6.0E-5 
               & 1.3E-5 \\
          & 20 & 3.7E-7
               & 2.2E-3  
               & 1.7E-4 
               & 9.2E-6 \\
          & 40 & 1.9E-8
               & 1.7E-4  
               & 6.6E-5 
               & 8.5E-6 \\
        \multicolumn{2}{c}{EOC:} & 4.2 & 0.0 & 0.9 & 4.2 \\ 
        \hline
    \end{tabular}  
    \end{adjustbox}
    \caption{Burgers' equation at $t=1$ \& scattered points.}
    \label{tab:burgers_t1_rand}
}
\end{table}
We note that the DGDLS method for $N=K$ does not yield stable numerical solutions. 
Yet, we again observe that increasing the number of collocation points in the DGDLS method results in 
more accurate numerical solutions. 
At least for equidistant collocation points, the DGDLS method even yields more accurate numerical solutions than the 
DGSEM on Gauss--Lobatto points.

Further, in both cases, we are able to recover or even exceed the EOC of the DGSEM when a sufficiently great number of 
collocation points is used.

\subsection{Extension to systems: The wave equation}
\label{sub:system}

\begin{figure}[!htb]
  \centering
  \begin{subfigure}[b]{0.33\textwidth}
    \includegraphics[width=\textwidth]{%
      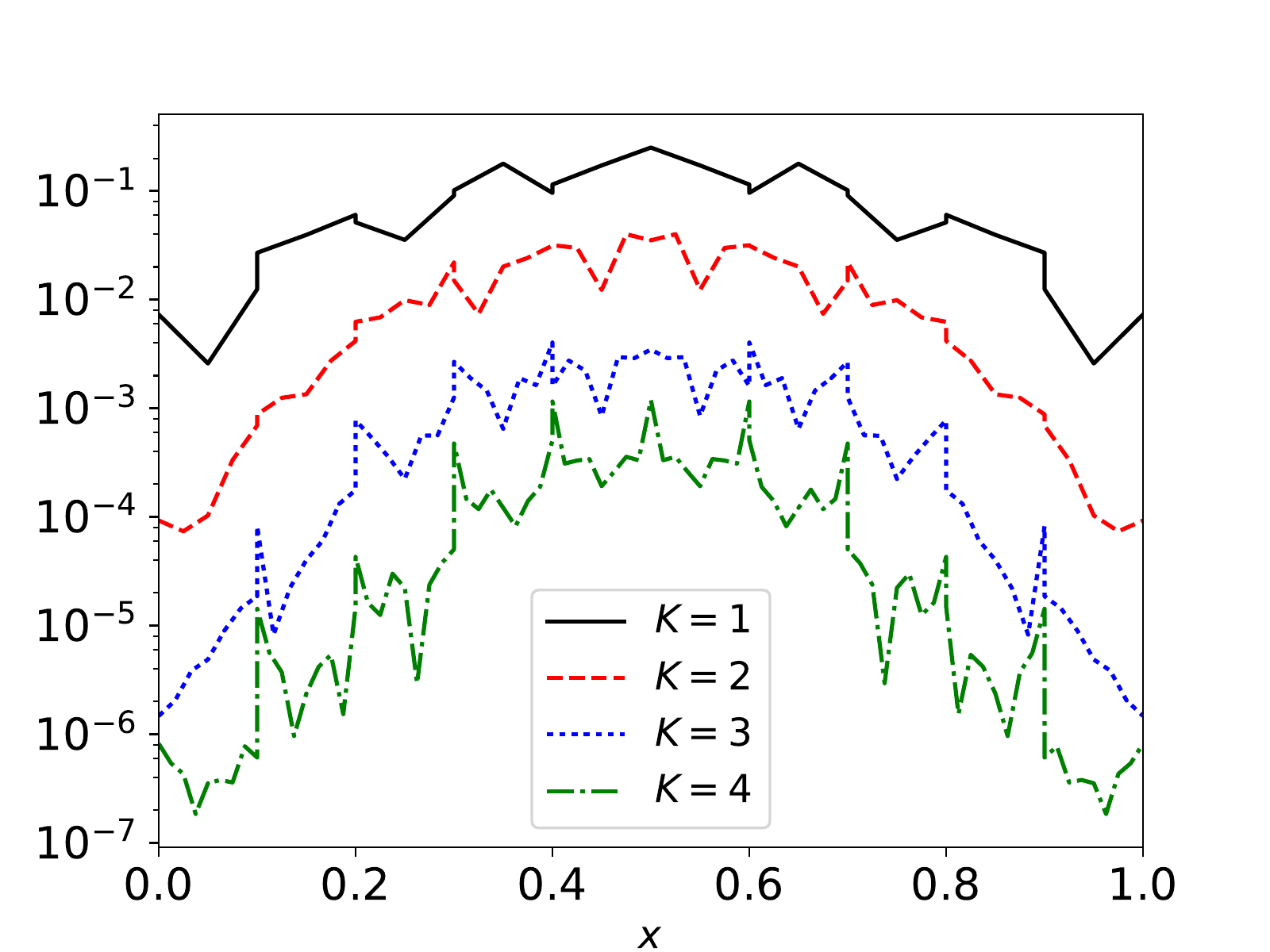}
    \caption{Errors at time $t=1$.}
    \label{fig:wave_equi_t1}
  \end{subfigure}%
  ~
  \begin{subfigure}[b]{0.33\textwidth}
    \includegraphics[width=\textwidth]{%
      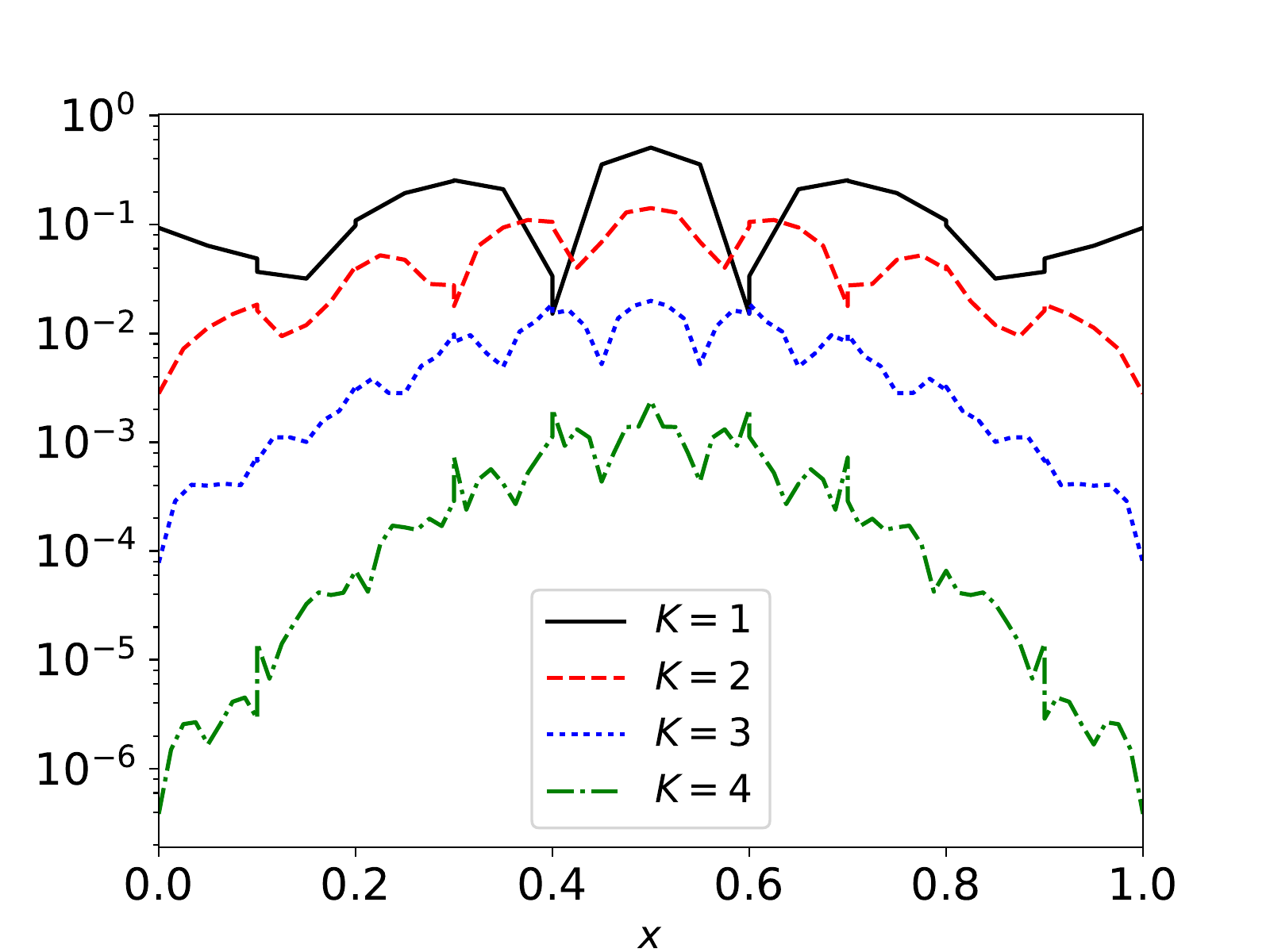}
    \caption{Errors at time $t=10$.}
    \label{fig:wave_equi_t10}
  \end{subfigure}%
  ~
  \begin{subfigure}[b]{0.33\textwidth}
    \includegraphics[width=\textwidth]{%
      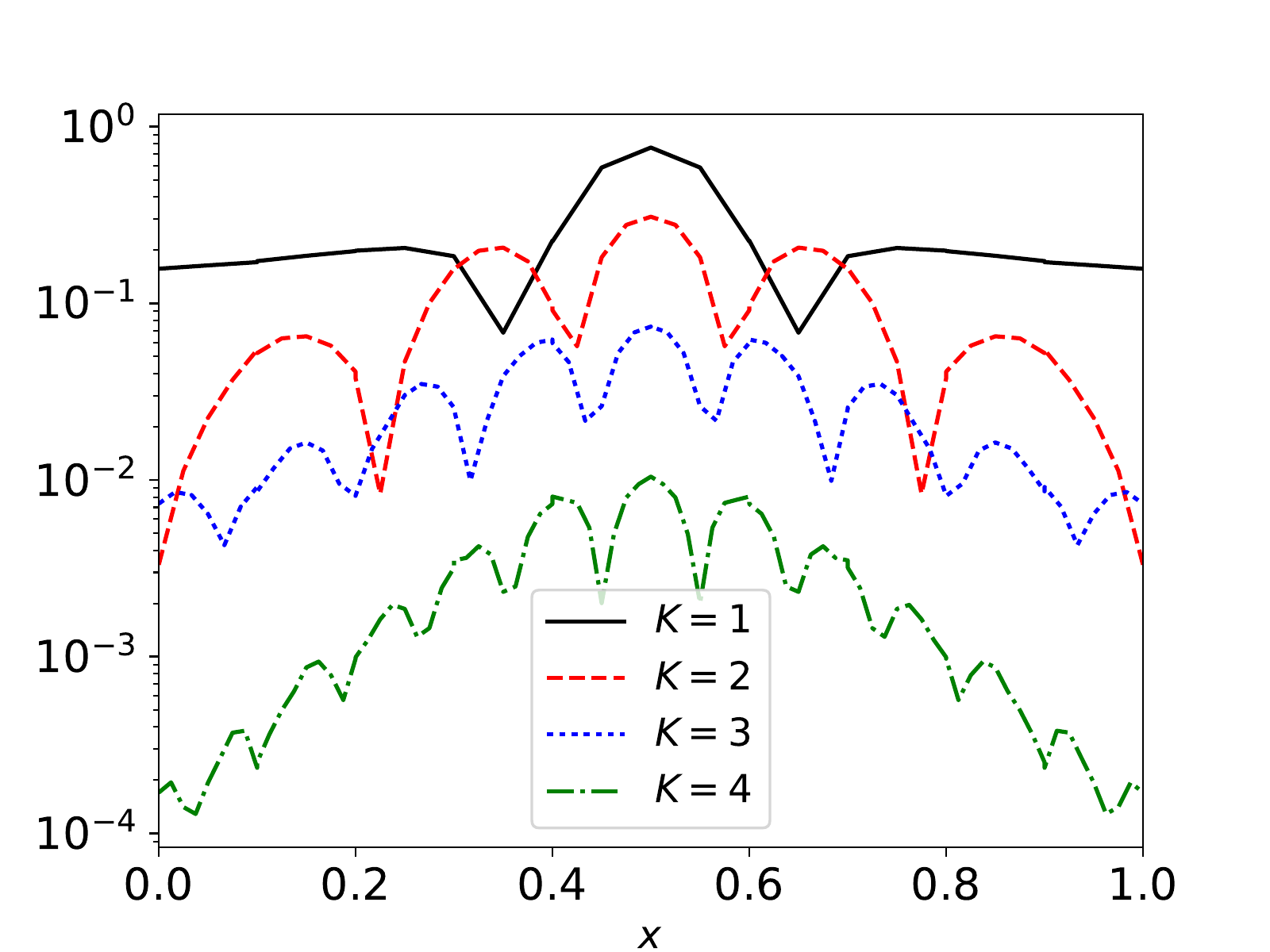}
    \caption{Errors at time $t=100$.}
    \label{fig:wave_equi_t100}
  \end{subfigure}%
  \caption{Pointwise errors over time for the wave equation. DGDLS method for $I=10$, $N=2K$, $K=1,2,3,4$, and 
equidistant collocation points.}
  \label{fig:wave_equi}
\end{figure}

\begin{figure}[!htb]
  \centering
  \begin{subfigure}[b]{0.33\textwidth}
    \includegraphics[width=\textwidth]{%
      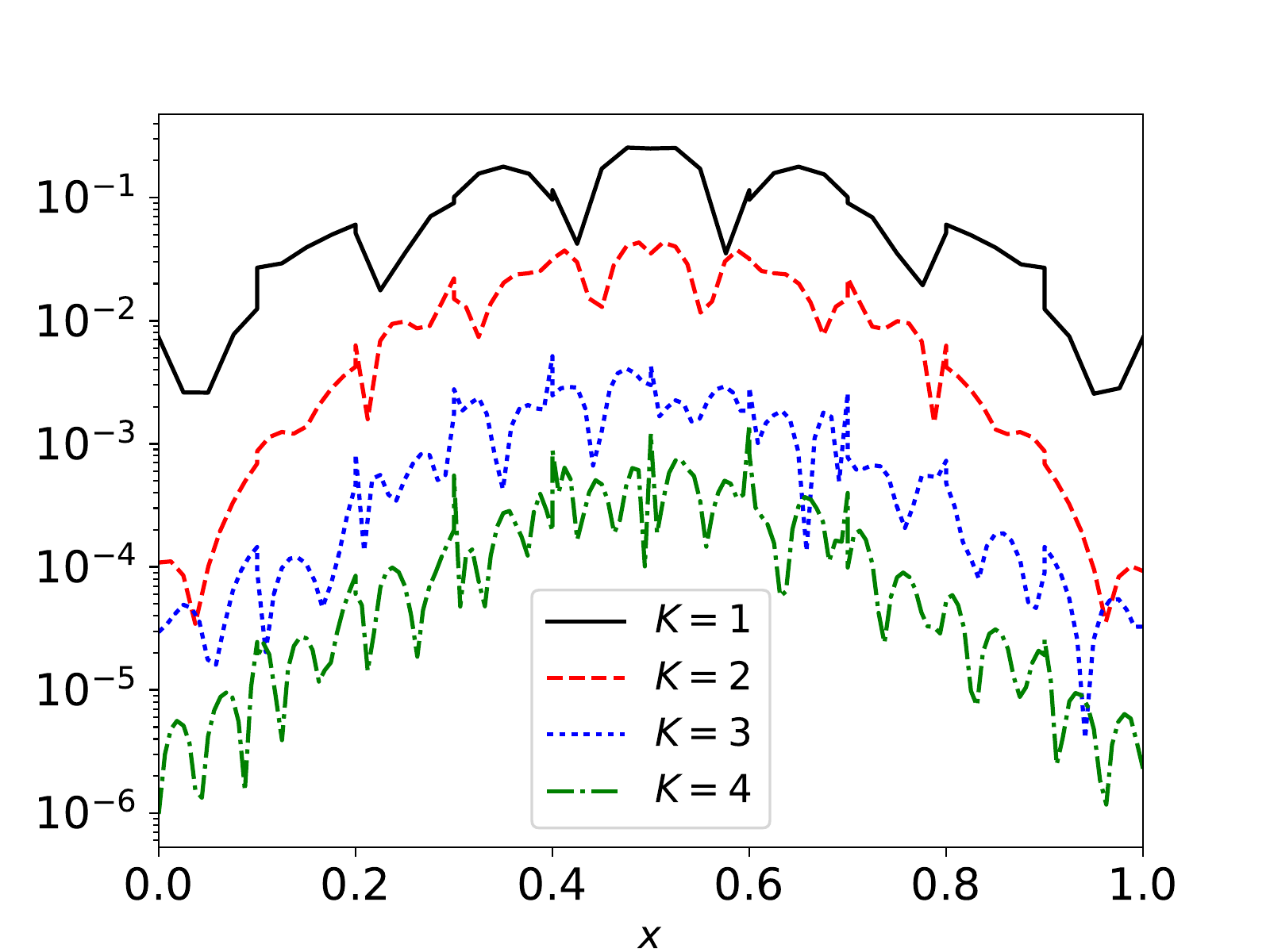}
    \caption{Errors at time $t=1$.}
    \label{fig:wave_rand_t1}
  \end{subfigure}%
  ~
  \begin{subfigure}[b]{0.33\textwidth}
    \includegraphics[width=\textwidth]{%
      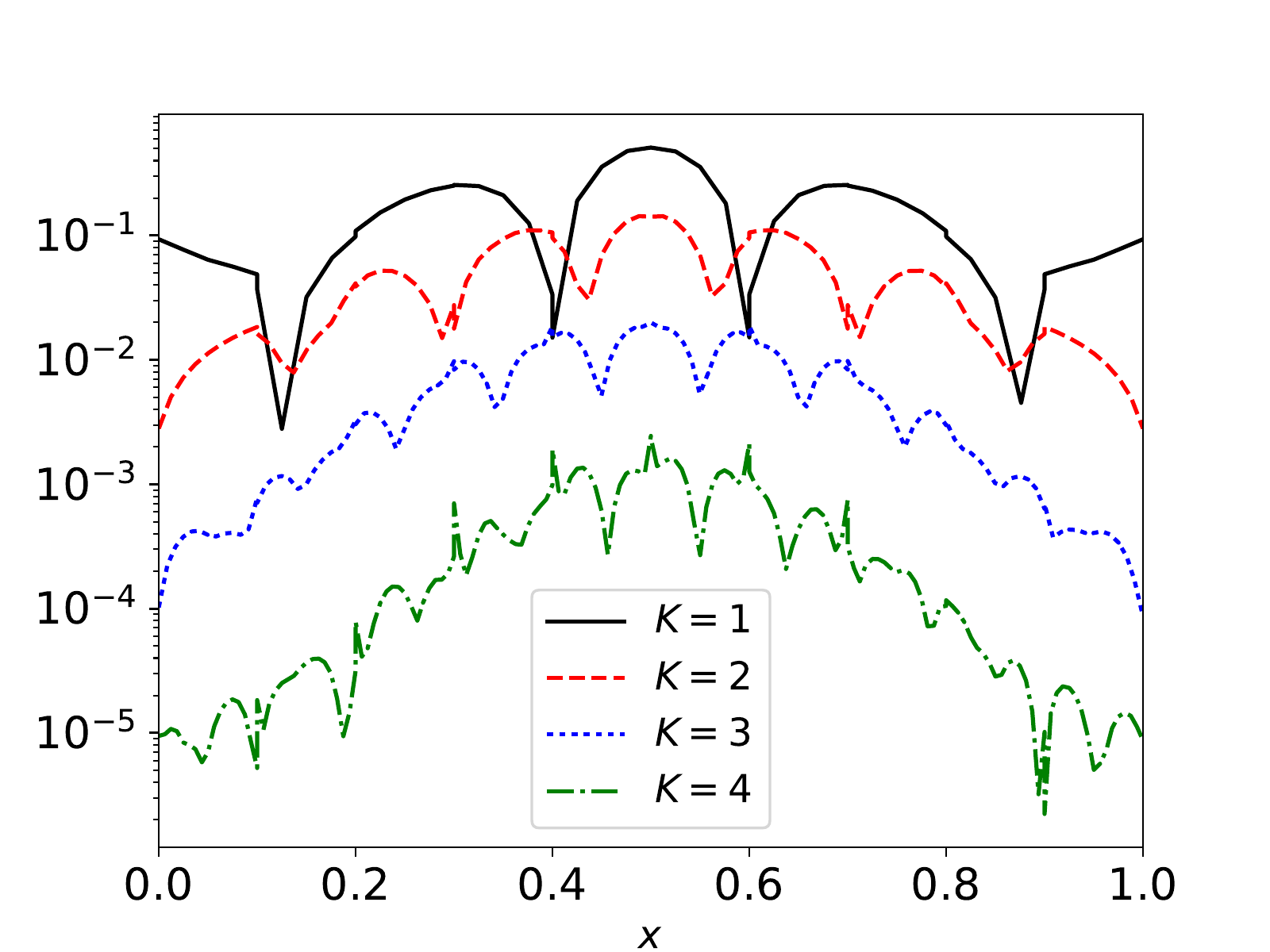}
    \caption{Errors at time $t=10$.}
    \label{fig:wave_rand_t10}
  \end{subfigure}%
  ~
  \begin{subfigure}[b]{0.33\textwidth}
    \includegraphics[width=\textwidth]{%
      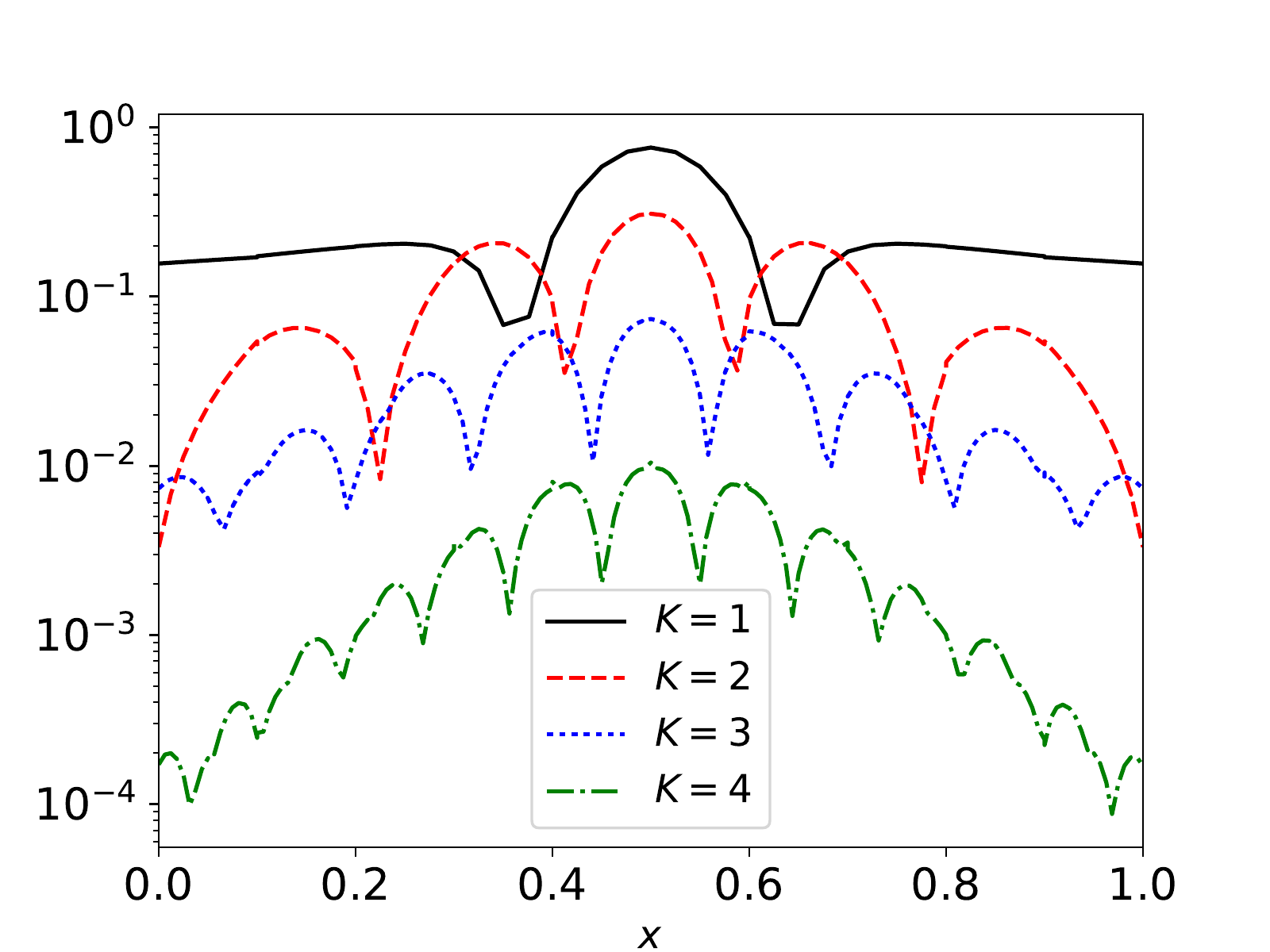}
    \caption{Errors at time $t=100$.}
    \label{fig:wave_rand_t100}
  \end{subfigure}%
  \caption{Pointwise errors over time for the wave equation. DGDLS method for $I=10$, $N=4K$, $K=1,2,3,4$, and 
scattered collocation points.}
  \label{fig:wave_rand}
\end{figure}

The extension of the DGDLS method to systems of conservation laws is the same as for most discretisations of the DG 
method, we simply apply the 
discretisation proposed in \S \ref{sec:proposed} to every component of the system separately. 
As a representative, we consider the second order scalar wave equation in one dimension 
\begin{equation}\label{eq:wave-eq}
  u_{tt} - c^2 u_{xx} = 0
\end{equation} 
on $\Omega=[0,1]$ with periodic boundary conditions. 
The wave equation \eqref{eq:wave-eq} can be rewritten as a first-order system of conservation laws 
\begin{equation}
\begin{aligned}
  u_t + c v_x & = 0, \\ 
  v_t + c u_x & = 0,
\end{aligned} 
\end{equation}
which is sometimes referred to as the one dimensional acoustic problem \cite{gelb2008discrete}. 
Given initial conditions ${u(0,x) = u_0(x)}$ and ${v(0,x) = v_0(x)}$, the solution is given by 
\begin{equation}
\begin{aligned}
  u(t,x) & = \frac{1}{2} \left[ u_0(x-ct) + u_0(x+ct) \right] + \frac{1}{2} \left[ v_0(x-ct) - v_0(x+ct) \right], \\
  v(t,x) & = \frac{1}{2} \left[ u_0(x-ct) - u_0(x+ct) \right] + \frac{1}{2} \left[ v_0(x-ct) + v_0(x+ct) \right]. 
\end{aligned}
\end{equation}
In the subsequent numerical tests we choose $c=1$ and consider the initial conditions ${u_0(x) = e^{-20(2x-1)^2}}$ and 
${v_0(x) = 0}$. 
For the numerical flux, we have used an upwind flux 
\begin{equation}
  \fnum = \frac{1}{2}
  \begin{pmatrix}
    (v^- + v^+) - (u^+ - u^-) \\ 
    (u^- + u^+) - (v^+ - v^-)
  \end{pmatrix}.
\end{equation}
Figures \ref{fig:wave_equi} and \ref{fig:wave_rand} illustrate the pointwise errors of the DGDLS method for a 
longtime simulation.
We have used $I=20$ elements and an increasing polynomial degree $K=1,2,3,4$ in all computations. 
Further, for equidistant points in Figure \ref{fig:wave_equi} we have used $N=2K$ collocation points and for the case 
of scattered points in Figure \ref{fig:wave_rand} we have used $N=4K$ collocation points. 
The figures show the pointwise errors 
\begin{equation}
  E(t,x) := | u(t,x) - u_{\mathrm{num}}(t,x) | + | v(t,x) - v_{\mathrm{num}}(t,x) |
\end{equation}
for times $t=1,10,100$. 
We note that the DGDLS method, regardless of whether equidistant or scattered collocation points are used, yields 
accurate results even for longtime simulations, especially when higher polynomial degrees are used. 
Future work might include a more detailed investigation of long-time error behaviour as, for instance, performed in 
\cite{kopriva2017error,offner2019error}.

\subsection{Extension to multiple dimensions: Linear advection equation with variable coefficients}
\label{sub:2d}
\begin{figure}[!htb]
  \centering
  \begin{subfigure}[b]{0.32\textwidth}
    \includegraphics[width=\textwidth]{%
      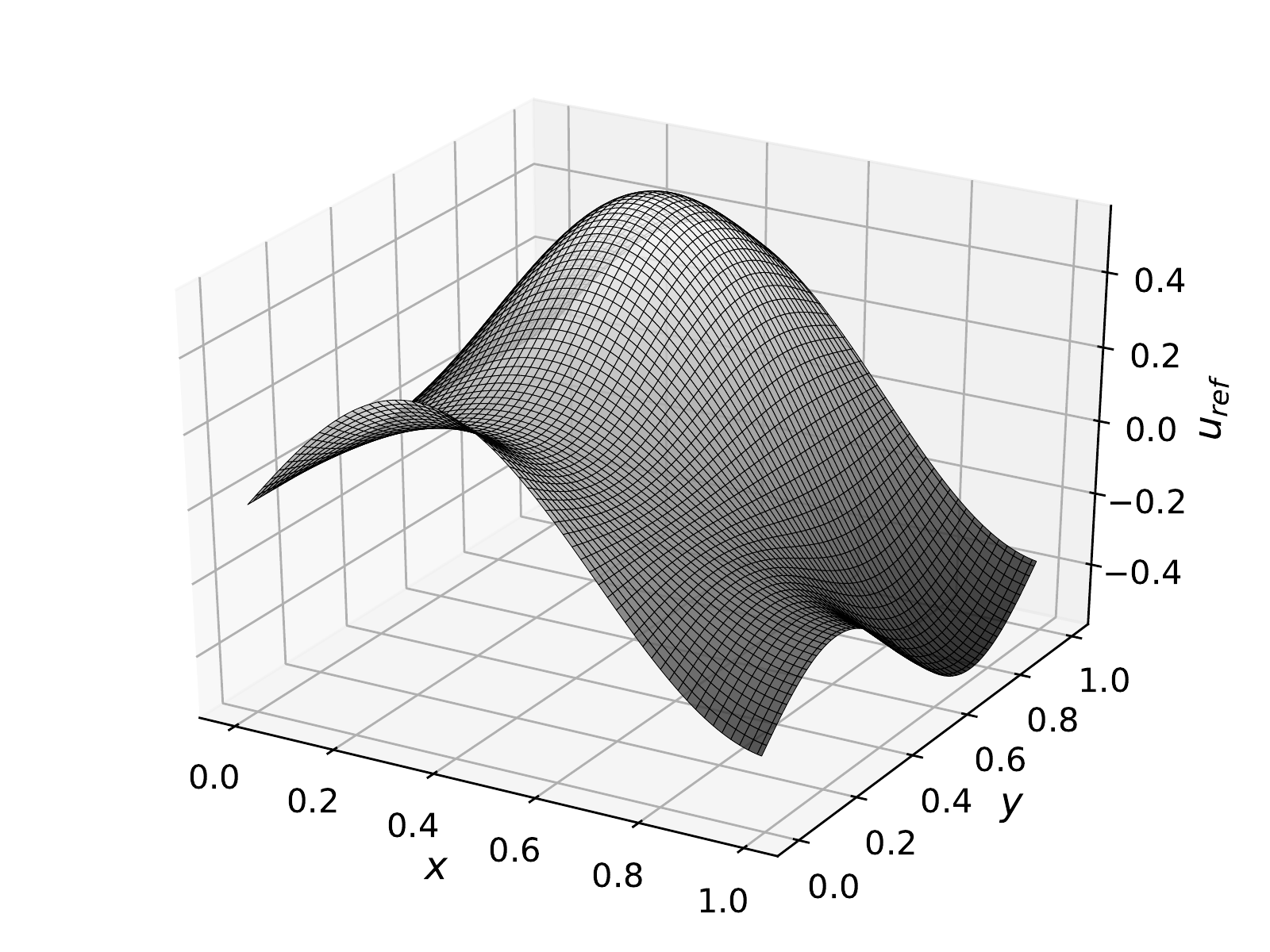}
    \caption{Reference solution.}
    \label{fig:2d_surface_ref}
  \end{subfigure}%
  ~
  \begin{subfigure}[b]{0.32\textwidth}
    \includegraphics[width=\textwidth]{%
      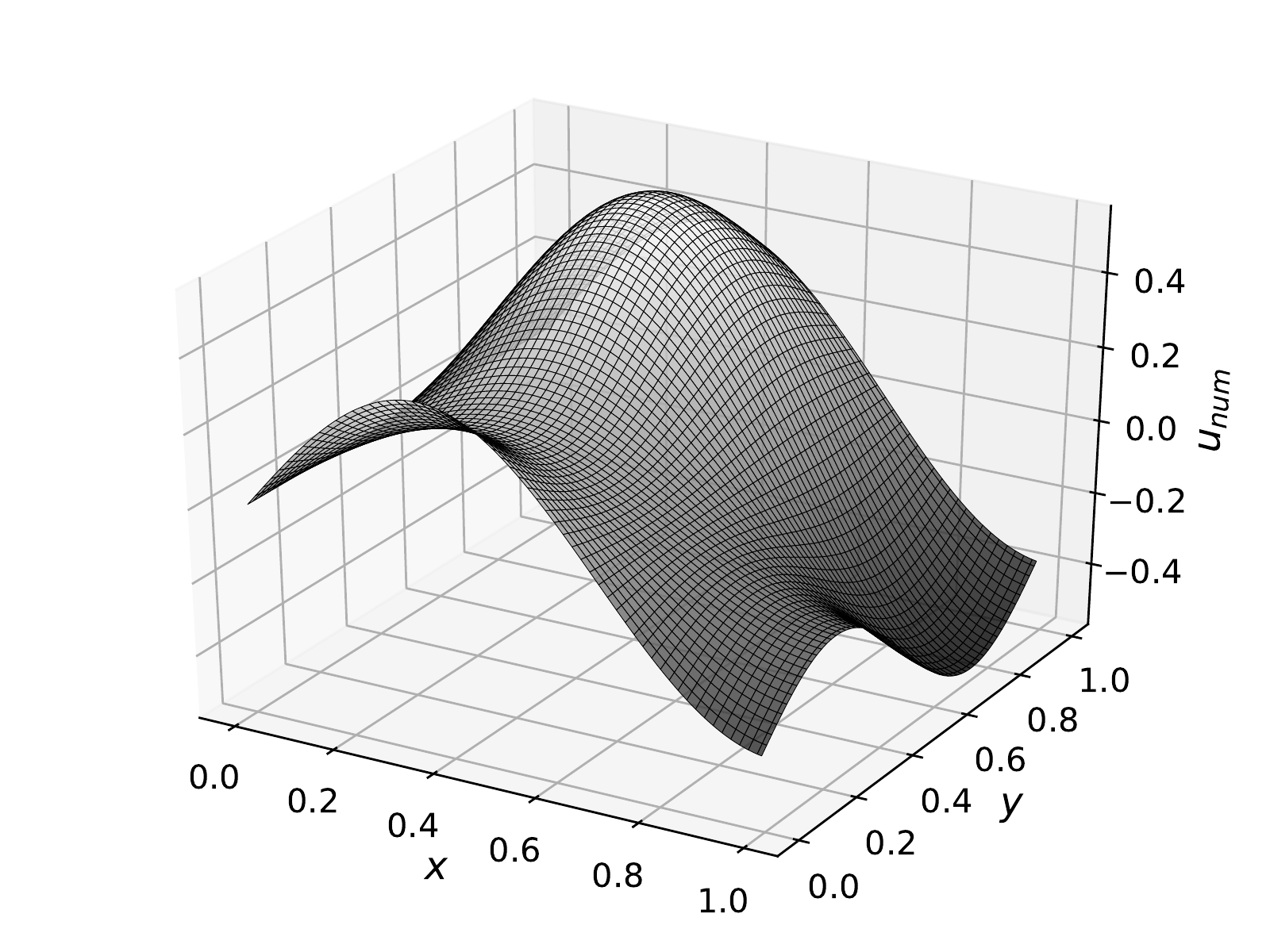}
    \caption{DGDLS, equid. points, $N=2K$.}
    \label{fig:2d_surface_DGDLS_equi_N=2K}
  \end{subfigure}%
  ~
  \begin{subfigure}[b]{0.32\textwidth}
    \includegraphics[width=\textwidth]{%
      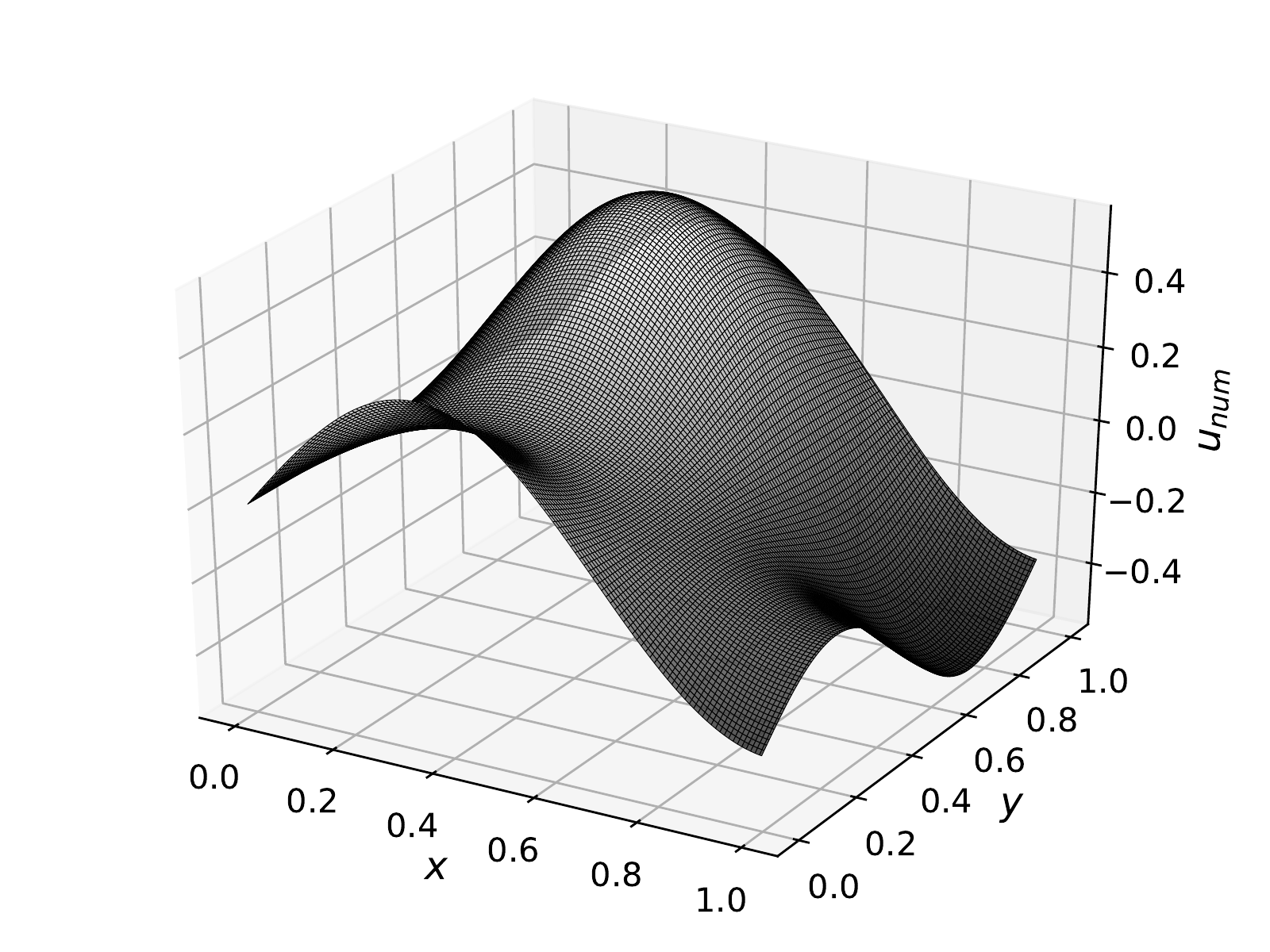}
    \caption{DGDLS, equid. points, $N=4K$.}
    \label{fig:2d_surface_DGDLS_equi_N=4K}
  \end{subfigure}%
  \\ 
  \begin{subfigure}[b]{0.32\textwidth}
    \includegraphics[width=\textwidth]{%
      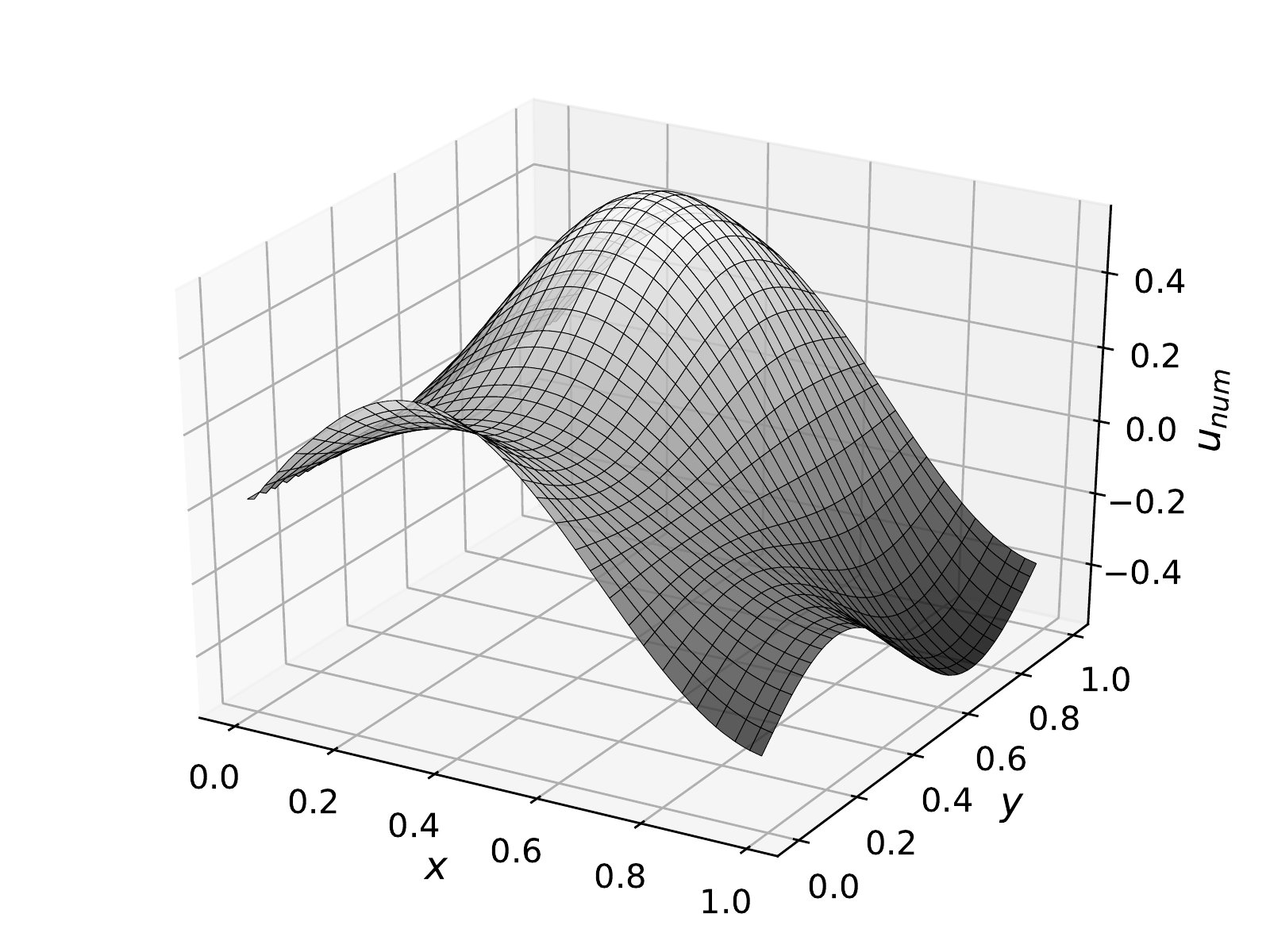}
    \caption{DGSEM, Gauss--Lobatto points.}
    \label{fig:2d_surface_DG_GL}
  \end{subfigure}%
  ~
  \begin{subfigure}[b]{0.32\textwidth}
    \includegraphics[width=\textwidth]{%
      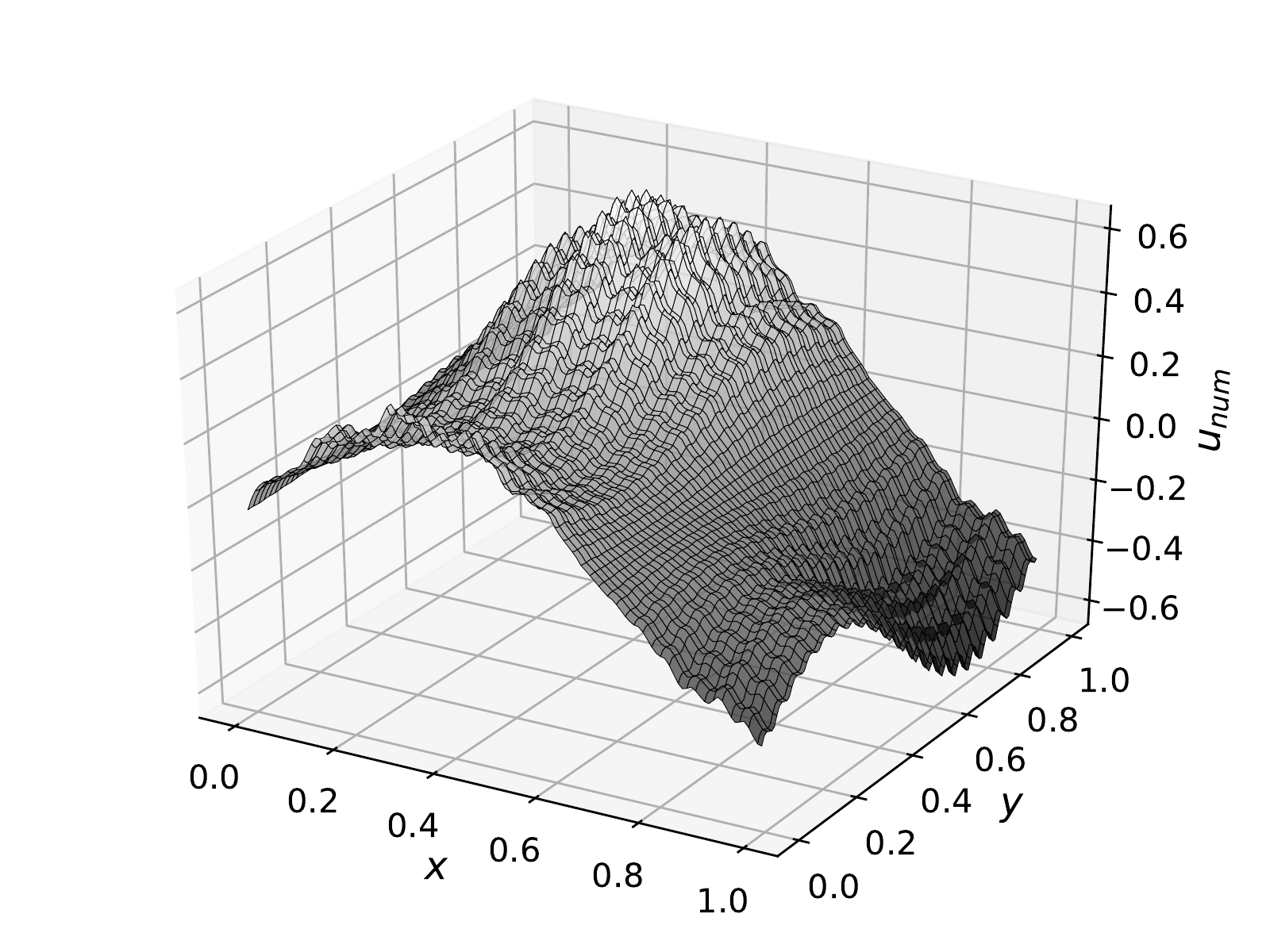}
    \caption{DGDLS, scat. points, $N=2K$.}
    \label{fig:2d_surface_DGDLS_scat_N=2K}
  \end{subfigure}%
  ~
  \begin{subfigure}[b]{0.32\textwidth}
    \includegraphics[width=\textwidth]{%
      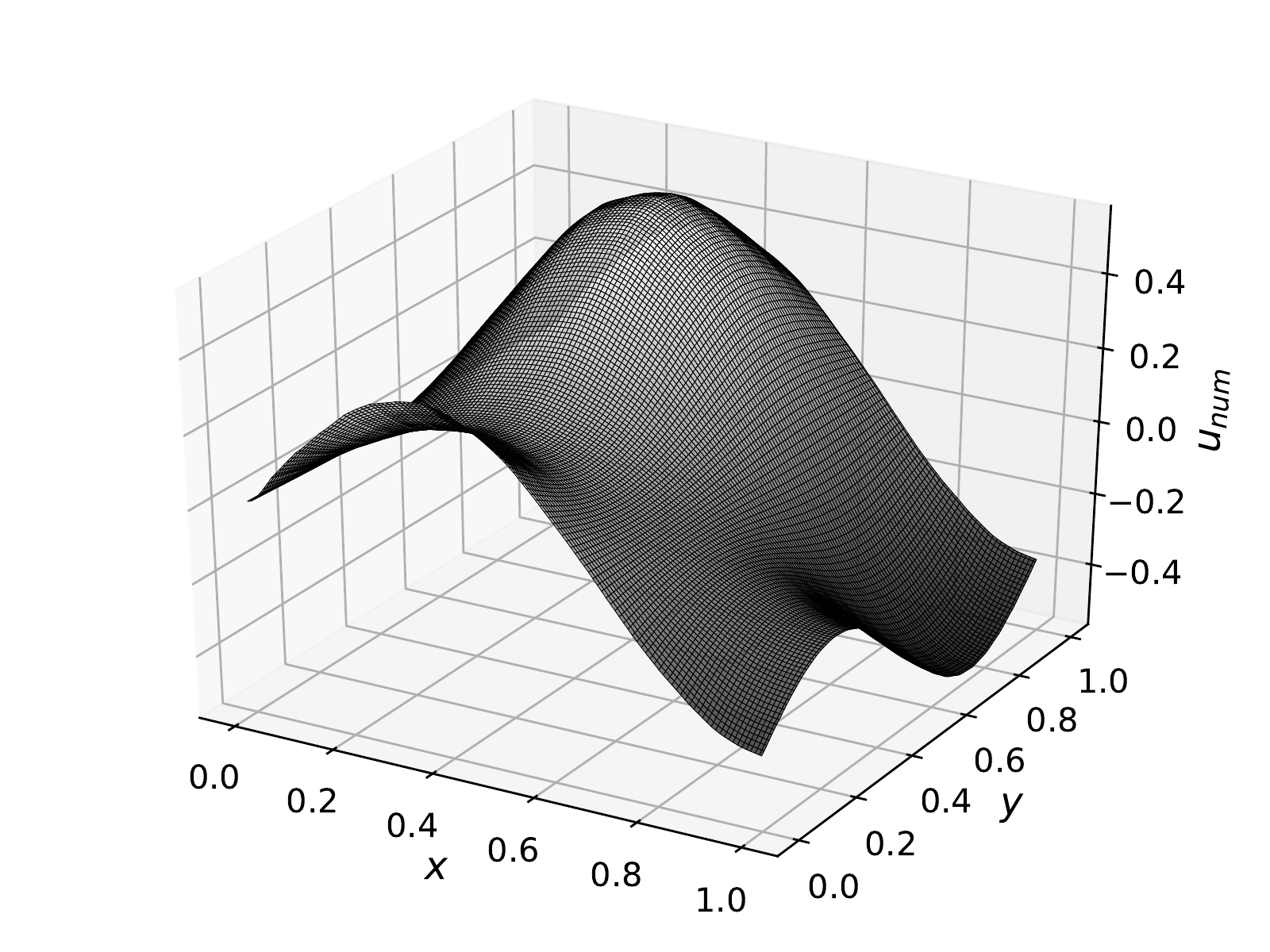}
    \caption{DGDLS, scat. points, $N=4K$.}
    \label{fig:2d_surface_DGDLS_scat_N=4K}
  \end{subfigure}%
  \caption{Surface plots of the solution at time $t=1$. 
  For the DGSEM and DGDLS method, a polynomial degree of $K=3$ and $I=20$ rectangular elements have been used in each 
direction.}
  \label{fig:2d_surface}
\end{figure}

\begin{figure}[!htb]
  \centering
  \begin{subfigure}[b]{0.32\textwidth}
    \includegraphics[width=\textwidth]{%
      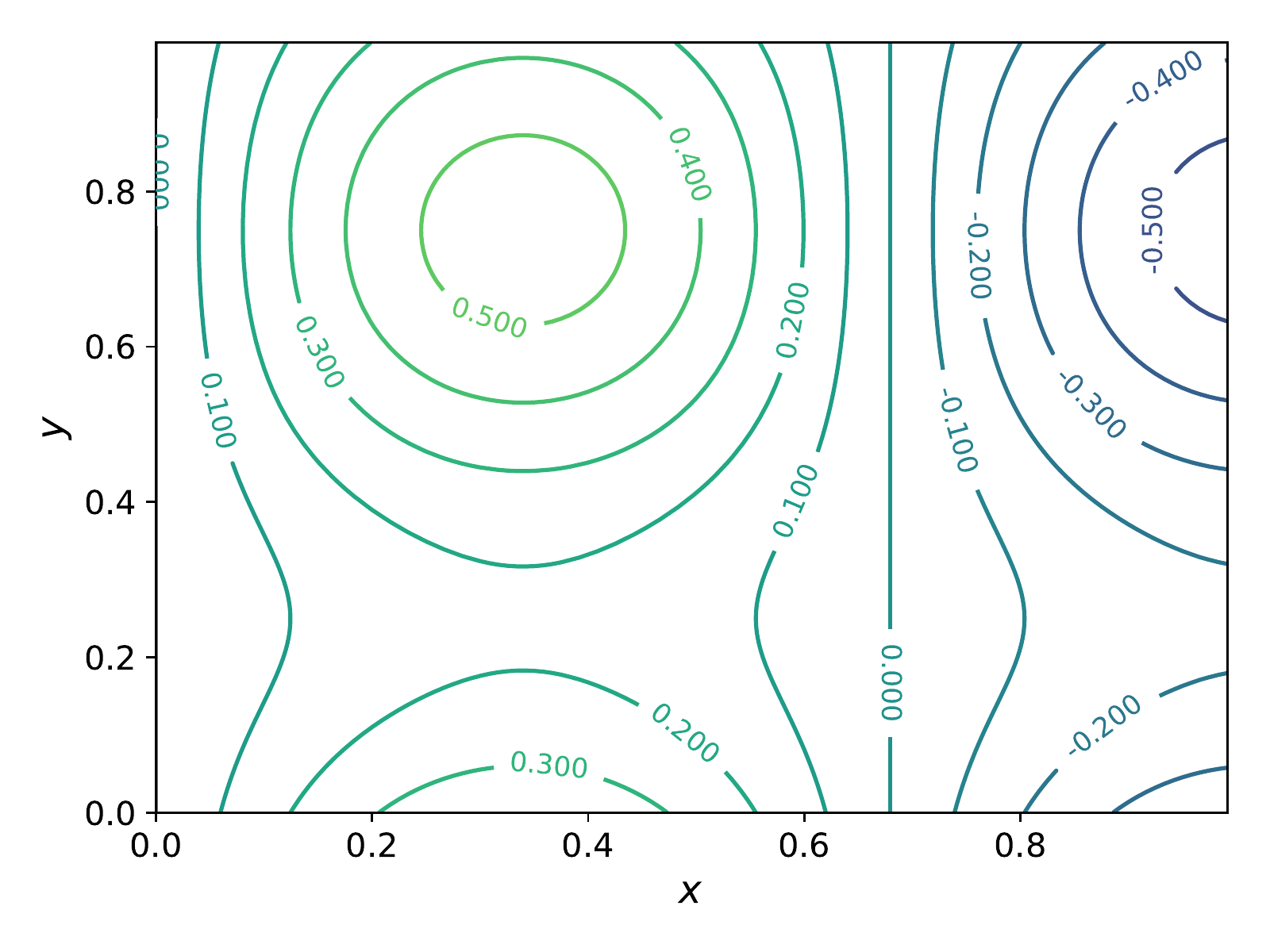}
    \caption{Reference solution.}
    \label{fig:2d_contour_ref}
  \end{subfigure}%
  ~
  \begin{subfigure}[b]{0.32\textwidth}
    \includegraphics[width=\textwidth]{%
      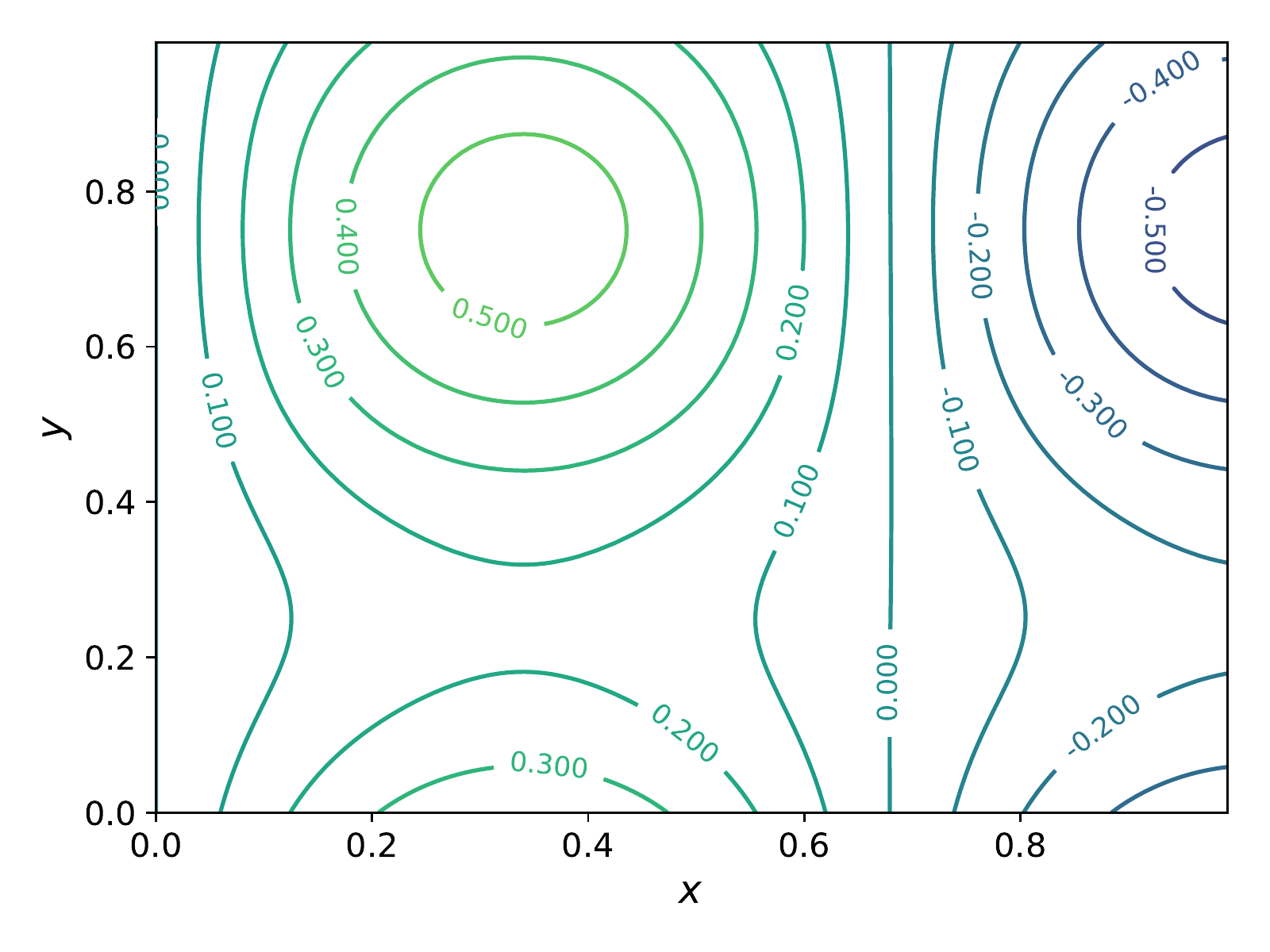}
    \caption{DGDLS, equid. points, $N=2K$.}
    \label{fig:2d_contour_DGDLS_equi_N=2K}
  \end{subfigure}%
  ~
  \begin{subfigure}[b]{0.32\textwidth}
    \includegraphics[width=\textwidth]{%
      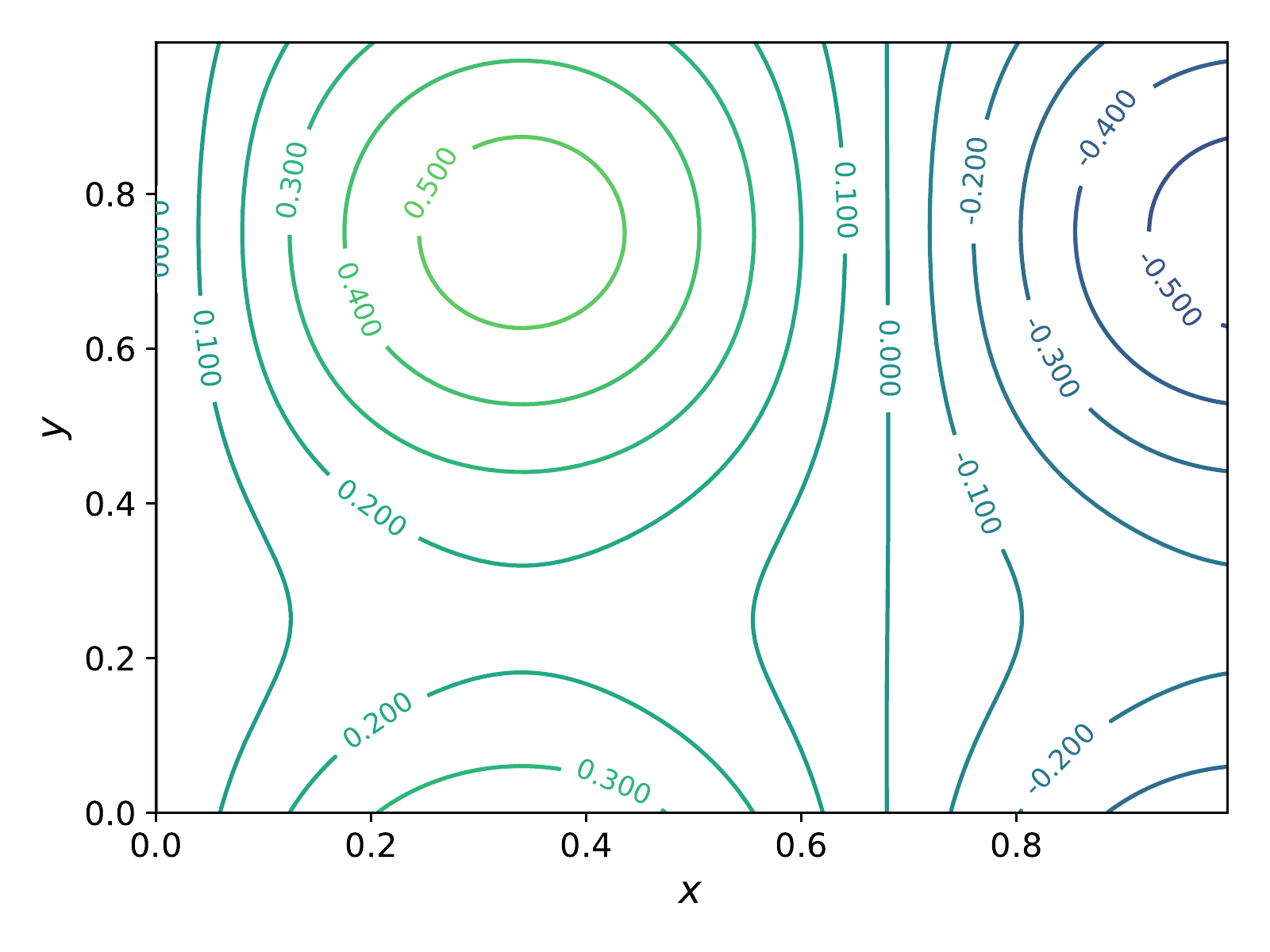}
    \caption{DGDLS, equid. points, $N=4K$.}
    \label{fig:2d_contour_DGDLS_equi_N=4K}
  \end{subfigure}%
  \\ 
  \begin{subfigure}[b]{0.32\textwidth}
    \includegraphics[width=\textwidth]{%
      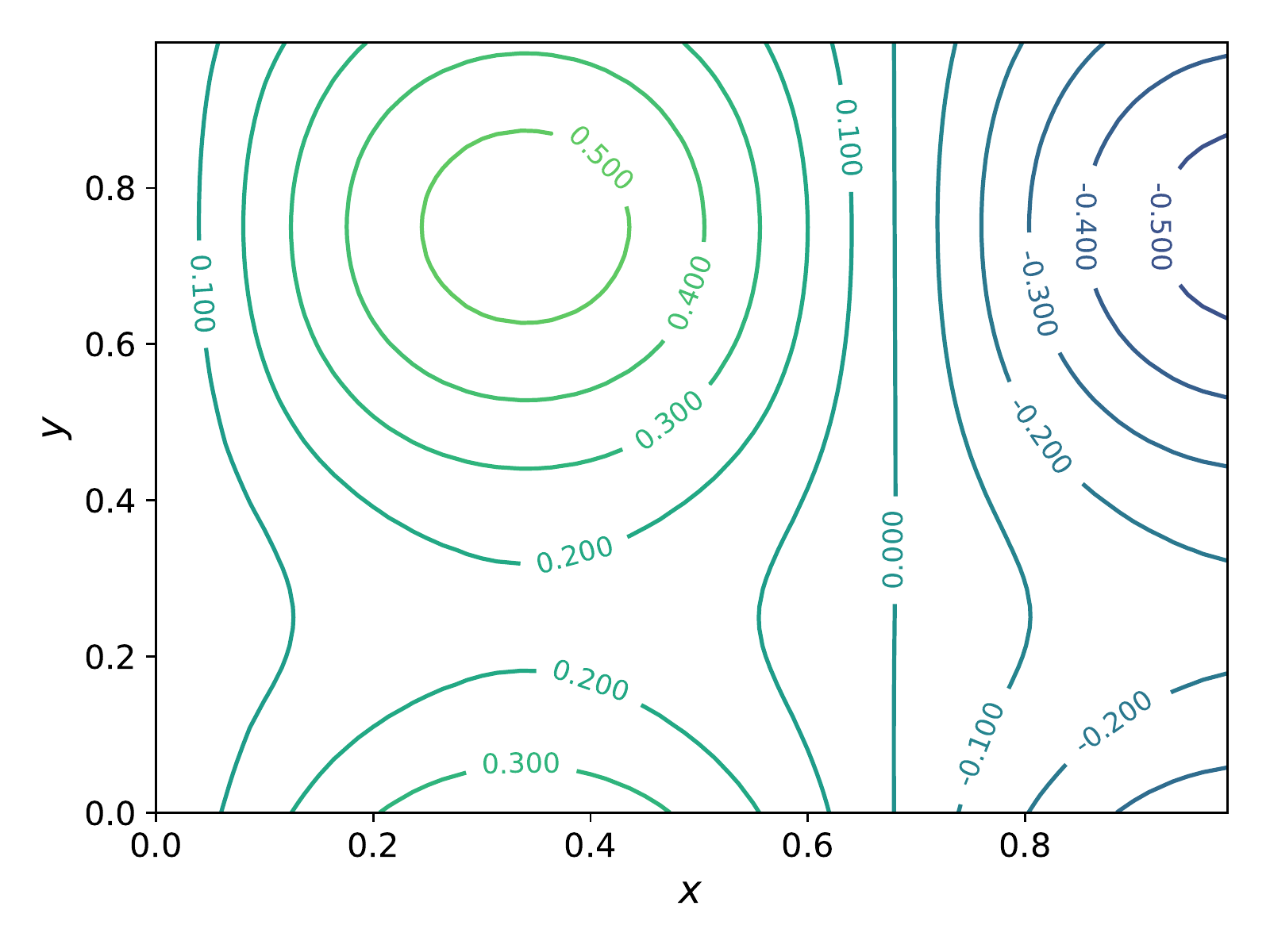}
    \caption{DGSEM, Gauss--Lobatto points.}
    \label{fig:2d_contour_DG_GL}
  \end{subfigure}%
  ~
  \begin{subfigure}[b]{0.32\textwidth}
    \includegraphics[width=\textwidth]{%
      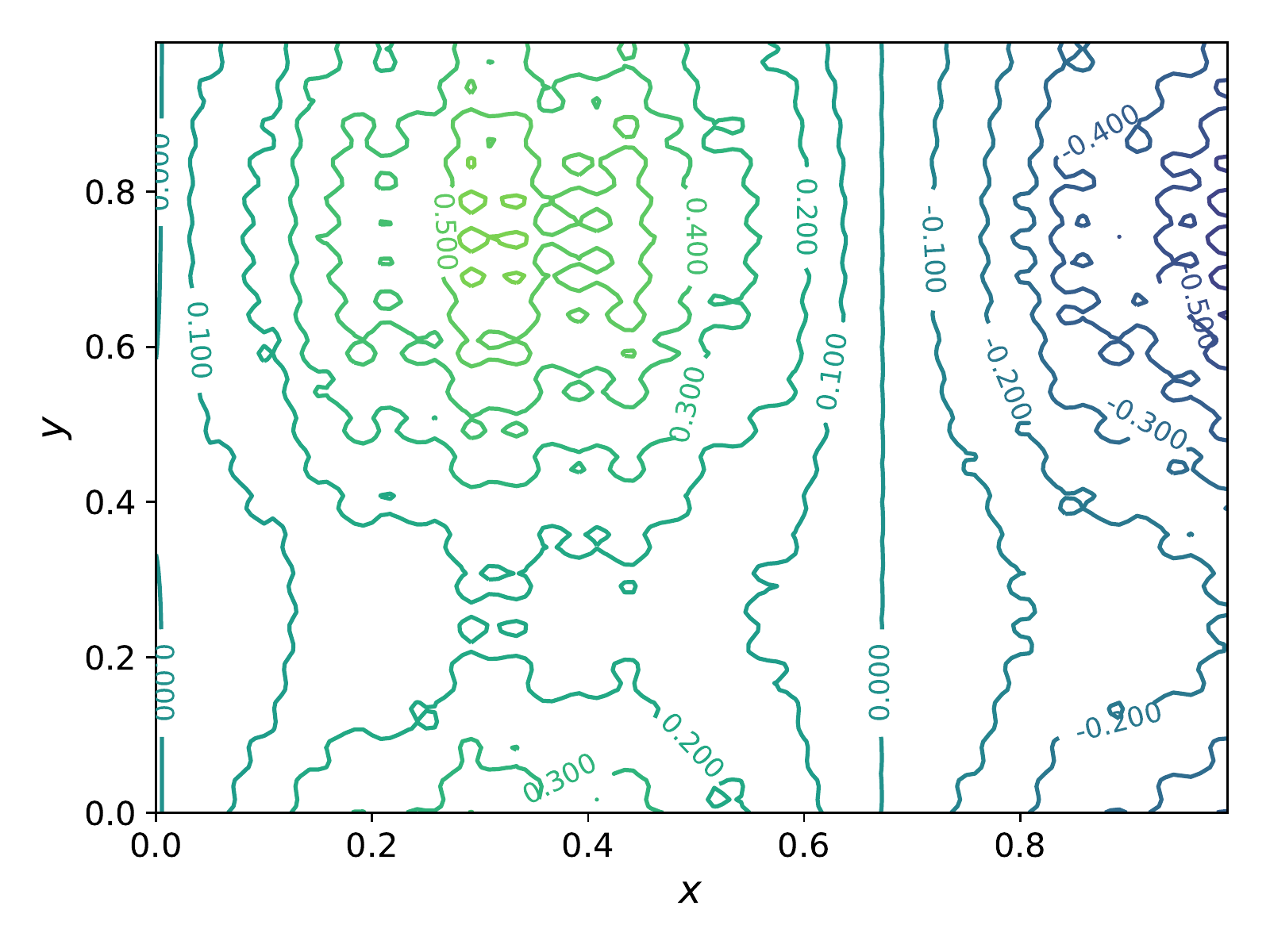}
    \caption{DGDLS, scat. points, $N=2K$.}
    \label{fig:2d_contour_DGDLS_scat_N=2K}
  \end{subfigure}%
  ~
  \begin{subfigure}[b]{0.32\textwidth}
    \includegraphics[width=\textwidth]{%
      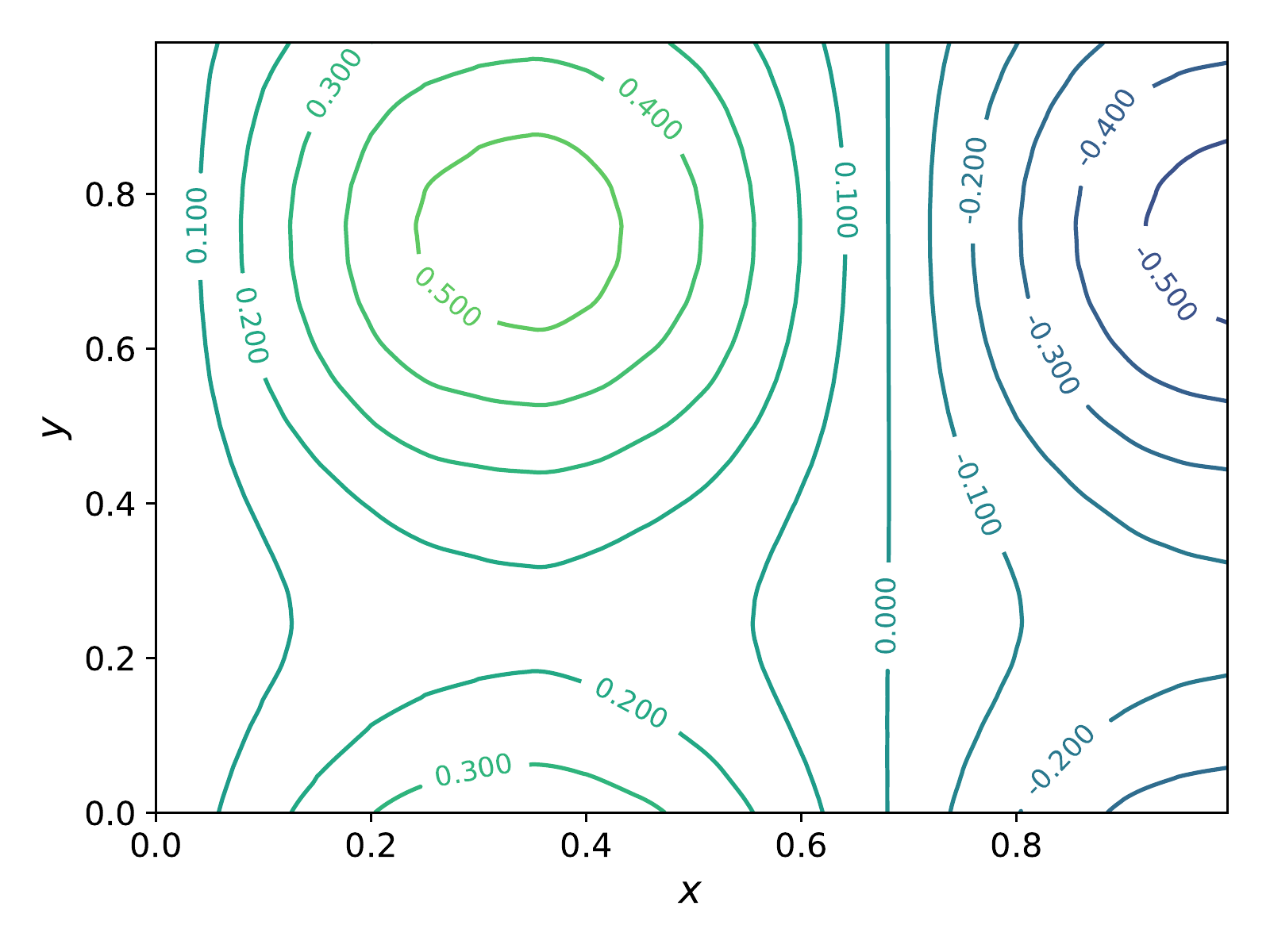}
    \caption{DGDLS, scat. points, $N=4K$.}
    \label{fig:2d_contour_DGDLS_scat_N=4K}
  \end{subfigure}%
  \caption{Contour plots of the solution at time $t=1$. 
  For the DGSEM and DGDLS method, a polynomial degree of $K=3$ and $I=20$ rectangular elements have been used in each 
direction.}
  \label{fig:2d_contour}
\end{figure}

In this subsection, we address the extension of the proposed DGDLS method to multiple dimensions by extending the ideas 
from the one dimensional case using a tensor product grid; see \cite[Chapter 7.1.6]{trangenstein2009numerical}.
Let us consider the linear advection equation in two dimensions 
\begin{equation}\label{eq:2d-problem}
  u_t + (au)_x + (bu)_y = 0
\end{equation}
with variable coefficients $a=a(x,y)$ and $b=b(x,y)$ on ${\Omega=[0,1]^2 \subset \R^2}$.  
Here, we choose ${a(x) = x}$ and $b(x) = 1$ together with initial condition 
\begin{equation}
  u_0(x,y) = \sin(4\pi x) \left( 1 - \frac{1}{2} \sin(2\pi y) \right)
\end{equation}
and boundary conditions 
\begin{align}
  u(t,x,0) & = u(t,x,1), \label{eq:BC1} \\
  u(t,0,y) & = 0. \label{eq:BC2}
\end{align}
Note that \eqref{eq:BC1} corresponds to periodic boundary conditions at the upper and lower boundary of $\Omega$, while 
\eqref{eq:BC2} corresponds to a physical inflow boundary condition at the left boundary of $\Omega$. 
The solution of the corresponding Cauchy problem can be calculated by the method of characteristics  
(see \cite[Chapter 3]{bressan2000hyperbolic}) and is given by 
\begin{equation}\label{eq:ref-sol} 
\begin{aligned}
  u(t,x,y) 
    & = \exp(-t) u_0(x\exp(-t),y) \\ 
    & = \exp(-t) \sin( 4 \pi x \exp(-t) ) \left( 1 - \frac{1}{2} \sin(2\pi y) \right). 
\end{aligned}
\end{equation}
Figure \ref{fig:2d_surface} and \ref{fig:2d_contour} respectively illustrate different surface and contour plots.   
Both figures demonstrate that the DGDLS method using $I=20$ rectangular elements, a polynomial degree of $K=3$ and 
$N=2K$ equidistant collocation points in each direction yields numerical solutions which are in good agreement with the 
reference solution \eqref{eq:ref-sol}. 
When using $N=2K$ scattered collocation points in each direction, the DGDLS method produces numerical solutions with 
slight oscillations; see Figure \ref{fig:2d_surface_DGDLS_scat_N=2K} for the surface plot and Figure 
\ref{fig:2d_contour_DGDLS_scat_N=2K} for the corresponding contour plot. 
Yet, when going over to a greater number of $N=4K$ scattered collocation points in each direction, the DGDLS method 
again yields numerical solutions which are in good agreement with the reference solution. 
 
\section{Summary}
\label{sec:summary} 

In this work, we have proposed and investigated stable collocation-type discretisations of the DG method on equidistant 
and scattered collocation points. 
We have done so by utilising DLS approximations instead of usual polynomial interpolation and 
LS-QRs, providing stable high-order numerical integration even on equidistant and scattered points. 
In \S \ref{sec:cons_stab}, we have proved conservation and linear $L^2$-stability of the proposed DGDLS method.  
In several numerical tests we have observed that the DGDLS method on equidistant points is able to recover --- 
sometimes even to exceed --- the accuracy and EOC of the usual DGSEM on Gauss--Lobatto points. 
Finally, the extension to the nonlinear viscous Burgers' equation, systems of conservation laws, and long time 
simulations, and a variable coefficient problem in two dimensions using a tensor product approach have been 
demonstrated. 
Future work will address more general entropy stability and the extension to unstructured meshes, such as  
triangles and nonconvex polygons. 

\section*{Acknowledgements}
The first author would like to thank the Max Planck Institute for Mathematics (MPIM) Bonn for 
wonderful working conditions.
Further, the first author was supported by the German Research Foundation (DFG, Deutsche 
Forschungsgemeinschaft) under Grant SO 363/15-1.

The second author was supported by SNF project (Number 175784) ``Solving advection dominated
problems with high order schemes with polygonal meshes: application to compressible
and incompressible flow problems''.

\bibliographystyle{alphaabbr}
\bibliography{literature}

\newcommand{\etalchar}[1]{$^{#1}$}
\begin{thebibliography}{CDRFC19}

\bibitem[Abg17]{abgrall2017high}
R.~Abgrall.
\newblock High order schemes for hyperbolic problems using globally continuous
  approximation and avoiding mass matrices.
\newblock {\em Journal of Scientific Computing}, 73(2-3):461--494, 2017.

\bibitem[Abg18]{abgrall2018general}
R.~Abgrall.
\newblock A general framework to construct schemes satisfying additional
  conservation relations. Application to entropy conservative and entropy
  dissipative schemes.
\newblock {\em Journal of Computational Physics}, 372:640--666, 2018.

\bibitem[ABT16]{abgrall2016avoid}
R.~Abgrall, P.~Bacigaluppi, and S.~Tokareva.
\newblock How to avoid mass matrix for linear hyperbolic problems.
\newblock In {\em Numerical Mathematics and Advanced Applications ENUMATH
  2015}, pages 75--86. Springer, 2016.

\bibitem[AMO18]{abgrall2018connection}
R.~Abgrall, E.~l. Meledo, and P.~Oeffner.
\newblock On the connection between residual distribution schemes and flux
  reconstruction.
\newblock {\em arXiv preprint arXiv:1807.01261}, 2018.

\bibitem[Bre00]{bressan2000hyperbolic}
A.~Bressan.
\newblock {\em Hyperbolic systems of conservation laws: the one-dimensional
  Cauchy problem}, volume~20.
\newblock Oxford University Press on Demand, 2000.

\bibitem[CDRFC19]{chan2019efficient}
J.~Chan, D.~C. Del Rey~Fern{\'a}ndez, and M.~H. Carpenter.
\newblock Efficient entropy stable Gauss collocation methods.
\newblock {\em SIAM Journal on Scientific Computing}, 41(5):A2938--A2966, 2019.

\bibitem[CHS90]{cockburn1990runge}
B.~Cockburn, S.~Hou, and C.-W. Shu.
\newblock The {R}unge--{K}utta local projection discontinuous {G}alerkin finite
  element method for conservation laws. IV. The multidimensional case.
\newblock {\em Mathematics of Computation}, 54(190):545--581, 1990.

\bibitem[CLS89]{cockburn1989tvb2}
B.~Cockburn, S.-Y. Lin, and C.-W. Shu.
\newblock TVB {R}unge--{K}utta local projection discontinuous {G}alerkin finite
  element method for conservation laws III: one-dimensional systems.
\newblock {\em Journal of Computational Physics}, 84(1):90--113, 1989.

\bibitem[CS89]{cockburn1989tvb}
B.~Cockburn and C.-W. Shu.
\newblock TVB {R}unge--{K}utta local projection discontinuous {G}alerkin finite
  element method for conservation laws. II. General framework.
\newblock {\em Mathematics of Computation}, 52(186):411--435, 1989.

\bibitem[CS91]{cockburn1991runge}
B.~Cockburn and C.-W. Shu.
\newblock The {R}unge--{K}utta local projection
  {$P^1$}-discontinuous-{G}alerkin finite element method for scalar
  conservation laws.
\newblock {\em ESAIM: Mathematical Modelling and Numerical Analysis},
  25(3):337--361, 1991.

\bibitem[CS98]{cockburn1998runge}
B.~Cockburn and C.-W. Shu.
\newblock The {R}unge--{K}utta discontinuous {G}alerkin method for conservation
  laws V: multidimensional systems.
\newblock {\em Journal of Computational Physics}, 141(2):199--224, 1998.

\bibitem[DZLD14]{dumbser2014posteriori}
M.~Dumbser, O.~Zanotti, R.~Loub{\`e}re, and S.~Diot.
\newblock A posteriori subcell limiting of the discontinuous {G}alerkin finite
  element method for hyperbolic conservation laws.
\newblock {\em Journal of Computational Physics}, 278:47--75, 2014.

\bibitem[Gas13]{gassner2013skew}
G.~J. Gassner.
\newblock A skew-symmetric discontinuous {G}alerkin spectral element
  discretization and its relation to {SBP-SAT} finite difference methods.
\newblock {\em SIAM Journal on Scientific Computing}, 35(3):A1233--A1253, 2013.

\bibitem[Gau04]{gautschi2004orthogonal}
W.~Gautschi.
\newblock {\em Orthogonal polynomials: Computation and approximation}.
\newblock Oxford University Press on Demand, 2004.

\bibitem[Gau11]{gautschi2011numerical}
W.~Gautschi.
\newblock {\em Numerical analysis}.
\newblock Springer Science \& Business Media, 2011.

\bibitem[GG19]{glaubitz2019high}
J.~Glaubitz and A.~Gelb.
\newblock High order edge sensors with $\ell^1$ regularization for enhanced
  discontinuous {G}alerkin methods.
\newblock {\em SIAM Journal on Scientific Computing}, 41(2):A1304--A1330, 2019.

\bibitem[GKS11]{gottlieb2011strong}
S.~Gottlieb, D.~I. Ketcheson, and C.-W. Shu.
\newblock {\em Strong stability preserving {R}unge--{K}utta and multistep time
  discretizations}.
\newblock World Scientific, 2011.

\bibitem[Gla19]{glaubitz2019shock}
J.~Glaubitz.
\newblock Shock capturing by {B}ernstein polynomials for scalar conservation
  laws.
\newblock {\em Applied Mathematics and Computation}, 363:124593, 2019.

\bibitem[GNA{\etalchar{+}}19]{glaubitz2019smooth}
J.~Glaubitz, A.~Nogueira, J.~Almeida, R.~Cant{\~a}o, and C.~Silva.
\newblock Smooth and compactly supported viscous sub-cell shock capturing for
  discontinuous {G}alerkin methods.
\newblock {\em Journal of Scientific Computing}, 79(1):249--272, 2019.

\bibitem[G{\"O}S18]{glaubitz2018application}
J.~Glaubitz, P.~{\"O}ffner, and T.~Sonar.
\newblock Application of modal filtering to a spectral difference method.
\newblock {\em Mathematics of Computation}, 87(309):175--207, 2018.

\bibitem[GPR08]{gelb2008discrete}
A.~Gelb, R.~B. Platte, and W.~S. Rosenthal.
\newblock The discrete orthogonal polynomial least squares method for
  approximation and solving partial differential equations.
\newblock {\em Communications in Computational Physics}, 3(3):734--758, 2008.

\bibitem[GS98]{gottlieb1998total}
S.~Gottlieb and C.-W. Shu.
\newblock Total variation diminishing {R}unge--{K}utta schemes.
\newblock {\em Mathematics of Computation of the American Mathematical
  Society}, 67(221):73--85, 1998.

\bibitem[GST01]{gottlieb2001strong}
S.~Gottlieb, C.-W. Shu, and E.~Tadmor.
\newblock Strong stability-preserving high-order time discretization methods.
\newblock {\em SIAM review}, 43(1):89--112, 2001.

\bibitem[GVL12]{golub2012matrix}
G.~H. Golub and C.~F. Van~Loan.
\newblock {\em Matrix computations}, volume~3.
\newblock JHU Press, 2012.

\bibitem[GWK16a]{gassner2016split}
G.~J. Gassner, A.~R. Winters, and D.~A. Kopriva.
\newblock Split form nodal discontinuous {G}alerkin schemes with
  summation-by-parts property for the compressible Euler equations.
\newblock {\em Journal of Computational Physics}, 327:39--66, 2016.

\bibitem[GWK16b]{gassner2016well}
G.~J. Gassner, A.~R. Winters, and D.~A. Kopriva.
\newblock A well balanced and entropy conservative discontinuous {G}alerkin
  spectral element method for the shallow water equations.
\newblock {\em Applied Mathematics and Computation}, 272:291--308, 2016.

\bibitem[HCP12]{huerta2012simple}
A.~Huerta, E.~Casoni, and J.~Peraire.
\newblock A simple shock-capturing technique for high-order discontinuous
  {G}alerkin methods.
\newblock {\em International Journal for Numerical Methods in Fluids},
  69(10):1614--1632, 2012.

\bibitem[HK08]{hesthaven2008filtering}
J.~Hesthaven and R.~Kirby.
\newblock Filtering in {L}egendre spectral methods.
\newblock {\em Mathematics of Computation}, 77(263):1425--1452, 2008.

\bibitem[Huy09]{huybrechs2009stable}
D.~Huybrechs.
\newblock Stable high-order quadrature rules with equidistant points.
\newblock {\em Journal of Computational and Applied Mathematics},
  231(2):933--947, 2009.

\bibitem[HW07]{hesthaven2007nodal}
J.~S. Hesthaven and T.~Warburton.
\newblock {\em Nodal discontinuous {G}alerkin methods: algorithms, analysis,
  and applications}.
\newblock Springer Science \& Business Media, 2007.

\bibitem[JS94]{jiang1994cell}
G.~S. Jiang and C.-W. Shu.
\newblock On a cell entropy inequality for discontinuous {G}alerkin methods.
\newblock {\em Mathematics of Computation}, 62(206):531--538, 1994.

\bibitem[Ket08]{ketcheson2008highly}
D.~I. Ketcheson.
\newblock Highly efficient strong stability-preserving {R}unge--{K}utta methods
  with low-storage implementations.
\newblock {\em SIAM Journal on Scientific Computing}, 30(4):2113--2136, 2008.

\bibitem[KG14]{kopriva2014energy}
D.~A. Kopriva and G.~J. Gassner.
\newblock An energy stable discontinuous {G}alerkin spectral element
  discretization for variable coefficient advection problems.
\newblock {\em SIAM Journal on Scientific Computing}, 36(4):A2076--A2099, 2014.

\bibitem[KK03]{kirby2003aliasing}
R.~M. Kirby and G.~E. Karniadakis.
\newblock De-aliasing on non-uniform grids: algorithms and applications.
\newblock {\em Journal of Computational Physics}, 191(1):249--264, 2003.

\bibitem[KNG17]{kopriva2017error}
D.~A. Kopriva, J.~Nordstr{\"o}m, and G.~J. Gassner.
\newblock Error boundedness of discontinuous Galerkin spectral element
  approximations of hyperbolic problems.
\newblock {\em Journal of Scientific Computing}, 72(1):314--330, 2017.

\bibitem[Kop09]{kopriva2009implementing}
D.~A. Kopriva.
\newblock {\em Implementing spectral methods for partial differential
  equations: Algorithms for scientists and engineers}.
\newblock Springer Science \& Business Media, 2009.

\bibitem[KS74]{kreiss1974finite}
H.-O. Kreiss and G.~Scherer.
\newblock Finite element and finite difference methods for hyperbolic partial
  differential equations.
\newblock In {\em Mathematical aspects of finite elements in partial
  differential equations}, pages 195--212. Elsevier, 1974.

\bibitem[KS06]{krylov2006approximate}
V.~I. Krylov and A.~H. Stroud.
\newblock {\em Approximate calculation of integrals}.
\newblock Courier Corporation, 2006.

\bibitem[KWH11]{klockner2011viscous}
A.~Kl{\"o}ckner, T.~Warburton, and J.~S. Hesthaven.
\newblock Viscous shock capturing in a time-explicit discontinuous {G}alerkin
  method.
\newblock {\em Mathematical Modelling of Natural Phenomena}, 6(3):57--83, 2011.

\bibitem[LeV02]{leveque2002finite}
R.~J. LeVeque.
\newblock {\em Finite volume methods for hyperbolic problems}, volume~31.
\newblock Cambridge university press, 2002.

\bibitem[LT98]{levy1998semidiscrete}
D.~Levy and E.~Tadmor.
\newblock From Semidiscrete to Fully Discrete: Stability of {R}unge--{K}utta
  Schemes by The Energy Method.
\newblock {\em SIAM review}, 40(1):40--73, 1998.

\bibitem[MO16]{meister2016positivity}
A.~Meister and S.~Ortleb.
\newblock A positivity preserving and well-balanced {DG} scheme using finite
  volume subcells in almost dry regions.
\newblock {\em Applied Mathematics and Computation}, 272:259--273, 2016.

\bibitem[MOSW13]{meister2013extended}
A.~Meister, S.~Ortleb, T.~Sonar, and M.~Wirz.
\newblock An extended discontinuous {G}alerkin and spectral difference method
  with modal filtering.
\newblock {\em ZAMM-Journal of Applied Mathematics and Mechanics/Zeitschrift
  f{\"u}r Angewandte Mathematik und Mechanik}, 93(6-7):459--464, 2013.

\bibitem[NC99]{nordstrom1999boundary}
J.~Nordstr{\"o}m and M.~H. Carpenter.
\newblock Boundary and interface conditions for high-order finite-difference
  methods applied to the Euler and {N}avier--{S}tokes equations.
\newblock {\em Journal of Computational Physics}, 148(2):621--645, 1999.

\bibitem[{\"O}GR19]{offner2019stability}
P.~{\"O}ffner, J.~Glaubitz, and H.~Ranocha.
\newblock Stability of correction procedure via reconstruction with
  summation-by-parts operators for {B}urgers' equation using a polynomial chaos
  approach.
\newblock {\em ESAIM: Mathematical Modelling and Numerical Analysis (ESAIM:
  M2AN)}, 52(6):2215--2245, 02 2019.

\bibitem[Ols95a]{olsson1995summation}
P.~Olsson.
\newblock Summation by parts, projections, and stability. I.
\newblock {\em Mathematics of Computation}, 64(211):1035--1065, 1995.

\bibitem[Ols95b]{olsson1995summation2}
P.~Olsson.
\newblock Summation by parts, projections, and stability. II.
\newblock {\em Mathematics of Computation}, 64(212):1473--1493, 1995.

\bibitem[{\"O}R19]{offner2019error}
P.~{\"O}ffner and H.~Ranocha.
\newblock Error boundedness of discontinuous {G}alerkin methods with variable
  coefficients.
\newblock {\em Journal of Scientific Computing}, pages 1--36, 2019.

\bibitem[PP06]{persson2006sub}
P.-O. Persson and J.~Peraire.
\newblock Sub-cell shock capturing for discontinuous {G}alerkin methods.
\newblock In {\em 44th AIAA Aerospace Sciences Meeting and Exhibit}, page 112,
  2006.

\bibitem[PTK11]{platte2011impossibility}
R.~B. Platte, L.~N. Trefethen, and A.~B. Kuijlaars.
\newblock Impossibility of fast stable approximation of analytic functions from
  equispaced samples.
\newblock {\em SIAM review}, 53(2):308--318, 2011.

\bibitem[RG{\"O}S18]{ranocha2018stability}
H.~Ranocha, J.~Glaubitz, P.~{\"O}ffner, and T.~Sonar.
\newblock Stability of artificial dissipation and modal filtering for flux
  reconstruction schemes using summation-by-parts operators.
\newblock {\em Applied Numerical Mathematics}, 128:1--23, 2018.

\bibitem[RH73]{reed1973triangular}
W.~H. Reed and T.~Hill.
\newblock Triangular mesh methods for the neutron transport equation.
\newblock Technical report, Los Alamos Scientific Lab., N. Mex.(USA), 1973.

\bibitem[R{\"O}S16]{ranocha2016summation}
H.~Ranocha, P.~{\"O}ffner, and T.~Sonar.
\newblock Summation-by-parts operators for correction procedure via
  reconstruction.
\newblock {\em Journal of Computational Physics}, 311:299--328, 2016.

\bibitem[Run01]{runge1901empirische}
C.~Runge.
\newblock {\"U}ber empirische {F}unktionen und die {I}nterpolation zwischen
  {\"a}quidistanten {O}rdinaten.
\newblock {\em Zeitschrift f{\"u}r Mathematik und Physik}, 46(224-243):20,
  1901.

\bibitem[SM14]{sonntag2014shock}
M.~Sonntag and C.-D. Munz.
\newblock Shock capturing for discontinuous {G}alerkin methods using finite
  volume subcells.
\newblock In {\em Finite Volumes for Complex Applications VII-Elliptic,
  Parabolic and Hyperbolic Problems}, pages 945--953. Springer, 2014.

\bibitem[Str94]{strand1994summation}
B.~Strand.
\newblock Summation by parts for finite difference approximations for d/dx.
\newblock {\em Journal of Computational Physics}, 110(1):47--67, 1994.

\bibitem[TBI97]{trefethen1997numerical}
L.~N. Trefethen and D.~Bau~III.
\newblock {\em Numerical linear algebra}, volume~50.
\newblock Siam, 1997.

\bibitem[Tor13]{toro2013riemann}
E.~F. Toro.
\newblock {\em Riemann solvers and numerical methods for fluid dynamics: a
  practical introduction}.
\newblock Springer Science \& Business Media, 2013.

\bibitem[Tra09]{trangenstein2009numerical}
J.~A. Trangenstein.
\newblock {\em Numerical solution of hyperbolic partial differential
  equations}.
\newblock Cambridge University Press, 2009.

\bibitem[Van91]{vandeven1991family}
H.~Vandeven.
\newblock Family of spectral filters for discontinuous problems.
\newblock {\em Journal of Scientific Computing}, 6(2):159--192, 1991.

\bibitem[Wil70a]{wilson1970discrete}
M.~W. Wilson.
\newblock Discrete least squares and quadrature formulas.
\newblock {\em Mathematics of Computation}, 24(110):271--282, 1970.

\bibitem[Wil70b]{wilson1970necessary}
M.~W. Wilson.
\newblock Necessary and sufficient conditions for equidistant quadrature
  formula.
\newblock {\em SIAM Journal on Numerical Analysis}, 7(1):134--141, 1970.

\end{thebibliography}

\end{document}